\date{}
\documentclass[12pt]{article}

\usepackage{isomath}
\usepackage{amsmath}
\usepackage{amsbsy}
\usepackage{amssymb}
\usepackage{amscd}
\usepackage{amsfonts}
\usepackage{graphicx,color}
\usepackage{verbatim}
\usepackage{euscript}
\usepackage[parfill]{parskip}
\usepackage{alltt}
\usepackage{stmaryrd}
\usepackage[utf8]{inputenc}
\usepackage[english]{babel}
\usepackage{fancyhdr}
\usepackage{xy}
\usepackage[margin=1in]{geometry}
\usepackage{caption}
\usepackage{subcaption}
\usepackage{url}
\usepackage[titletoc,title]{appendix}
\numberwithin{equation}{section}
\usepackage{float}
\restylefloat{table}
\usepackage{setspace}
\usepackage{xcolor}
\usepackage{enumitem}

\title{Computing Singularly Perturbed Differential Equations\footnote{{\bf To appear in Journal of Computational Physics.}}}

\author{Sabyasachi Chatterjee\thanks{Dept.\ of Civil \& Environmental Engineering, Carnegie Mellon University, Pittsburgh, PA 15213. sabyasac@andrew.cmu.edu.} \and Amit Acharya\thanks{Dept.\ of Civil \& Environmental Engineering, and Center for Nonlinear Analysis, Carnegie Mellon University, Pittsburgh, PA 15213. acharyaamit@cmu.edu.} \and Zvi Artstein\thanks{Dept.\ of Mathematics, The Weizmann Institute of Science, Rehovot, Israel, 7610001. zvi.artstein@weizmann.ac.il.}}

\begin{document}

\colorlet{change}{blue}
\colorlet{ref1}{magenta}
\colorlet{ref2}{brown}

\maketitle
\begin{abstract}
\noindent A computational tool for coarse-graining nonlinear systems of ordinary differential equations in time is discussed. Three illustrative model examples are worked out that demonstrate the range of capability of the method. This includes the averaging of Hamiltonian as well as dissipative microscopic dynamics whose `slow' variables, defined in a precise sense, can often display mixed slow-fast response as in relaxation oscillations, and dependence on initial conditions of the fast variables. Also covered is the case where the quasi-static assumption in solid mechanics is violated. The computational tool is demonstrated to capture all of these behaviors in an accurate and robust manner, with significant savings in time. A practically useful strategy for accurately initializing short bursts of microscopic runs for the evolution of slow variables is integral to our scheme, without the requirement that the slow variables determine a unique invariant measure of the microscopic dynamics.
\end{abstract}
\section{Introduction}
This paper is concerned with a computational tool for understanding the behavior of systems of evolution, governed by (nonlinear) ordinary differential equations, on a time scale that is much slower than the time scales of the intrinsic dynamics. A paradigmatic example is a molecular dynamic assembly under loads, where the characteristic time of the applied loading is very much larger than the period of atomic vibrations. We examine appropriate theory for such applications and devise a computational algorithm.  The singular perturbation problems we address contain  a small parameter $\epsilon$ that reflects the ratio between the slow  and the fast time scales. In many cases, the solutions of the problem obtained by setting the small parameter to zero matches solutions to the full problem with small $\epsilon$, except in a small region - a boundary/initial layer. But, there are situations, where the  limit of solutions of the original problem as $\epsilon$ tends to zero does not match the solution of the problem obtained by setting the small parameter to zero. Our paper covers this aspect as well. In the next section we present the framework of the present study, and its sources. Before displaying our algorithm in Section \ref{sec:theory_alg}, we display previous approaches to the computational challenge. It allows us to pinpoint our contribution. Our algorithm is demonstrated through computational examples on three model problems that have been specially designed to contain the complexities in temporal dynamics expected in more realistic systems. The implementation is shown to perform robustly in all cases. These cases include the averaging of fast oscillations as well as of exponential decay, including problems where the evolution of slow variables can display fast, almost-discontinuous, behavior in time. The problem set is designed to violate any ergodic assumption, and the computational technique deals seamlessly with situations that may or may not have a unique invariant measure for averaging fast response for fixed slow variables. Thus, it is shown that initial conditions for the fast dynamics matter critically in many instances, and our methodology allows for the modeling of such phenomena. The method also deals successfully with conservative or dissipative systems. In fact, one example on which we demonstrate the efficacy of our computational tool is a linear, spring-mass, \emph{damped} system that can display permanent oscillations depending upon delicate conditions on masses and spring stiffnesses and initial conditions; we show that our methodology does not require \emph{a-priori} knowledge of such subtleties in producing the correct response.

\section{The framework}
A particular case of the differential equations we deal with is of the form
\begin{align}\label{eq:1.1}
\frac{dx}{dt} &= \frac{1}{\epsilon} F(x)+ G(x),
\end{align}
with $ \epsilon > 0 $ a small real parameter, and $x \in \mathbb{R}^{n}$. For reasons that will become clear in the sequel we refer to the component $G(x)$ as the drift component.
\par
Notice that the dynamics in (\ref{eq:1.1}) does not exhibit a prescribed split into a fast and a slow dynamics. We are
interested in the case where such a split is either not tractable or does not exist.
\par
Another particular case where a split into a fast and slow dynamics can be identified, is also of interest to us, as follows.
\begin{align}\label{eq:1.2}
\frac{dx}{dt} &= \frac{1}{\epsilon} F(x,l)\\
\frac{dl}{dt} &= L(x,l),\nonumber
\end{align}
\par
with $x \in \mathbb{R}^n$ and $l \in \mathbb{R}^m$. We think of the variable $l$ as a load. Notice that the dynamics of the load is determined by an external ``slow" equation, that, in turn, may be affected by the ``fast” variable $x$.
\par
The general case we study is a combination of the previous two cases, namely,
\begin{align}\label{eq:1.3}
\frac{dx}{dt} &= \frac{1}{\epsilon} F(x,l)+ G(x,l)\\
\frac{dl}{dt} &= L(x,l),\nonumber
\end{align}
\par
which accommodates both a drift and a load. In the theoretical discussion we address
the general case. We display the two particular cases, since there are many interesting
examples of the type (\ref{eq:1.1}) or (\ref{eq:1.2}).
\par
An even more general setting would be the case where the right hand side of (\ref{eq:1.3})
is of the form $H(x, l, \epsilon)$, namely, there is no a priori split of the right hand side of the
equation into fast component and a drift or a slow component. A challenge then would
be to identify, either analytically or numerically, such a split. We do not address this case
here, but our study reveals what could be promising directions of such a general study.
\par
We recall that the parameter $\epsilon$ in the previous equation represents the ratio between the slow (or ordinary) part in the equation and the fast one. In Appendix B we examine one of our examples, and demonstrate how to derive the dimensionless equation with the small parameter, from the raw mechanical equation. In real world situations, $\epsilon$ is small yet it is not infinitesimal. Experience teaches us, however, that the limit behavior, as $\epsilon$ tends to 0, of the solutions is quite helpful in understanding of the physical phenomenon and in the computations. This is, indeed, demonstrated in the examples that follow.
\par
References that carry out a study of equations of the form (\ref{eq:1.1}) are, for instance, Tao,
Owhadi and Marsden \cite{Tao_2010_Noninstructive}, Artstein, Kevrekidis, Slemrod and Titi \cite{artstein2007slow}, Ariel, Engquist and Tsai \cite{ariel2009multiscale, ariel2009numerical},  Artstein, Gear, Kevrekidis, Slemrod and Titi \cite{Artstein_Gear_2011_KdV_burgers}, Slemrod and Acharya \cite{Slemord_Acharya_2012_Time_average}; conceptually similar questions implicitly arise in the work of Kevrekidis et al. \cite{kevrekidis2003equation}. The form (\ref{eq:1.2}) coincides with the Tikhonov model, see, e.g., O'Malley \cite{O_Malley_2014_Singular_perturb}, Tikhonov, Vasileva and Sveshnikov \cite{Tikhonov_1985_DE}, Verhulst \cite{Verhulst_2005_App_of_sing_perturb}, or Wasow \cite{Wasow_1965_Asym}. The literature concerning this case followed, mainly, the so called Tikhonov approach, namely, the assumption that the solutions of the
$x$-equation in (\ref{eq:1.2}), for $l$ fixed, converge to a point $x(l)$ that solves an algebraic equation,
namely, the second equation in (\ref{eq:1.3}) where the left hand side is equal to 0. The limit dynamics then is a trajectory $(x(t), l(t))$, evolving on the manifold of stationary points $x(l)$. We are interested, however, in the case where the limit dynamics may not be determined by such a manifold, and may exhibit infinitely rapid oscillations. A theory and applications alluding to such a case are available, e.g., in Artstein and Vigodner \cite{Artstein_Vigonder_1996_dynamic_limit}, Artstein \cite{Artstein_2002_perturbed}, Acharya \cite{Acharya2007coarse_variable, Acharya_2010_Coarsegraining_autonomous},
Artstein, Linshiz and Titi \cite{Artstein_Linshiz_2007_Young_measure}, Artstein and Slemrod \cite{Artstein_Slemord_2001_sigular_perturb_limit}.
\par
\section{The goal}
A goal of our study is to suggest efficient computational tools that help revealing the limit
behavior of the system as $\epsilon$ gets very small, this on a prescribed, possibly long, interval.
The challenge in such computations stems from the fact that, for small $\epsilon$, computing the
ordinary differential equation takes a lot of computing time, to the extent that it becomes
not practical. Typically, we are interested in a numerical description of the full solution,
namely, the progress of the coupled slow/fast dynamics. At times, we may be satisfied with
partial information, say in the description of the progress of a slow variable, reflecting a
measurement of the underlying dynamics. To that end we first identify the mathematical
structure of the limit dynamics on the given interval. The computational algorithm will
reveal an approximation of this limit dynamics, that, in turn, is an approximation of the
full solution for arbitrarily small $\epsilon$. If only a slow variable is of interest, it can be derived
from the established approximation.
\section{The limit dynamics}
In order to achieve the aforementioned goal, we display the limit structure, as $\epsilon \rightarrow 0,$ of the dynamics of (\ref{eq:1.3}). To this end we identify the $\textit{fast time equation}$
\begin{align}\label{eq:3.1}
\frac{dx}{d\sigma} = F(x, l),
\end{align}
when $l$ is held fixed (recall that $l$ may not show up at all, as in (\ref{eq:1.1})). The equation (\ref{eq:3.1})
is the $\textit{fast part}$ of (\ref{eq:1.3}) (as mentioned, $G(x)$ is the $\textit{drift}$ and the solution $l(t)$ of the load
equation is the $\textit{load}$).
\par
Notice that when moving from (\ref{eq:1.3}) to (\ref{eq:3.1}), we have changed the time scale, with
$t = \epsilon \sigma.$ We refer to $\sigma$ as the fast time scale.
\par
In order to describe the limit dynamics of (\ref{eq:1.3}) we need the notions of: Probability
measures and convergence of probability measures, Young measures and convergence in the Young measures sense, invariant measures and limit occupational measures. In particular, we shall make frequent use of the fact that when occupational measures of solutions of (\ref{eq:3.1}), on long time intervals, converge, the limit is an invariant measure of (\ref{eq:3.1}). A concise explanation of these notions can be found, e.g., in \cite{Artstein_Vigonder_1996_dynamic_limit, artstein2007slow}.
\par
It was proved in \cite{Artstein_Vigonder_1996_dynamic_limit} for (\ref{eq:1.2}) and in \cite{artstein2007slow} for (\ref{eq:1.1}), that under quite general conditions, the dynamics converge, as $\epsilon\rightarrow0,$ to a Young measure, namely, a probability measure-valued map, whose values are invariant measures of (\ref{eq:3.1}). These measures are drifted in the case of (\ref{eq:1.1}) by the drift component of the equation, and in the case (\ref{eq:1.2}) by the load. We display the result in the general case after stating the assumptions under which the result holds.
\par
\textbf{Assumption 4.1.} The functions $F(.,.), G(.,.)$ and $L(.,.)$ are continuous. The solutions, say $x(.)$, of the fast equation (\ref{eq:3.1}), are determined uniquely by the initial data, say $x(\sigma_{0}) = x_{0}$, and stay bounded for $\sigma \geq \sigma_{0}$, uniformly for $x_{0}$ and for $l$ in bounded sets.
\par
Here is the aforementioned result concerning the structure of the limit dynamics.
\par
\textbf{Theorem 4.2.} For every sequence $\epsilon_{i}\rightarrow0$ and solutions $(x_{\epsilon_{i}}(t), l_{\epsilon_{i}}(t))$ of the perturbed equation (\ref{eq:1.3}) defined on $[0, T]$, with $(x_{\epsilon_{i}}(0), l_{\epsilon_{i}}(0))$ in a bounded set, there exists a subsequence $\epsilon_{j}$ such that $(x_{\epsilon_{j}}(.), l_{\epsilon_{j}}(.))$ converges as $j\rightarrow\infty$, where the convergence in the $x$-coordinates is in the sense of Young measures, to a Young measure, say $\mu(.)$, whose values are invariant measures of the fast equation (\ref{eq:3.1}), and the convergence in the $l$-coordinates is uniform on the interval, with a limit, say $l_{0}(.)$, that solves the differential equation
\begin{align}\label{eq:3.2}
\frac{dl}{dt}=\int_{\mathbb{R}^{n}}L(x,l)\mu(t)dx.
\end{align}
\par
The previous general result has not been displayed in the literature, but the arguments in \cite{artstein2007slow} in regard to (\ref{eq:1.1}) or the proof given in \cite{Artstein_Vigonder_1996_dynamic_limit} for the case (\ref{eq:1.2}), apply to the present setting as well.
\section{Measurements and slow observables}
A prime role in our approach is played by $\textit{slow observables}$, whose dynamics can be followed. The intuition behind the notion is that the observations which the observable reveals, is a physical quantity on the macroscopic level, that can be detected. Here we identify some candidates for such variables. The role they play in the computations is described in the next section.
\par
In most generality, an $\textit{observable}$ is a mapping that assigns to a probability measure $\mu(t)$ arising as a value of the Young measure in the limit dynamics of (\ref{eq:1.3}), a real number, or a vector, say in $\mathbb{R}^{k}$. Since the values of the Young measure are obtained as limits of occupational measures (that in fact we use in the computations), we also demand that the observable be defined on these occupational measures, and be continuous when passing from the occupational measures to the value of the Young measure.
\par
An observable $v(.)$ is a $\textit{slow observable}$ if when applied to the Young measure $\mu(.)$ that determines the limit dynamics in Theorems 4.2, the resulting vector valued map $v(t) = v(\mu(t), l(t))$ is continuous at points where the measure $\mu(.)$ is continuous.
\par
An $\textit{extrapolation rule}$ for a slow observable $v(.)$ determines an approximation of the value $v(t + h)$, based on the value $v(t)$ and, possibly, information about the value of the Young measure $\mu(t)$ and the load $l(t)$, at the time $t$. A typical extrapolation rule would be generated by the derivative, if available, of the slow observable. Then $v(t + h) = v(t) + h\frac{dv}{dt}(t).$
\par
A trivial example of a slow observable of (\ref{eq:1.3}) with an extrapolation rule is the variable $l(t)$ itself. It is clearly slow, and the right hand side of the differential equation $\eqref{eq:3.2}$ determines the extrapolation rule, namely : 
\begin{align}\label{eq:4.1}
l(t + h) = l(t) + h\frac{dl}{dt}(t).
\end{align}
\par
An example of a slow observable possessing an extrapolation rule in the case of (\ref{eq:1.1}), is
an $\textit{orthogonal observable}$, introduced in \cite{artstein2007slow}. It is based on a mapping $m(x, l):\mathbb{R}^{n}\rightarrow \mathbb{R}$
which is a first integral of the fast equation (\ref{eq:3.1}) (with $l$ fixed), namely, it is constant along solutions of ($\ref{eq:3.1}$). Then we define the observable $v(\mu) = m(x,l)$ with $x$ any point in the support of $\mu$. But in fact, it will be enough to assume that the mapping $m(x,l)$ is constant on the supports of the invariant measures arising as values of a Young measure. The definition of $v(\mu) = m(x, l)$ with $x$ any point in the support of $\mu$ stays the same, that is, $m(x,l)$ may not stay constant on solutions away from the support of the limit invariant measure. It was shown in \cite{artstein2007slow} for the case (\ref{eq:1.1}), that if $m(.)$ is continuously differentiable, then $v(t)$
satisfies, almost everywhere, the differential equation
\begin{align}\label{eq:4.2}
\frac{dv}{dt}= \int_{\mathbb{R}^{n}}\nabla m(x)G(x)\mu(t)dx.
\end{align}
\par
It is possible to verify that the result holds also when the observable satisfies the weaker condition just described, namely, it is a first integral only on the invariant measures that arise as values of the limit Young measure. The differential equation ($\ref{eq:4.2}$) is not in a closed form, in particular, it is not an ordinary differential equation. Yet, if one knows $\mu(t)$ and $v(t)$ at time $t$, the differentiability expressed in ($\ref{eq:4.2}$) can be employed to get an extrapolation of the form $v(t + h) = v(t) + h \frac{dv}{dt}(t)$ at points of continuity of the Young measure, based on the right hand side of (\ref{eq:4.2}). A drawback of an  orthogonal observable for practical purposes is that finding first integrals of the fast motion is, in general, a non-trivial matter.
\par
A natural generalization of the orthogonal observable would be to consider a moment
or a generalized moment, of the measure $\mu(t)$. Namely, to drop the orthogonality from the definition, allowing a general $m$ : $\mathbb{R}^{n}\rightarrow \mathbb{R}$ be a measurement (that may depend, continuously though, on $l$ when $l$ is present), and define
\begin{align}\label{eq:4.3}
v(\mu) = \int_{\mathbb{R}^n} m(x) \mu(dx).
\end{align}
Thus, the observable is an average, with respect to the probability measure, of the bounded continuous measurement $m(.)$ of the state. If one can verify, for a specific problem, that $\mu(t)$ is piecewise continuous, then the observable defined in (\ref{eq:4.3}) is indeed slow. The drawback of such an observable is the lack of an apparent extrapolation rule. If, however, in a given application, an extrapolation rule for the moment can be identified, it will become a useful tool in the analysis of the equation.
\par
A generalization of (\ref{eq:4.3}) was suggested in \cite{acharya2006computational,ariel2009multiscale} in the form of running time-averages as slow variables, and was made rigorous in the context of
delay equations in \cite{Slemord_Acharya_2012_Time_average}. Rather than considering the average of the bounded and continuous function $m(x)$ with respect $\mu(t)$, we suggest considering the average with respect to the values of the Young measure over an interval [$t-\Delta,t$], i.e,
\begin{align}\label{eq:4.4}
v(t)=\frac{1}{\Delta}\int_{t-\Delta}^{t} \int_{\mathbb{R}^n} m(x)\mu(s)(dx)ds.
\end{align}
\par
Again, the measurement $m$ may depend on the load. Now the observable (\ref{eq:4.4}) depends not only on the value of the measure at $t$, but on the ``history" of the Young measure, namely its values on [$t-\Delta,t$]. The upside of the definition is that $v(t)$ is a Lipschitz function of $t$ (the Lipschitz constant may be large when $\Delta$ is small) and, in particular, is almost everywhere differentiable. The almost everywhere derivative of the slow variable is expressed at the points $t$ where $\mu(.)$ is continuous at $t$ and at $t-\Delta$, by
\begin{align}\label{eq:4.5}
\frac{dv}{dt}=\frac{1}{\Delta} \left( \int_{\mathbb{R}^n} m(x)\mu(t)(dx)-\int_{\mathbb{R}^n} m(x)\mu(t-\Delta)(dx)\right).
\end{align}
This derivative induces an extrapolation rule.
\par
For further reference we call an observable that depends on the values of the Young measure over an interval prior to $t$, an $\textit{H-observable}$ (where the $H$ stands for history).
\par
An $H$-observable need not be an integral of generalized moments, i.e., of integrals. For instance, for a given measure $\mu$ let
\begin{align}\label{eq:4.6}
r(\mu) = max\{x \cdot e_{1}: x \in supp(\mu)\},
\end{align}
where $e_{1}$ is a prescribed unit vector and supp($\mu$) is the support of $\mu$. Then, when supp($\mu$) is continuous in $\mu$, (and recall Assumption 4.1)  the expression
\begin{align}\label{eq:4.7}
v(t) = \frac{1}{\Delta}\int_{t-\Delta}^{t} r(\mu(\tau))d\tau,
\end{align}
is a slow observable, and
\begin{align}\label{eq:4.8}
\frac{dv}{dt} = \frac{1}{\Delta} \left(r(\mu(t))- r(\mu(t-\Delta))\right)
\end{align}
determines its extrapolation rule.
\par
The strategy we display in the next section applies whenever slow observables with valid extrapolation rules are available. The advantage of the $H$-observables as slow variables is that any smooth function $m(.)$ generates a slow observable and an extrapolation rule. Plenty of slow variables arise also in the case of generalized moments of the measure, but then it may be difficult to identify extrapolation rules. The reverse situation occurs with orthogonal observables. It may be difficult to identify first integrals of (\ref{eq:3.1}), but once such an integral is available, its extrapolation rule is at hand.
\par
Also note that in all the preceding examples the extrapolation rules are based on derivatives. We do not exclude, however, cases where the extrapolation is based on a different argument. For instance, on information of the progress of some given external parameter, for instance, a control variable. All the examples computed in the present paper will use $H$-observables.
\section{The algorithm}\label{sec:theory_alg}
Our strategy is a modification of a method that has been suggested in the literature and applied in some specific cases. We first describe these, as it will allow us to pinpoint our contribution.
\par
A computational approach to the system (\ref{eq:1.2}) has been suggested in Vanden-Eijnden \cite{Vanden_2003_Numtech} and applied in Fatkullin and Vanden-Eijnden \cite{Fatkullin_2004_Lorenz}. It applies in the special case where the fast process is stochastic, or chaotic, with a unique underlying measure that may depend on the load. The underlying measure is then the invariant measure arising in
the limit dynamics in our Theorem 4.2. It can be computed by solving the fast equation, initializing it at an arbitrary point. The method carried out in \cite{Vanden_2003_Numtech} and \cite{Fatkullin_2004_Lorenz} is, roughly, as follows. Suppose the value of the slow variable (the load in our terminology) at time $t$ is known. The fast equation is run then until the invariant measure is revealed. The measure is then employed in the averaging that determines the right hand side of the slow equation at $t$, allowing to get a good approximation of the load variable at $t + h$. Repeating this scheme results in a good approximation of the limit dynamics of the system. The method relies on the property that the value of the load determines the invariant measure. The
latter assumption has been lifted in \cite{Artstein_Linshiz_2007_Young_measure}, analyzing an example where the dynamics is not
ergodic and the invariant measure for a given load is not unique, yet the invariant measure
appearing in the limit dynamics can be detected by a good choice of an initial condition
for the fast dynamics. The method is, again, to alternate between the computation of
the invariant measure at a given time, say $t$, and using it then in the averaging needed to determine the slow equation. Then determine the value of the load at $t+h$. The structure
of the equation allows to determine a good initial point for computing the invariant measure
at $t + h$, and so on and so forth, until the full dynamics is approximated.
\par
The weakness of the method described in the previous paragraph is that it does not
apply when the split to fast dynamics and slow dynamics, i.e. the load, is not available,
and even when a load is there, it may not be possible to determine the invariant measure
using the value of the load.
\par
Orthogonal observables were employed in \cite{artstein2007slow} in order to analyze the system (\ref{eq:1.1}), and a computational scheme utilizing these observables was suggested. The scheme was applied in \cite{Artstein_Gear_2011_KdV_burgers}. The suggestion in these two papers was to identify orthogonal observables whose values determine the invariant measure, or at least a good approximation of it. Once such observables are given, the algorithm is as follows. Given an invariant
measure at $t$, the values of the observables at $t + h$ can be determined based on their
values at $t$ and the extrapolation rule based on (\ref{eq:4.2}). Then a point in $\mathbb{R}^n$ should be found, which
is compatible with the measurements that define the observables. This point would be
detected as a solution of an algebraic equation. Initiating the fast equation at that point
and solving it on a long fast time interval, reveals the invariant measure. Repeating the
process would result in a good approximation of the full dynamics.
\par
The drawback of the previous scheme is the need to find orthogonal observables,
namely measurements that are constant along trajectories within the invariant measures
in the limit of the fast flow, and verify that their values determine the value of the Young
measure, or at least a good approximation of it. Also note that the possibility to determine
the initialization point with the use of measurements, rather than the observables, relies
on the orthogonality. Without orthogonality, applying the measurements to the point of
initialization of the fast dynamics, yields no indication. In this connection, it is important to mention the work of Ariel, Engquist, and Tsai \cite{ariel2009multiscale,ariel2009numerical} that demonstrates theory, and computations utilizing the Hetergeneous Multiscale Modeling (HMM) scheme of E and Engquist \cite{weinan2003heterognous}, for defining a complete set of slow variables that determine the unique invariant measure for microscopic systems equipped with such.

The scheme that we suggest employs general observables with extrapolation rules. As
mentioned, there are plenty of these observables, in particular $H$-observables, as 
(\ref{eq:4.4}) indicates. The scheme shows how to use them in order to compute the full dynamics,
or a good approximation of it.
\par
It should be emphasized that none of our examples satisfy the ergodicity assumption placed in the aforementioned literature. Also, in none of the examples it is apparent, if possible at all, to find orthogonal observables. In addition, our third example exhibits extremely rapid variation in the slow variable (resembling a jump), a phenomenon not treated in the relevant literature so far.
\par
We provide two versions of the scheme. One for observables determined by the values
$\mu(t)$ of the Young measure, and the second for $H$-observables, namely, observables depending
on the values of the Young measure over an interval [$t-\Delta,t$]. The modifications needed
in the latter case are given in parentheses.
\par
\textbf{The scheme}. Consider the full system (\ref{eq:1.3}), with initial conditions $x(t_{0}) = x_{0}$ and
$l(t_{0}) = l_{0}$. Our goal is to produce a good approximation of the limit solution, namely the
limit Young measure $\mu(.)$, on a prescribed interval [$t_{0}, T_{0}$].
\par
\textbf{A general assumption}. The function of time defined by the closest-point projection, in some appropriate metric, of any fixed point in state-space on the support of $\mu(t)$ for each $t$ in any interval in which $\mu(\cdot)$ is continuous, is smooth. A number of slow observables, say $v_1,\ldots, v_k$ can be identified, each of them equipped with an extrapolation rule, valid at all continuity points of the Young measure.
\par
\textbf{Initialization of the algorithm}. Solve the fast equation (\ref{eq:3.1}) with initial condition
$x_{0}$ and a fixed initial load $l_{0}$, long enough to obtain a good approximation of the value
of the Young measure at $t_{0}$. (In the case of an $H$-observable solve the full equation with
$\epsilon$ small, on an interval [$t_{0}, t_{0} + \Delta$], with $\epsilon$ small enough to get a good approximation
of the Young measure on the interval). In particular the initialization produces good
approximations of $\mu(t)$ for $t = t_{0}$ (for $t = t_{0} +\Delta$ in case of an $H$-observable). Compute the
values $v_{1}(\mu(t)), ..., v_{k}(\mu(t))$ of the observables at this time $t$.
\par
\textbf{The recursive part.}
\par
\textbf{Step 1:} We assume that at time $t$ the values of the observables $v_{1}(\mu(t)), ..., v_{k}(\mu(t))$,
applied to the value $\mu(t)$ of the Young measure, are known. We also assume that we have enough information to invoke the extrapolation rule to these observables (for instance, we can compute $\frac{dv}{dt}$ which determines the extrapolation when done via a derivative). If the model has a load variable, it should be one of the observables.
\par
\textbf{Step 2:} Apply the extrapolation rule to the observables and get an approximation of the values $v_{1}(\mu(t + h)),..., v_{k}(\mu(t + h))$, of the observables at time $t + h$ (time $t + h -\Delta$ in the case of an $H$-observable). Denote the resulting approximation by ($v_{1}$, ..., $v_{k}$).
\par
\textbf{Step 3:} Make an intelligent guess of a point $x(t + h)$ (or $x(t + h -\Delta)$ in the case of
$H$-observable), that is in the basin of attraction of $\mu(t + h)$ (or $\mu(t + h -\Delta)$ in the case of $H$-observable). See a remark below concerning the intelligent guess.
\par
\textbf{Step 4:} Initiate the fast equation at $(x(t + h), l(t + h))$ and run the equation until a good approximation of an invariant measure arises (initiate the full equation, with small $\epsilon$
at $(x(t + h -\Delta), l(t + h -\Delta))$, and run it on [$t + h -\Delta, t + h$], with $\epsilon$ small enough such that a good approximation for a Young measure on the interval is achieved).
\par
\textbf{Step 5:} Check if the invariant measure $\mu(t + h)$ revealed in the previous step is far
from $\mu(t)$ (this step, and consequently step 6.1, should be skipped if the Young measure
is guaranteed to be continuous).
\par
\textbf{Step 6.1:} If the answer to the previous step is positive, it indicates that a point
of discontinuity of the Young measure may exist between $t$ and $t + h$. Go back to $t$ and
compute the full equation (initiating it with any point on the support of $\mu(t)$) on [$t, t+h$],
or until the discontinuity is revealed. Then start the process again at Step 1, at the time
$t + h$ (or at a time after the discontinuity has been revealed).
\par
\textbf{Step 6.2:} If the answer to the previous step is negative, compute the values of the
observables $v_{1}(\mu(t + h)), ..., v_{k}(\mu(t + h))$ at the invariant measure that arises by running
the equation. If there is a match, or almost a match, with $(v_{1}, ..., v_{k})$, accept the value $\mu(t + h)$ (accept the computed values on [$t + h -\Delta,t + h$] in the case of an $H$-observable) as the value of the desired Young measure, and start again at Step 1, now at time $t + h$.
If the match is not satisfactory, go back to step 3 and make an improved guess.
\par
{\setlength{\parindent}{0cm}
\textbf{Conclusion of the algorithm.} Continue with steps 1 to 6 until an approximation of the
Young measure is computed on the entire time interval [$t_{0}, T_{0}$].
}
\par
\textbf{Making the intelligent guess in Step 3.} The goal in this step is to identify a point in the basin of attraction of $\mu(t+h)$. Once a candidate for such a point is suggested, the decision whether to accept it or not is based on comparing the observables computed on the invariant measure generated by running the equation with this initial point, to the values as predicted by the extrapolation rule. If there is no match, we may improve the initial suggestion for the point we seek.
\par
In order to get an initial point, we need to guess the direction in which the invariant measure is drifted. We may assume that the deviation of $\mu(t+h)$ from $\mu(t)$ is similar to the deviation of $\mu(t)$ from $\mu(t-h)$.
Then the first intelligent guess would be, say, a point $x_{t+h}$ such that
$x_{t+h}-x_{t}= x_{t}-x_{t-h}$ where $x_{t}$ is a point in the support of $\mu(t)$ and $x_{t-h}$ is the point in the support of $\mu(t-h)$ closest, or near closest, to $x_{t}$. At this point, in fact, we may try several such candidates, and check them in parallel. If none fits the criterion in Step 6.2, the process that starts again at Step 3, could use the results in the previous round, say by perturbing the point that had the best fit in a direction of, hopefully, a better fit. A sequence of better and better approximations can be carried out until the desired result is achieved.
\section{The expected savings in computer time}
\par
The motivation behind our scheme of computations is that the straightforward approach, namely, running the entire equation (\ref{eq:1.3}) on the full interval, is not feasible if $\epsilon$ is very small. To run (\ref{eq:3.1}) in order to compute the invariant measure at a single point $t$, or to run (\ref{eq:1.3}) on a short interval [$t -\Delta,t$], with $\Delta$ small, does not consume a lot of computing time. Thus, our scheme replaces the massive computations with computing the values of the Young measure at a discrete number of points, or short intervals, and the computation of the extrapolation rules to get an estimate of the progress of the observables. The latter step does not depend on $\epsilon$, and should not consume much computing time. Thus, if $h$ is large (and large relative to $\Delta$ in the case of $H$-observables), we achieve a considerable saving.
\par
These arguments are also behind the saving in the cited references, i.e., \cite{Vanden_2003_Numtech, Fatkullin_2004_Lorenz, Artstein_Linshiz_2007_Young_measure, Artstein_Gear_2011_KdV_burgers}. In our algorithm there is an extra cost of computing time, namely, the need to detect points in the basin of attraction of the respective invariant
measures, i.e., Step 3 in our algorithm. The extra steps amount to, possibly, an addition
of a discrete number of computations that reveal the invariant measures.
An additional computing time may be accrued when facing a discontinuity in the
Young measure. The cost is, again, a computation of the full Young measure around the
discontinuity. The possibility of discontinuity has not been address in the cited references.
\section{Error estimates}
\par
The considerations displayed and commented on previously were based on the heuristics
behind the computations. Under some strong, yet common, conditions on the smoothness
of the processes, one can come up with error estimate for the algorithm. We now
produce such an estimate, which is quite standard (see e.g., \cite{Num1978}). Of interest here is the nature of assumptions on our process, needed to guarantee
the estimate.
\par
Suppose the computations are performed on a time interval [$0, T$], with time steps
of length $h$. The estimate we seek is the distance between the computed values, say
$P(kh)$, $k = 0, 1, ...,N,$ of invariant measures, and the true limit dynamics, that is the
values $\mu(kh)$ of the limit Young measure, at the same mesh points (here $hN$ is close to $T$).
\par
Recall that the Young measure we compute is actually the limit as $\epsilon\rightarrow0$ of solutions of \eqref{eq:1.3}. Still, we think of $\mu(t)$ as reflecting a limiting dynamics, and consider the value $\mu(0)$ as its initial state, in the space of invariant measures (that may indeed be partially supported on points, for instance, in case the load is affected by the dynamics, we may wish to compute it as well). We shall denote by $\mu(h,\nu)$ the value of the Young measure obtained at the time $h$, had $\nu$ been the initial state. Likewise, we denote by $P(h,\nu)$ the result of the computations at time $h$, had $\nu$ been the invariant measure to which the algorithm is applied. We place the assumptions on $\mu(h,\nu)$ and $P(h,\nu)$. After stating the assumption we comment on the reflection of them on the limit dynamics and the algorithm. Denote by $\rho(.,.)$ a distance, say the Prohorov metric, between probability measures.
\par
\textbf{Assumption 8.1.} There exists a function $\delta(h)$, continuous at $0$ with $\delta(0) = 0$, and a constant $\eta$ such that
\begin{align}\label{eq:8.1}
\rho(P(h,\nu), \mu(h,\nu))\leq\delta(h)h,
\end{align}
and
\begin{align}\label{eq:8.2}
\rho(\mu(h,\nu_{1}),\mu(h,\nu_{2}))\leq (1+\eta h)\rho(\nu_{1},\nu_{2}).
\end{align}
\par
\textbf{Remark.} Both assumptions relate to some regularity of the limiting dynamics. The
inequality (\ref{eq:8.2}) reflects a Lipschitz property. It holds, for instance, in case $\mu(t)$ is a solution of an ordinary differential equations with Lipschitz right hand side. Recall that without
a Lipschitz type condition it is not possible to get reasonable error estimates even in the
ode framework. Condition (\ref{eq:8.1}) reflects the accuracy of our algorithm on a time step of
length $h$. For instance, If a finite set of observables determines the invariant measure, if
the extrapolation through derivatives are uniform, if the mapping that maps the measures
to the values of the observables is bi-Lipschitz, and if the dynamics $\mu(t)$ is Lipschitz, then
(\ref{eq:8.1}) holds. In concrete examples it may be possible to check these properties directly (we
comment on that in the examples below).
\par
{\bf Theorem 8.2.} Suppose the inequalities in Assumption 8.1 hold. Then
\begin{align}\label{eq:8.3}
\rho(P(kh),\mu(kh)) \leq \delta(h) \eta^{-1}(e^{\eta{T}}-1).
\end{align}
{\bf Proof.} Denote $E_{k} = \rho((P(kh),\mu(kh))$, namely, $E_{k}$ is the error accrued from $0$ to the time $kh$. Then $E_{0}= 0$. Since $P(kh) = P(h,P((k-1)h))$ and $\mu(kh) = \mu(h,\mu((k-1)h))$ it follows that
\begin{align}\label{eq:8.4}
 E_{k} \leq \rho(P(h,P((k-1)h)),\mu(h,P((k-1)h))) +\nonumber\\
  \rho(\mu(h,P((k-1)h)),\mu(h,\mu((k-1)))).
\end{align}
From inequalities (\ref{eq:8.1}) and (\ref{eq:8.2}) we get
\begin{align}\label{eq:8.5}
E_{k} \leq \delta(h)h + (1 + \eta h)E_{k-1}.
\end{align}
Spelling out the recursion we get
\begin{align}\label{eq:8.6}
E_{k} & \leq \delta(h)h(1+ (1+\eta h) + (1+\eta h)^{2} + ... + (1 + \eta h)^{k-1})\nonumber\\
          &= \delta(h) \eta^{-1} ((1+ \eta h)^{k}-1).
\end{align}
Since $hk \leq T$ the claimed inequality (\ref{eq:8.3}) follows.
\par
\textbf{Remark.} The previous estimate imply that the error tends to zero as the step size
$h$ tends to zero. Note, however, that the bigger the step size, the bigger the saving of
computational time is. This interplay is common in computing in general.
\par
\textbf{Remark.} Needless to say, the previous estimate refers to the possible errors on
the theoretical dynamics of measures. A different type of errors that occur relates to the
numerics of representing and approximating the measures. The same issue arises in any
numerical analysis estimates, but they are more prominent in our framework due to the
nature of the dynamics.

\section{Computational Implementation}\label{impl_algo}
We are given the initial conditions of the fine and the slow variables, $x(-\Delta) = x_0 $ and $l(-\Delta) = l_0 $. We aim to produce a good approximation, of the limit solution in the period $[0, T_0]$. Due to the lack of algorithmic specification for determining orthogonal observables, we concentrate on $H$-observables in this paper.

In the implementation that follows, we get to the specifics of how the calculations are carried out. Recall that the $H$-observables are averages of the form
\begin{align}\label{coarse_obs_impl}
v(t)=\frac{1}{\Delta}\int_{t-\Delta}^{t} \int_{\mathbb{R}^n} m(x)\mu(s)(dx)ds.
\end{align}
Namely, the $H$-observables are slow variables with a time-lag, i.e. the fine/microscopic system needs to have run on an interval of length $\Delta$ before the observable at time $t$ can be defined.

\emph{From here onwards, with some abuse of notation, whenever we refer to only a measure at some instant of time we mean the value of the Young measure of the fast dynamics at that time. When we want to refer to the Young measure, we mention it explicitly.}

We will also refer to any component of the list of variables of the original dynamics \eqref{eq:1.3} as \emph{fine} variables.

We think of the calculations marching forward in the slow time-scale in discrete steps of size $h$ with $T_0 = nh$.  Thus the variable $t$ below in the description of our algorithm takes values of $0h,1h,2h,\ldots,nh$.

\par 
\textbf{Step 1: Calculate the rate of change of slow variable} \\ We calculate the rate of change of the slow variable at time $t$ using the following: 
\begin{align}\label{eq:comp_impl_1}
\frac{dv}{dt}(t)=\frac{1}{\Delta} \left( \int_{\mathbb{R}^n} m(x)\mu(t)(dx)-\int_{\mathbb{R}^n} m(x)\mu(t-\Delta)(dx)\right).
\end{align}

Let us denote the term $\int_{\mathbb{R}^n} m(x)\mu(t)(dx)$ as $R_t^m$ and the term $\int_{\mathbb{R}^n} m(x)\mu(t-\Delta)(dx)$ as $R_{t-\Delta}^m$. The term $R_t^m$ is computed as
\begin{align}\label{eq:comp_impl_R1}
R_t^m = \frac {1}{N_t} \sum_{i=1}^{N_t} m(x_\epsilon(\sigma_i),l_\epsilon(\sigma_i)).
\end{align}

The successive values ($x_\epsilon(\sigma_i)$, $l_\epsilon(\sigma_i)$) are obtained by running the fine system
\begin{equation}\label{eq:alg_fast_2}
\begin{split}
\frac{dx_\epsilon}{d\sigma} &= F(x_\epsilon,l_\epsilon) + \epsilon G(x_\epsilon,l_\epsilon)\\
\frac{dl_\epsilon}{d\sigma}& = \epsilon L(x_\epsilon,l_\epsilon),
\end{split}
\end{equation}
 with initial condition $x_{guess}(\sigma = \frac{t}{\epsilon})$ and $l(\sigma=\frac{t}{\epsilon} )$. 
\par 
We discuss in Step 5 how we obtain $x_{guess}(\sigma)$. Here, $N_t$ is the number of increments taken for the value of $R_t^m$ to converge upto a specified value of tolerance. Also, $N_t$ is large enough such that the effect of the initial transient does not affect the value of $R_t^m$. 
\par 
Similarly, $R_{t-\Delta}^m$ is computed as: 

\begin{align}\label{eq:comp_impl_R2}
R_{t-\Delta}^m = \frac {1}{N_{t-\Delta}} \sum_{i=1}^{N_{t-\Delta}} m\left(x_\epsilon(\sigma_i), l_\epsilon(\sigma_i)\right),
\end{align}
where successive values $x_\epsilon(\sigma_i)$ are obtained by running the fine system (\ref{eq:alg_fast_2}) with  initial condition $x_{guess}(\sigma- \frac {\Delta} {\epsilon})$ and $l(\sigma- \frac {\Delta}{\epsilon})$. 

\par 
We discuss in Step 3 how we obtain $x_{guess}(\sigma- \frac {\Delta} {\epsilon})$. Here, $N_{t-\Delta}$ is the number of increments taken for the value of $R_{t-\Delta}^m$ to converge upto a specified value of tolerance. 

\par
\textbf{Step 2: Find the value of slow variable }\\ We use the extrapolation rule to obtain the predicted value of the slow variable at the time $t+h$: 
\begin{align}\label{eq:comp_impl_expl}
v(t+h) = v(t) + \frac{dv}{dt}(t) \, h, 
\end{align}
where $\frac{dv}{dt} (t)$ is obtained from (\ref{eq:comp_impl_1}).

\par
\textbf{Step 3: Determine the closest point projection} \\  
We assume that the closest-point projection of any fixed point in the fine state space, on the Young measure of the fine evolution, evolves slowly in any interval where the Young measure is continuous. We use this idea to define a guess  $x_{guess}(t+h-\Delta)$, that is in the basin of attraction of $\mu(t+h-\Delta)$.  The fixed point, denoted as $x^{arb}_{t-\Delta}$, is assumed  to belong to the set of points, $x_\epsilon(\sigma_i)$ for which the value of $R_{t-\Delta}^m$ in (\ref{eq:comp_impl_R2}) converged. Specifically, we make the choice of $x^{arb}_{t-\Delta}$ as $x_\epsilon(\sigma_{N_{t-\Delta}})$ where $x_\epsilon(\sigma_i)$ is defined in (\ref{eq:comp_impl_R2}) and $N_{t-\Delta}$ is defined in the discussion following it. Next, we compute the closest point projection of this point (in the Euclidean norm) on the support of the measure at $t-h-\Delta$. This is done as follows. 

Define $x^{conv}_{t-h-\Delta}$ as the point $x_\epsilon(\sigma_{N_{t-h-\Delta}})$, where  $x_\epsilon(\sigma_i)$ is defined in the discussion surrounding (\ref{eq:comp_impl_R2}) with $\sigma$ replaced by $\sigma - \frac{h}{\epsilon}$, and $N_{t-h-\Delta}$ is the number of increments taken for the value of $R_{t-h-\Delta}^m$ to converge (the value of $x^{conv}_{t-h-\Delta}$ is typically stored in memory during calculations for the time $t-h-\Delta$). The fine system (\ref{eq:alg_fast_2}) with $\sigma$ replaced by $\sigma - \frac{h}{\epsilon}$ is initiated from $x^{conv}_{t-h-\Delta}$, and we calculate the distance of successive points on this trajectory with respect to $x^{arb}_{t-\Delta}$ until a maximum number of increments have been executed. We set the maximum number of increments as $2\,N_{t-h-\Delta}$. The point(s) on this finite time trajectory that records the smallest distance from $x^{arb}_{t-\Delta}$ is defined as the closest point projection, $x^{cp}_{t-h-\Delta}$. 

Finally, the guess, $x_{guess}(t+h-\Delta)$ is given by 
\[
x_{guess}(t+h-\Delta)=2\,x^{arb}_{t-\Delta} - x^{cp}_{t-h-\Delta},
\]
for $t>0$ (see Remark associated with Step 3 in Sec. \ref{sec:theory_alg}). For $t = 0$, we set \[
x_{guess}(h-\Delta)= x^{arb}_{0} +  \frac { ( x^{arb}_{0} - x^{cp}_{-\Delta} ) } {\Delta} (h-\Delta),
\]
where $x^{arb}_{.}$ and $x^{cp}_{.}$ are defined in the first and the second paragraph respectively in Step 3 above with the time given by the subscripts. This is because the computations start at $t=-\Delta$ and we do not have a measure at $t=-h-\Delta$ and hence cannot compute $x^{cp}_{-h-\Delta}$ to be able to use the above formula to obtain $x_{guess}(t+h-\Delta)$. 

Thus, the implicit assumption is that $x_{guess}(t+h-\Delta)$ is the closest-point projection of $x^{arb}_{t-\Delta}$ on the support of the measure $\mu(t+h-\Delta)$, and that there exists a function $x^{cp}(s)$ for $t-h-\Delta \leq s\leq t+h-\Delta$ that executes slow dynamics in the time interval if the measure does not jump within it.

\par
\textbf{Step 4: Accept the measure}\\
We initiate the fine equation (\ref{eq:alg_fast_2}) at $(x_{guess}(t+h-\Delta), l(t+h-\Delta))$ and run the equation from $\sigma+\frac{h}{\epsilon}-\frac{\Delta}{\epsilon}$ to $\sigma+\frac{h}{\epsilon}$ (recall $\sigma = \frac{t}{\epsilon}$). We say that there is a match in the value of a slow variable \emph{if} the following equality holds (approximately):
\begin{equation}\label{eq:v(t+h)}
v(t+h)=\frac{1}{N'}\sum_{i=1}^{N'} m\left(x_\epsilon(\sigma_i),l_\epsilon(\sigma_i)\right),
\end{equation}
where $v(t + h)$ refers to the predicted value of the slow variable obtained from the extrapolation rule in Step 2 above. The successive values $(x_\epsilon(\sigma_i), l_\epsilon(\sigma_i))$ are obtained from system (\ref{eq:alg_fast_2}). Here, $N'=\frac {\Delta}{\epsilon \, \Delta \sigma}$ where $\Delta \sigma$ is the fine time step.
\par
If there is a match in the value of the slow variable, we accept the measure which is generated, in principle, by running the fine equation \eqref{eq:alg_fast_2} with the guess $x_{guess}(t+h-\Delta)$ and move on to the next coarse increment. 

If not, we check if there is a jump in the measure. We say that there is a jump in the measure if the value of $R_{t+h}^m$ is significantly different from the value of $R_t^m$. This can be stated as: 
\begin{align}
\left|\frac {R_{t+h}^m - R_t^m} {R_t^m}\right| \gg \frac{1}{N}\sum_n \left|\frac {R_{t-(n-1)h}^m - R_{t-nh}^m} {R_{t-nh}^m}\right|,
\end{align}
where $n$ is such that there is no jump in the measure between $t-nh$ and $t-(n-1)h$ and $N$ is the maximal number of consecutive (integer) values of such $n$. 

If there is no jump in the measure, we try different values of $x_{guess}(t+h-\Delta)$ based on different values of $x^{arb}_{t-\Delta}$ and repeat Steps 3 and 4. 

If there is a jump in the measure, we declare $v(t+h)$  to be the the right-hand-side of (\ref{eq:v(t+h)}). The rationale behind this decision is the assumption $x_{guess}(t+h-\Delta)$ lies in the basin of attraction of the measure at $t+h-\Delta$. 

{\bf Step 5: Obtain fine initial conditions for rate calculation} \\
Step 1 required the definition of $x_{guess}(t)$. We obtain it as  
\[
x_{guess}(t)=x^{arb}_{t-\Delta} + \frac {\left( x^{arb}_{t-\Delta} - x^{cp}_{t-h} \right)}{\left(h-\Delta \right)} \, \Delta,
\]
for $t>0$, and $x^{arb}_{.}$ and $x^{cp}_{.}$ are defined in the same way as in Step 3, but at different times given by the subscripts. For $t=0$, we obtain $x_{guess}(0)$, which is required to compute $R_{0}^m$, by running the fine equation (\ref{eq:alg_fast_2}) from ~$\sigma=-\frac{\Delta}{\epsilon}$ to $\sigma=0$. This is because the computations start at $t=-\Delta$ and we do not have a measure at $t=-h$ and hence cannot compute $x^{cp}_{-h}$ to be able to use the above formula to obtain $x_{guess}(t)$.

Another possible way to obtain $x_{guess}(t)$ is using the same extrapolation rule used in Step 3 to obtain $x_{guess}(t+h-\Delta)$ which means we can also obtain $x_{guess}(t)$ as $x_{guess}(t)=x^{arb}_{t-h-\Delta}+\frac{(x^{arb}_{t-h-\Delta}-x^{cp}_{t-2h-\Delta})}{h}(h+\Delta)$. But we have not used it in the computational results that follow.

We continue in this manner until an approximation of slow observables is computed on the entire time interval [$0, T_{0}$].

{\bf Discussion.} The use of the guess for fine initial conditions to initiate the fine system to compute $R^m_{t+h}$ and $R^m_{t+h-\Delta}$ is an integral part of this implementation. This allows us to systematically use the coarse evolution equation (\ref{eq:comp_impl_1}). This feature is a principal improvement over previous work \cite{tan2013coarse, tan2014md}. 

We will refer to this scheme, which is a mixture of rigorous and heuristic arguments, as Practical Time Averaging (PTA) and we will refer to results from the scheme by the same name. Results obtained solely by running the  complete system will be referred to as {\it fine} results, indicated by the superscript or subscript {\it f}  when in a formula. 

Thus, if $v$ is a scalar slow variable, then we denote the slow variable value obtained using PTA scheme as $v^{PTA}$ while the slow variable value obtained by running the fine system alone is called $v^f$.  

The speedup, $S(\epsilon)$, in compute time between the $fine$ and $PTA$ calculations is presented in the results that follow in subsequent sections. This is defined to be the ratio of the time taken by the $fine$ calculations to that by the $PTA$ calculations for an entire simulation, say consisting of $n$ steps of size $h$ on the slow time-scale. 

Let $T^{cpu}_f (\epsilon)$ and $T^{cpu}_{PTA} (\epsilon)$ be the compute times to obtain the {\it fine} and $PTA$ results  per jump on the slow time scale, respectively, for the specific value of $\epsilon$. The compute time to obtain the $PTA$ results for $n$ jumps in the slow time scale is $n T^{cpu}_{PTA}(\epsilon)$ which can be written as
\[
n T^{cpu}_{PTA}(\epsilon)=n T^{cpu}_{PTA,1}(\epsilon) + T^{cpu}_{PTA,2}(\epsilon),
\]
where $T^{cpu}_{PTA,1}(\epsilon)$ is the compute time to perform the computations mentioned in Step 1 to Step 5 for every jump in the slow time scale. Since we cannot use the formula for $x_{guess}(t)$ mentioned in Step 5 to obtain $x_{guess}(0)$ and we have to run the fine equation \eqref{eq:alg_fast_2} from $\sigma=-\frac{\Delta}{\epsilon}$ to $\sigma=0$, an additional overhead is incurred in the compute time for the $PTA$ computations which we denote as $T^{cpu}_{PTA,2}(\epsilon)$. Thus
\[
S(\epsilon) = \frac{n T^{cpu}_f (\epsilon)}{n T^{cpu}_{PTA}(\epsilon)} \approx \frac{T^{cpu}_{f}(\epsilon)}{T^{cpu}_{PTA,1}(\epsilon)}.
\]
for large $n$.

Error in the PTA result is defined as: 
\begin{align}
Error(\%)=\frac{v^{PTA}-v^f}{v^f} \times 100.
\end{align}

We obtain $v^f$ as follows: 

{\bf Step 1:} We run the fine system (\ref{eq:alg_fast_2}) from $\sigma=-\frac{\Delta}{\epsilon}$ to $\sigma=\frac{T_0}{\epsilon}$ using initial conditions ($x_0$, $l_0$) to obtain ($x_\epsilon(\sigma_i),l_\epsilon(\sigma_i)$) where $\sigma_i = i \, \Delta \sigma$ and $i \in \mathbb Z_+$ and $i \leq \frac{T_0+\Delta}{\epsilon \, \Delta \sigma}$. 

{\bf Step 2:} We calculate $v^f(t)$ using: 
\begin{align}\label{eq:v_fine}
v^f(t)=\frac{1}{N'}\sum_{i=N^0(t)}^{N^0(t)+ N'} m\left(x_\epsilon(\sigma_i),l_\epsilon(\sigma_i)\right),
\end{align} 
where $N'=\frac {\Delta}{\epsilon \, \Delta \sigma}$ and $N^0(t)=\frac{t+\Delta}{\epsilon \, \Delta \sigma}$ where $\Delta \sigma$ is the fine time step. 

{\bf Remark.} If we are aiming to understand the evolution of the slow variables in the slow time scale, we need to calculate them, which we do in Step 2. However, the time taken in computing the average of the state variables in Step 2 is much smaller compared to the time taken to run the fine system in Step 1. We will show this in the results sections that follow. 

{\bf Remark.} All the examples computed in this paper employ $H$-observables. When, however, orthogonal observables are used, the time taken to compute their values using the PTA scheme ($T^{cpu}_{PTA}$) will not depend on the value of $\epsilon$.  

\section{Example I: Rotating planes}
\noindent
Consider the following four-dimensional system, where we denote by $x$ the vector $x = (x_{1},x_{2},x_{3},x_{4})$.
\begin{align}\label{eq:9.1}
\frac{dx}{dt} = \frac{F(x)}{\epsilon} + G(x),
\end{align}
\noindent
where:
\begin{align}\label{eq:9.2}
F(x) = ((1-\lvert x\rvert)x+\gamma(x))
\end{align}
\noindent
with
\begin{align}\label{eq:9.3}
\gamma(x) = (x_{3},x_{4},-x_{1},-x_{2}).
\end{align}
\noindent
The drift may be determined by an arbitrary function $G(x)$. For instance, if we let
\begin{align}\label{eq:9.4}
G(x) = (-x_{2}, x_{1}, 0, 0),
\end{align}
\noindent
then we should expect nicely rotating two-dimensional planes. A more complex drift may result in a more complex dynamics of the invariant measures, namely the two dimensional limit cycles.

\par
\subsection{Discussion} 
The right hand side of the fast equation has two components. The first drives each point $x$ which is not the origin{\color{blue},} toward the sphere of radius 1. The second, $\gamma(x)$, is perpendicular to $x$. It is easy to see that
the sphere of radius 1 is invariant under the fast equation. For any initial condition $x_0$ on the sphere of radius 1, the fast time equation is
\begin{equation} 
\begin{aligned}
\dot{ {\bf x} } = \begin{pmatrix}  0 & 0 & 1 & 0 \\
								 0 & 0 & 0 & 1 \\
								-1 & 0 & 0 & 0 \\
								  0 &-1 & 0 & 0 \end{pmatrix} {\bf x} ~. 
\end{aligned}
\label{ex1_rate}
\end{equation}
It is possible to see that the solutions are periodic, each contained in a two dimensional subspace. An explicit solution (which we did not used in the computations) is
\begin{equation} 
{\bf x} = cos \left( t \right) \begin{pmatrix} x_{0,1} \\ x_{0,2} \\ x_{0,3} \\ x_{0,4} \end{pmatrix}
             + sin \left( t \right) \begin{pmatrix} x_{0,3} \\ x_{0,4} \\ - x_{0,1} \\ -x_{0,2} \end{pmatrix}.
\label{ex1_sol}
\end{equation}
Thus, the solution at any point of time is a linear combination of ${\bf x}_0$ and $ { \gamma} ( {\bf x}_0 )$ and stays in the plane defined by them. Together with the previous observation we conclude that the limit occupational measure of the fast dynamics should exhibit oscillations in a two-dimensional subspace of the four-dimensional space. The two dimensional subspace itself is drifted by the drift $G(x)$. The role of the computations is then to follow the evolution of the oscillatory two dimensional limit dynamics.

\par 
We should, of course, take advantage of the structure of the dynamics that was revealed in the previous paragraph. In particular, it follows that three observables of the form $r(\mu)$ given in (\ref{eq:4.6}), with $e_{1}$, $e_{2}$ and $e_{3}$ being unit vectors in $\mathbb{R}^{4}$, determine the invariant measure. They are not orthogonal (it may be very difficult to find orthogonal observables in this example), hence we may use, for instance, the $H$-observables introduced in (\ref{eq:4.7}).
\par
It is also clear that the circles that determine the invariant measures move  smoothly with the planes. Hence employing observables that depend smoothly on the planes would imply that conditions (\ref{eq:8.1}) and (\ref{eq:8.2}) hold, validating the estimates of Theorem 8.2.

\par
\subsection{Results} 
We chose the slow observables to be the averages over the limit cycles of the four rapidly oscillating variables and their squares since we want to know how they progress. We define the variables $w_i=x_i^2$ for $i=1,2,3$ and $4$. The slow variables are $x_1^f$, $x_2^f$, $x_3^f$, $x_4^f$ and $w_1^f $ , $w_2^f $ , $w_3^f$, $w_4^f $. The slow variable $x_1^f$ is given by (\ref{coarse_obs_impl}) with $m(x)=x_1$. The slow variables $x_2^f$, $x_3^f$ and $x_4^f$ are defined similarly. The slow variable $w_1^f$ is given by (\ref{coarse_obs_impl}) with $m(x)=w_1$. The slow variables $w_2^f$, $w_3^f$ and $w_4^f$ are defined similarly (we use the superscript $f$, that indicates the {\it fine} solution, since in order to compute these observables we need to solve the entire equation, though on a small interval). We refer to the $PTA$ variables as $x_1^{PTA}$, $x_2^{PTA}$, $x_3^{PTA}$, $x_4^{PTA}$ and $w_1^{PTA} $ , $w_2^{PTA} $ , $w_3^{PTA}$, $w_4^{PTA} $. A close look at the solution (\ref{ex1_sol}) reveals that the averages, on the limit cycles, of the fine variables, are all equal to zero, and we expect the numerical outcome to reflect that. The average of the squares of the fine variables evolve slowly in time.  
In our algorithm, non-trivial evolution of the slow variable does not play a role in tracking the evolution of the measure of the complete dynamics. Instead they are used only to accept the slow variable (and therefore, the measure) at any given discrete time as valid according to Step 4 of Section \ref{impl_algo}. It is the device of choosing the initial guess in Step 3 and Step 5 of Section \ref{impl_algo} that allows us to evolve the measure discretely in time.

\begin{figure*}[!h]
    \centering
    \begin{subfigure}[h]{0.4\textwidth}
        \centering
        \includegraphics[height=3.1in]{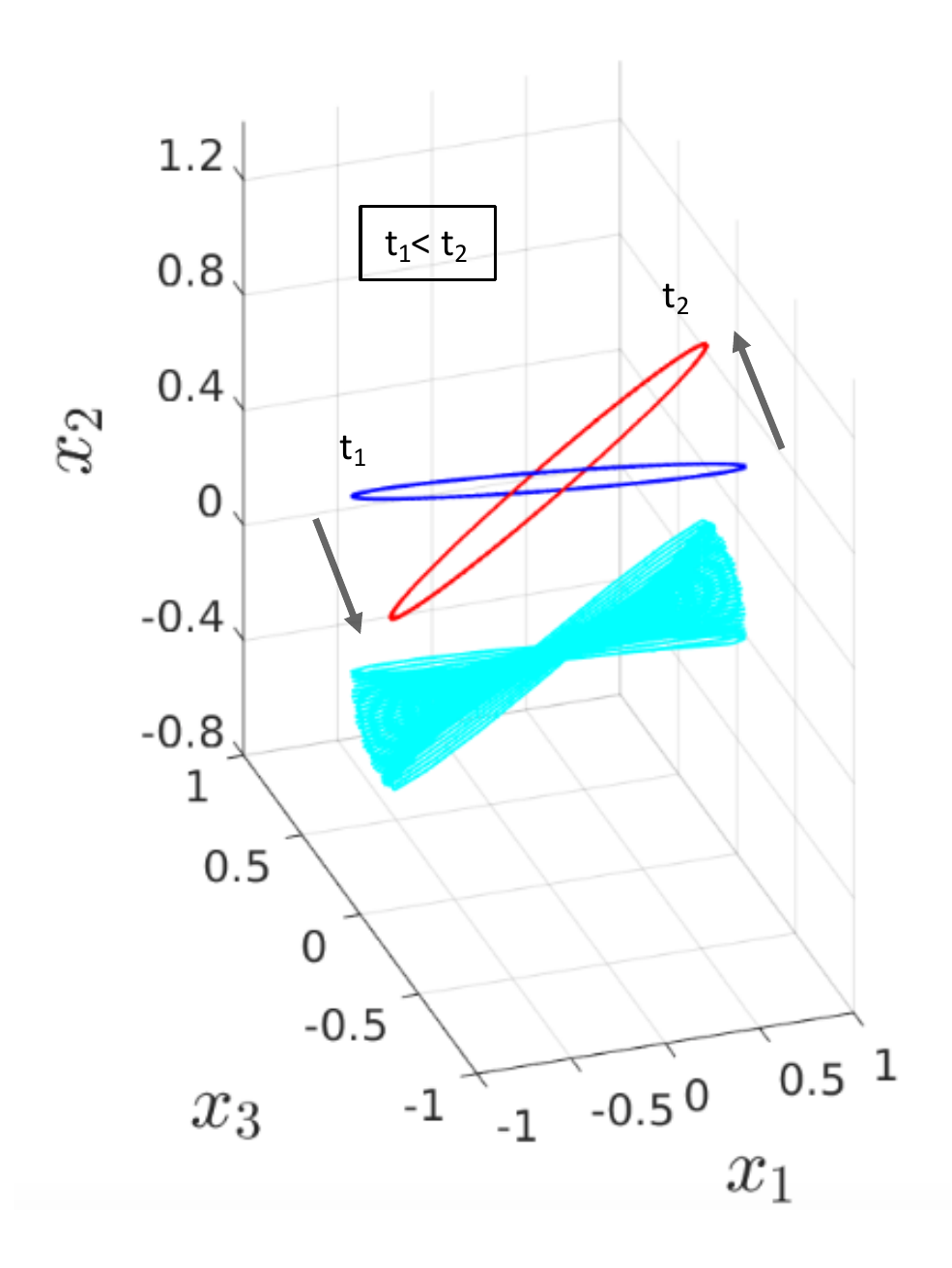}
        \caption{}
    \end{subfigure}%
    ~ 
    \begin{subfigure}[h]{0.6\textwidth}
        \centering
        \includegraphics[height=3.1in]{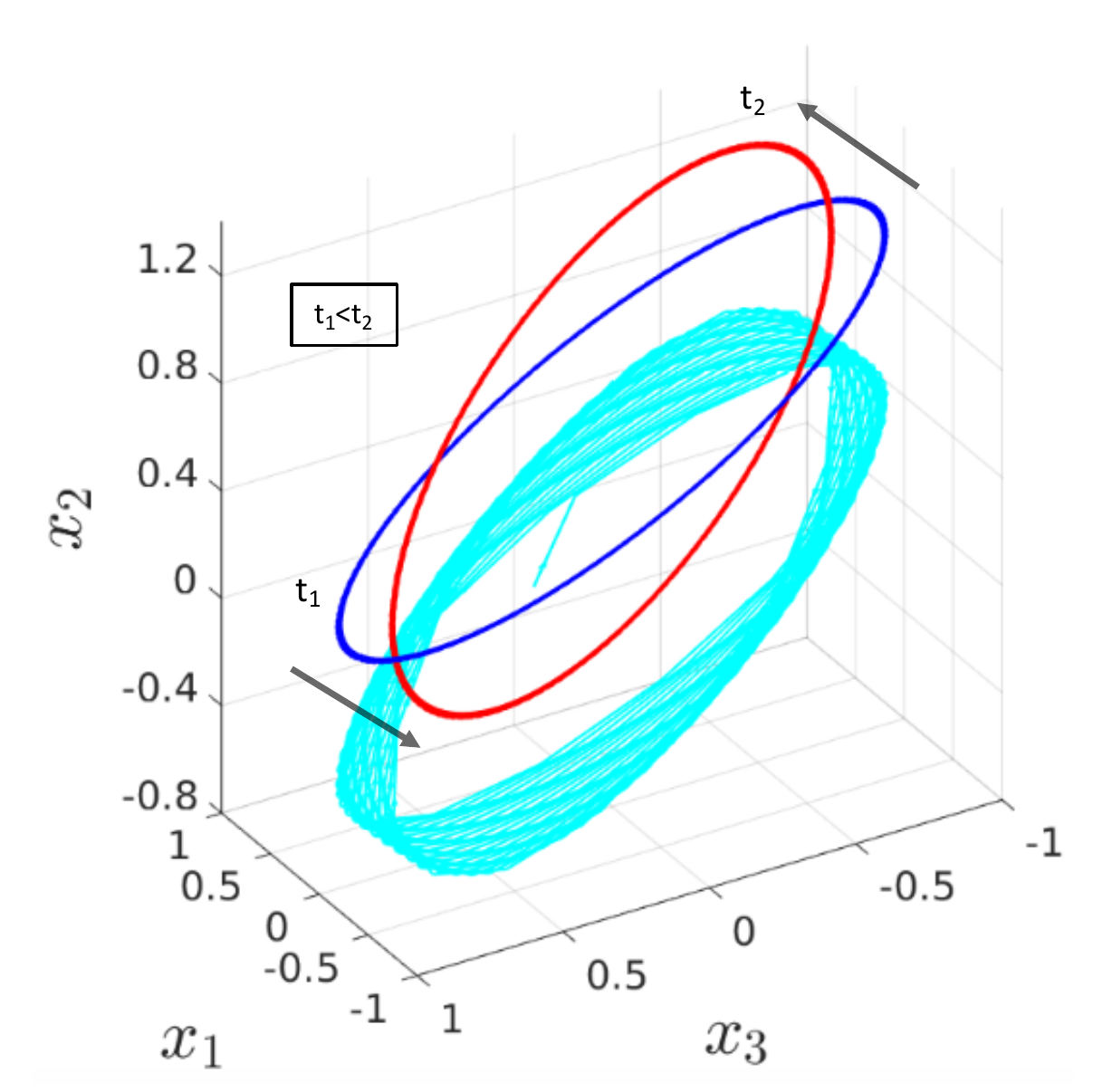}
        \caption{}
        \end{subfigure}
        \caption{\textit{ The curve in cyan shows the phase space diagram of $x_1$, $x_2$ and $x_3$ obtained by running the system \eqref{eq:9.1} to $t=2$ with $\epsilon=10^{-7}$. Part \textup{(}a\textup{)} and Part\textup{(}b\textup{)} show different views of the phase space diagram. The blue curve shows the portion of the phase portrait obtained around time $t_1$ while the red curve shows the portion around a later time $t_2$. \textup{(}For interpretation of the references to color in this figure legend, the reader is referred to the web version
of this article.\textup{)}}}   
\label{fig:ex1_phase}
\end{figure*}

%\begin{figure}[!h]
%\centering
%\begin{minipage}{.45\textwidth}
%  \centering
%  \includegraphics[width=\linewidth]{figures/ex1_x1_fine_evol1.pdf}
%  \caption{\textit{Phase space diagram of $x_1$, $x_2$ and $x_3$}}
%  \label{fig:ex1_phase}
%\end{figure}

Fig. \ref{fig:ex1_phase} shows the phase space diagram of $x_1$, $x_2$ and $x_3$.

\begin{figure}[!h]
\centering
\begin{minipage}{.45\textwidth}
  \centering
  \includegraphics[width=\linewidth]{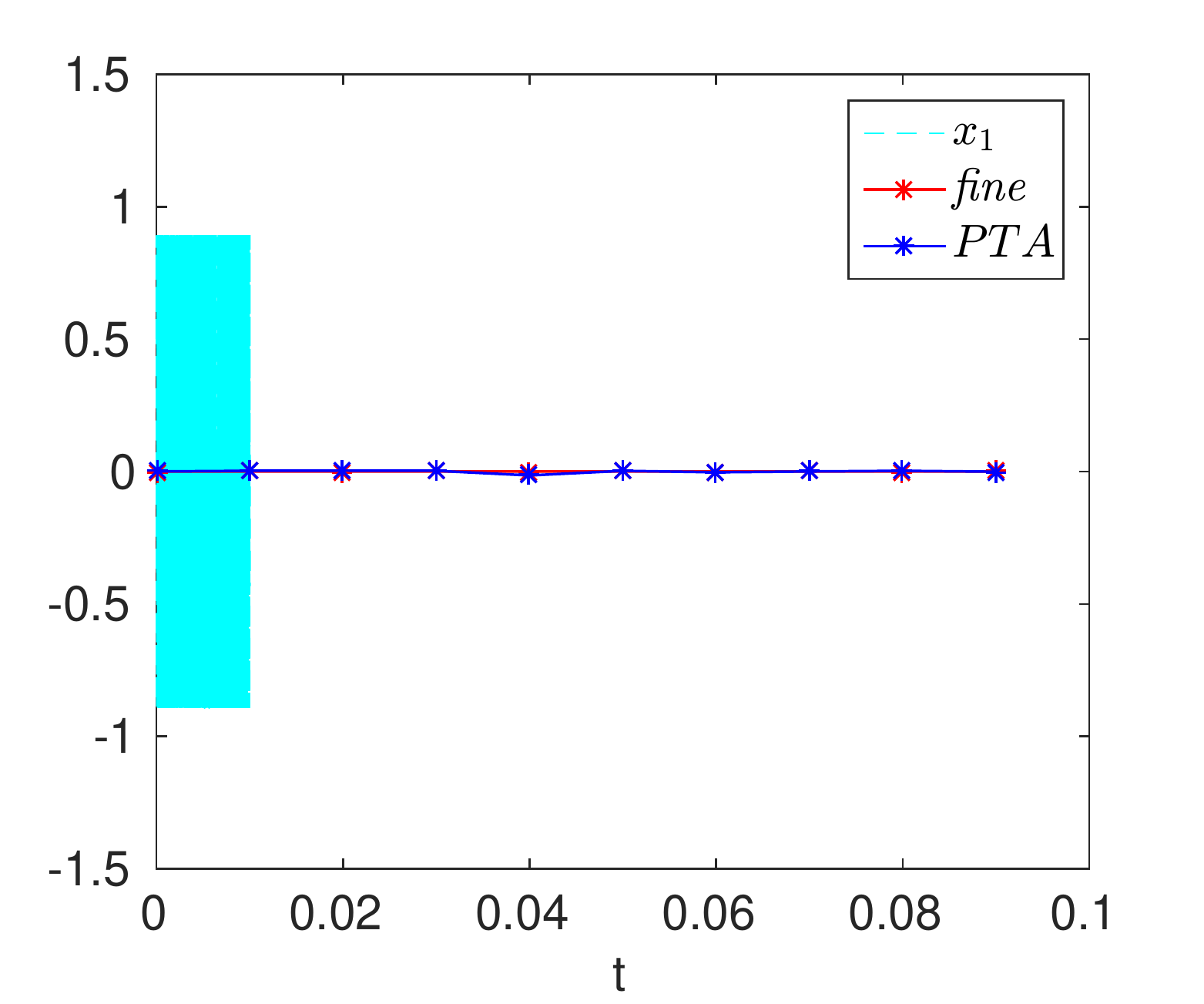}
  \caption{\textit{ The rapidly oscillating solution of the full equation} of $x_1$ is given by the plot marked $x_1$ which shows rapid oscillations around the {\it fine} and $PTA$ values (which is, as expected, equal to 0). The $PTA$ and the {\it fine} results overlap.}  
  \label{fig:ex1_x1_fine}
\end{minipage}%
\hfill
\begin{minipage}{.45\textwidth}
  \centering
  \includegraphics[width=\linewidth]{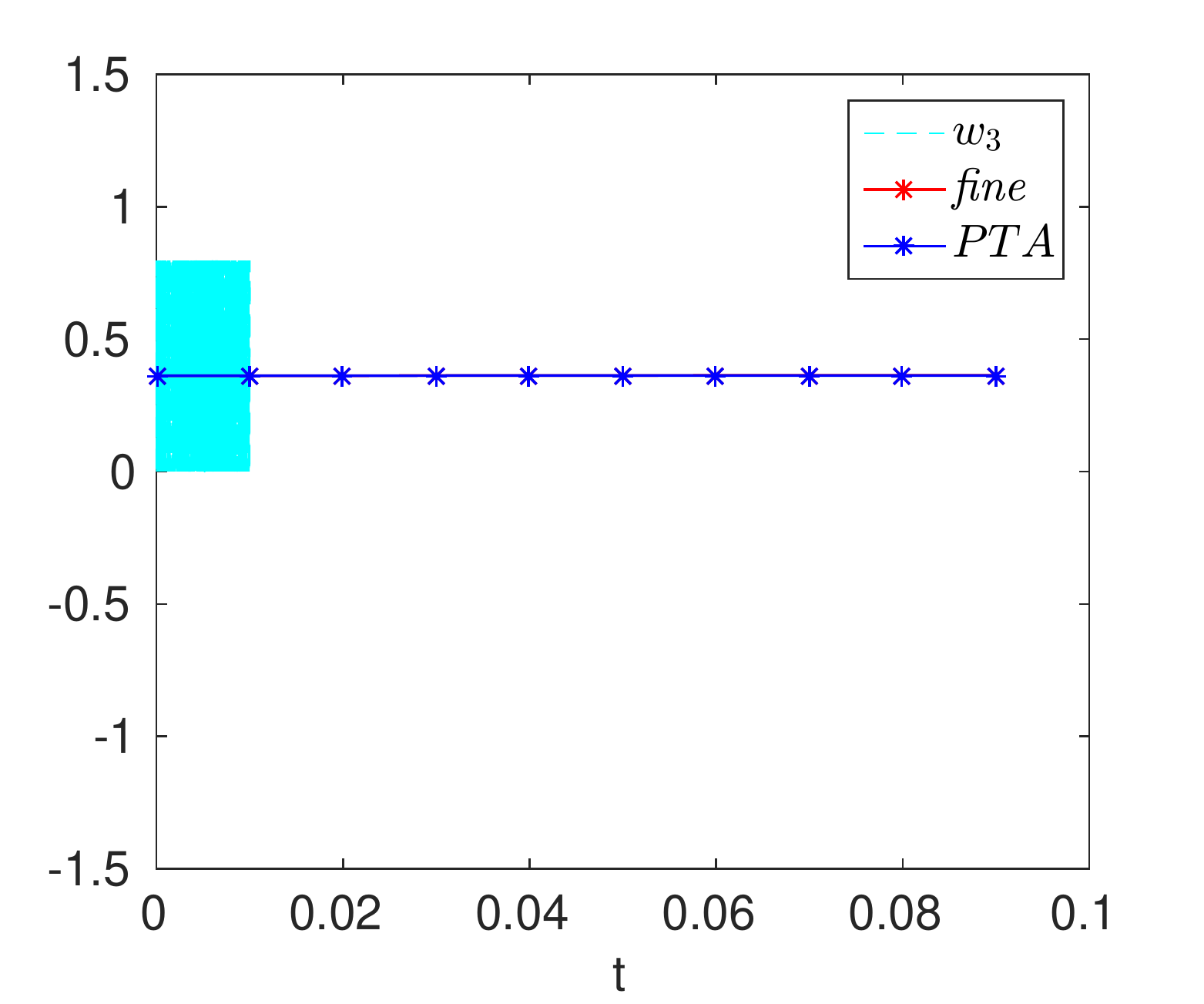}
  \caption{\textit{  The rapidly oscillating solution of the full equation of $w_3$ is given by the plot marked $w_3$. The drift in the {\it fine} and $PTA$ values cannot be seen on the given scale. But the drift is visible in Fig. \ref{fig:ex1_x3_sq}}. The $PTA$ and the {\it fine} results overlap.}  
  \label{fig:ex1_x3sq_fine}
\end{minipage}
\end{figure}

Fig. \ref{fig:ex1_x1_fine} shows the rapid oscillations of the rapidly oscillating variable $x_1$ and the evolution of the slow variable $x_1^f$. Fig. \ref{fig:ex1_x3sq_fine} shows the rapid oscillations of $w_3$ and the evolution of the slow variable $w_3^f$. We find that $x_3$ and $x_4$  evolve exactly in a similar way as $x_1$ and $x_2$ respectively. We find from the results that $x_3^f$, $x_4^f$, $w_3^f$ and $w_4^f$ evolve exactly similarly as $x_1^f$, $x_2^f$, $w_1^f$ and $w_2^f$ respectively.

\begin{figure}[!h]
\centering
\begin{minipage}{.5\textwidth}
  \centering
  \includegraphics[width=\linewidth]{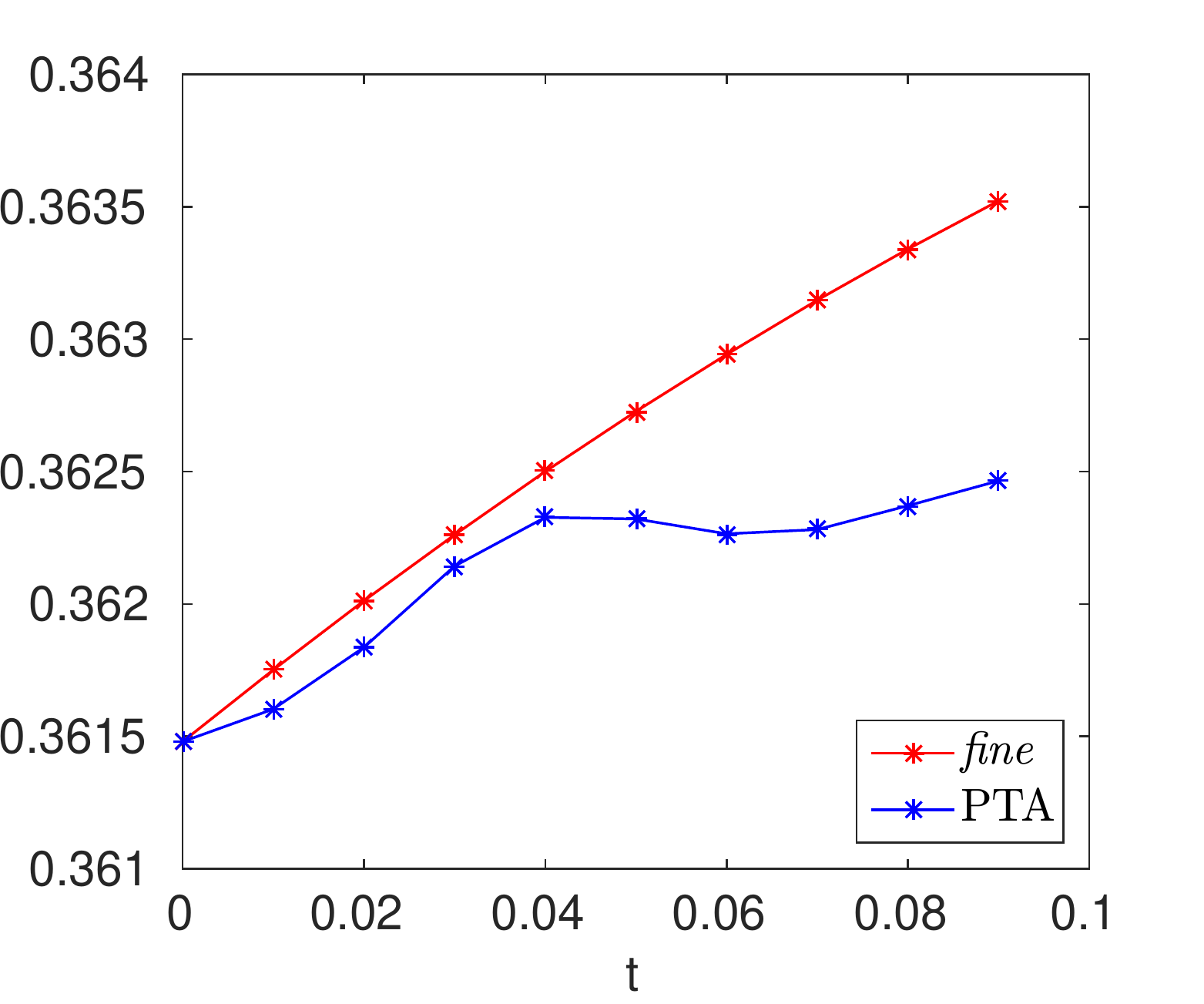}
  \caption{\textit{Evolution of $w_3^f$.} }
  \label{fig:ex1_x3_sq}
\end{minipage}%
\hfill
\begin{minipage}{.5\textwidth}
  \centering
  \includegraphics[width=\linewidth]{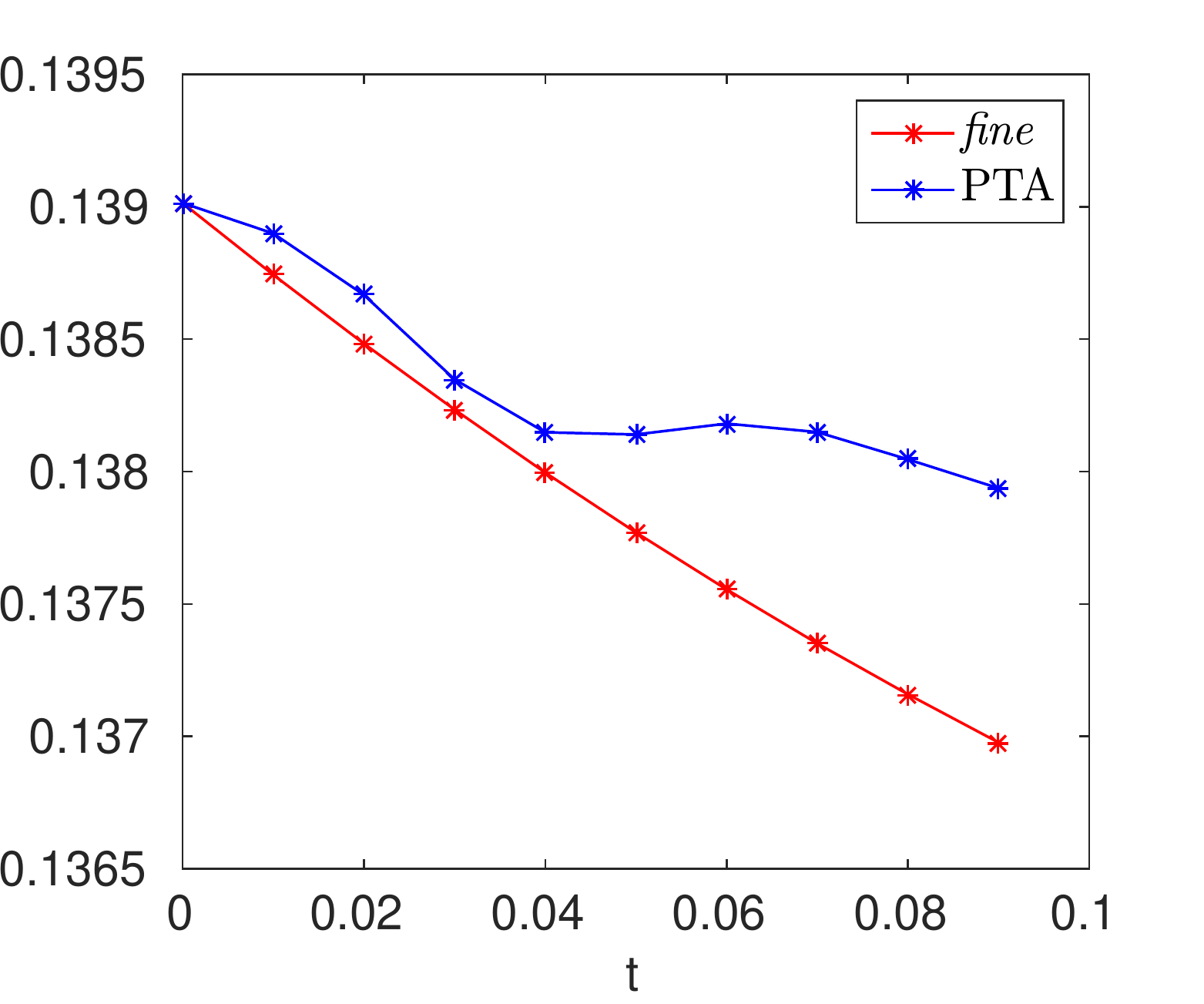}
  \caption{\textit{Evolution of $w_4^f$.} }
  \label{fig:ex1_x4_sq}
\end{minipage}
\end{figure}

\begin{figure}[!h]
\centering
\includegraphics[scale = .4]{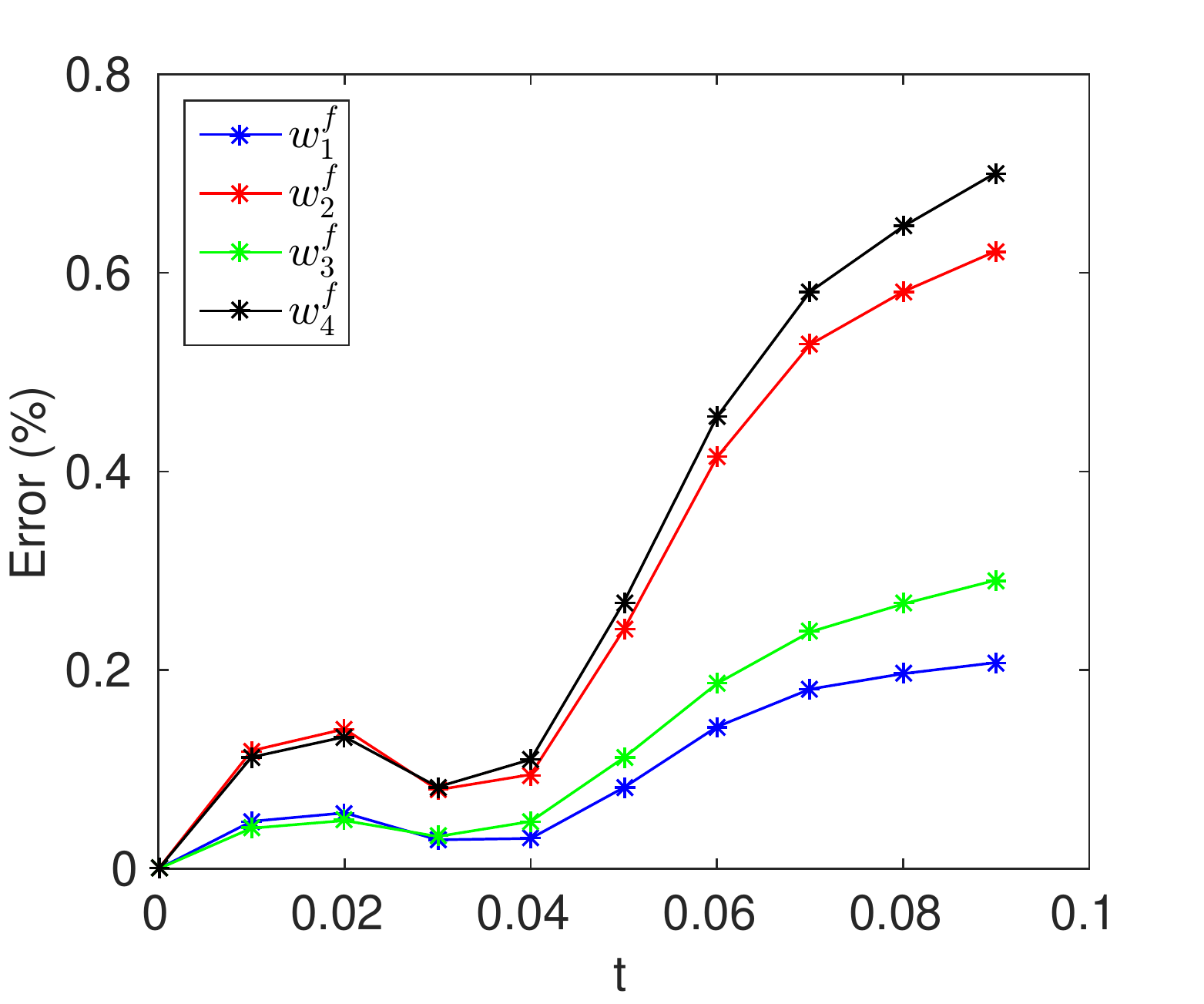}
\caption{\textit{Example I - Error}. }
\label{fig:ex1_rel_err}
\end{figure}
The comparison between the {\it fine} and the PTA results of the slow variables $w_3^f$ and $w_4^f$ are shown in Fig. \ref{fig:ex1_x3_sq} and Fig. \ref{fig:ex1_x4_sq} (we have not shown the evolution of $w_1^f$ and $w_2^f$ since they evolve exactly similarly to $w_3^f$ and $w_4^f$ respectively). The error in the PTA results are shown in Fig. \ref{fig:ex1_rel_err}. Since the values of $x_1^f$, $x_2^f$, $x_3^f$ and $x_4^f$ are very close to $0$, we have not provided the error in PTA results for these slow variables. 
\par 
{\bf Savings in computer time}
\begin{figure}[!h]
\centering
\includegraphics[scale = .4]{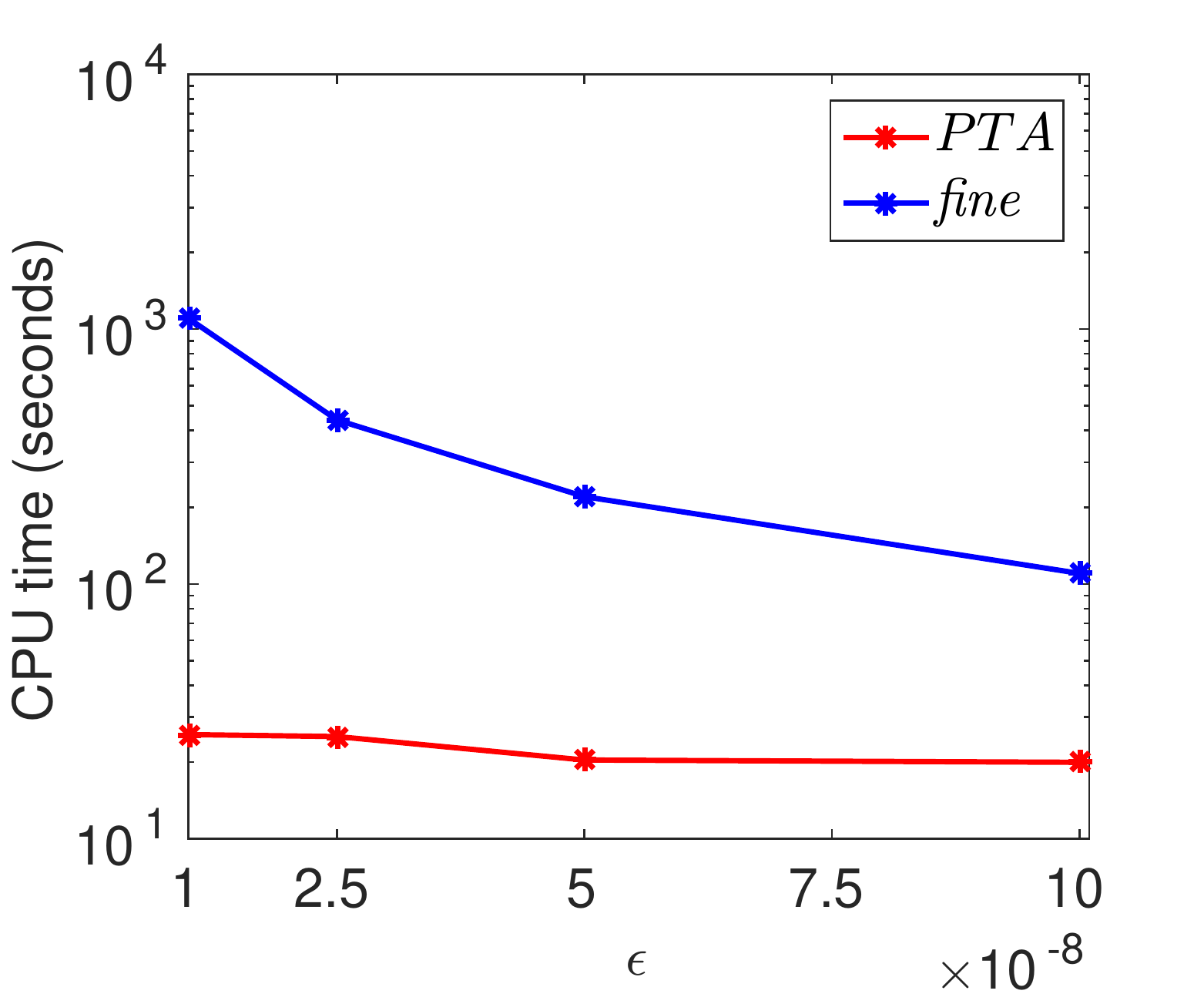}
\caption{\textit{Example I - Compute time comparison for simulations spanning $t=0.01$ to $t=0.02$.} }
\label{fig:ex1_cputime}
\end{figure}

In Fig. \ref{fig:ex1_cputime}, we see that as $\epsilon$ decreases, the compute time for the fine run increases very quickly while the compute time for the PTA run increases relatively slowly. The compute times correspond to simulations spanning $t=0.01$ to $t=0.02$ with $\Delta=0.001$. The speedup in compute time, $S$, obtained as a function of $\epsilon$, is given by the following polynomial: 
\begin{align}
S(\epsilon)= 73.57 - 3.70 \times 10^{9} \, \epsilon + 6.76 \times 10^{16} \, \epsilon^2 - 3.74\times 10^{23} \, \epsilon^3.
\end{align}

The function $S(\epsilon)$ is an interpolation of the computationally obtained data to a cubic polynomial. A more efficient calculation yielding higher speedup is to calculate the slow variable $v$ using Simpson's rule instead of using \eqref{eq:v(t+h)}, by employing the procedures outlined in section \ref{prb2:res_num_case1}, section \ref{prb2:res_num_case2} and the associated \textit{Remark} in Appendix \ref{prb2:res_case2}. In this problem, we took the datapoint of $\epsilon=10^{-8}$ and obtained $S(10^{-8}) = 43$. This speedup corresponds to an accuracy of $0.7\%$ error. However, as $\epsilon$ decreases and approaches zero, the asymptotic value of $S$ becomes 74. 

\section{Example II: Vibrating springs}\label{vibrate}
%%
%\begin{figure}[t]
%\centering
%\includegraphics[width=6.5in, height=1.5in]{figure_1_cspde.pdf}
%\caption{Sketch of mechanical system for Example II}\label{fig:1}
%\end{figure}
%%
%%

Consider the mass-spring system in Fig. \ref{fig:prb2}. The governing system of equations for the system in dimensional time is given in Appendix \ref{prb2:equations}. The system of equations posed in the slow time scale is
\begin{align}\label{eq:10.3}
\epsilon \frac{dx_{1}}{d t} &=  T_f \, y_{1}\nonumber\\
\epsilon \frac{dy_{1}}{d t}&= - T_f \left(\frac{k_{1}}{m_{1}}(x_{1} - w_{1})-\frac{\eta}{m_{1}}(y_{2}-y_{1})\right)\nonumber\\
\epsilon \frac{dx_{2}}{d t}&=  T_f \, y_{2} \nonumber\\
\epsilon \frac{dy_{2}}{d t}&= - T_f \left(\frac{k_{2}}{m_{2}}(x_{2}-w_{2}) + \frac{\eta}{m_{2}}(y_{2}-y_{1})\right)\nonumber\\
 \frac{dw_{1}}{d t}&= T_s \, L_{1}(w_{1})\nonumber\\
 \frac{dw_{2}}{d t}&= T_s \, L_{2}(w_{2})~.
\end{align}

\begin{figure}[H]
\centering
\includegraphics[scale = .5]{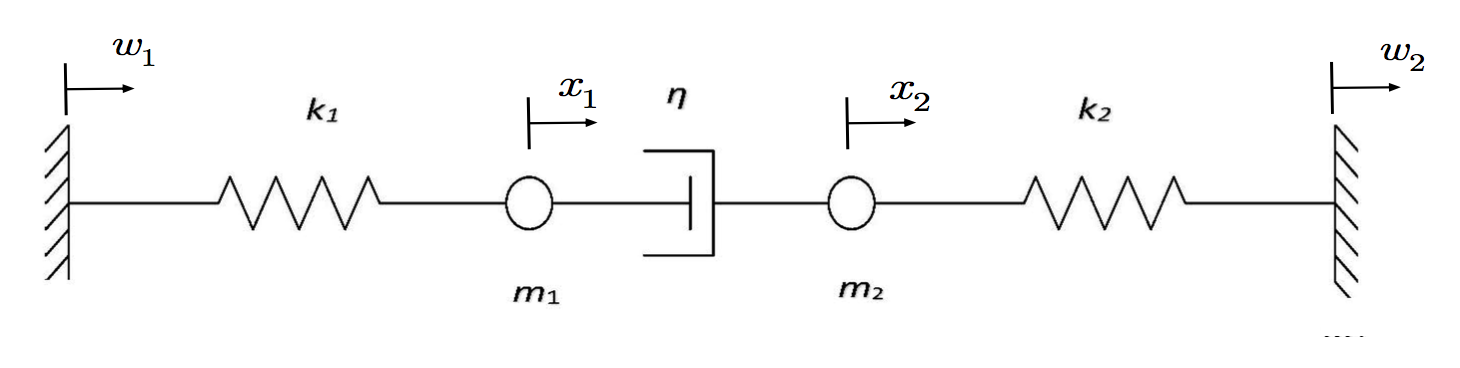}
\caption{Sketch of the mechanical system for problem II.}
\label{fig:prb2}
\end{figure}

The derivation of (\ref{eq:10.3}) from the system in dimensional time is given in Appendix \ref{prb2:equations}.  The small parameter $\epsilon$ arises from the ratio of the fast oscillation of the springs to the slow application of the load. A closed-form solution to (\ref{eq:10.3}) can be computed, but in the present study the solution will be used just for verifying the computations, and will not be used in computations themselves. The closed-form solutions are presented in Appendices \ref{prb2:res_case1} and \ref{prb2:res_case2}. 
%\end{minipage}

%%
%%
\par
\subsection{Discussion}\label{ex2:discussion}
\begin{itemize} 
\item 
For a fixed value of the slow dynamics, that is, for fixed positions $w_{1}$ and $w_{2}$ of the walls, the dynamics of the fine equation is as follow. If $\frac{k_{1}}{m_{1}} \neq {\frac{k_{2}}{m_{2}}}$, then all the energy is dissipated, and the trajectory converges to the origin (the reason behind this behavior is explained in the Remark of Appendix \ref{prb2:res_case1}). If the equality holds, only part of the energy possessed by the initial conditions is dissipated, and the trajectory converges to a periodic one (in rare cases it will be the origin), whose energy is determined by the initial condition (the reason behind this behavior is explained in Case 2.1 and Case 2.2 of Section \ref{prb2:res_num_case2} and in the Remark of Appendix \ref{prb2:res_case2}). The computational challenge is when fast oscillations persist. Then the limiting periodic solution determines an invariant measure for the fast flow. When the walls move, slowly, the limit invariant measure moves as well. The computations should detect this movement. However, if the walls move very slowly, there is a possibility that in the limit the energy does not change at all as the walls move. 
\par
Notice that the invariant measure is not determined by the position of the walls, and additional slow observables should be incorporated. A possible candidate is the total energy stored in the invariant measure. Since the total energy is constant on the limit cycle, it forms an orthogonal observable as described in section 4. Its extrapolation rule is given by (\ref{eq:4.1}). In order to apply (\ref{eq:4.1}) one has to derive the effect of the movement of the walls on the observable, namely, on the total energy.
\par
Two other observables could be the average kinetic energy and the average potential energy on the invariant measure. In both cases, the form of $H$-observables should be employed, as it is not clear how to come up with an extrapolation rules for these observables.
\par
It is clear that with the three observables just mentioned, if the two forcing elements $L_{1}$ and $L_{2}$ in \eqref{eq:10.3} are Lipschitz, then conditions (\ref{eq:8.1}) and (\ref{eq:8.2}) are satisfied, and, consequently,
the conclusion of Theorem 8.2 holds.

\item 
We define \emph{kinetic energy} ($K$), \emph{potential energy} ($U$) and \emph{reaction force} on the right wall ($R_2$) as:
\begin{align}\label{ex2_statefn}
K(\sigma)&=\frac{1}{2} \left(m_1 \, {y_{1,\epsilon}(\sigma)}^2 + m_2 \, {y_{2,\epsilon}(\sigma)}^2\right) \nonumber\\
U(\sigma)&=\frac{1}{2} k_1 {(x_{1,\epsilon}(\sigma)-w_{1,\epsilon}(\sigma))}^2 + \frac{1}{2} k_2 {(x_{2,\epsilon}(\sigma)-w_{2,\epsilon}(\sigma))}^2 \\
R_2(\sigma) &= -k_2\left(x_{2,\epsilon}(\sigma) - w_{2,\epsilon}(\sigma)\right)\nonumber.
\end{align}

\item
The H-observables that we obtained in this example are the \emph{average kinetic energy} ($K^f$), \emph{average potential energy} ($U^f$) and \emph{average reaction force} on the right wall ($R_2^f$) which are calculated as:
\begin{align}\label{ex2_hobs}
{K^f}(t)&=\frac{1}{N'} \sum_{i=1}^{N'} K(\sigma_i) \nonumber\\
{U^f}(t)&=\frac{1}{N'} \sum_{i=1}^{N'} U(\sigma_i) \\
{R_2^f}(t)&=\frac{1}{N'} \sum_{i=1}^{N'} R_2(\sigma_i), \nonumber  
\end{align}
where $N'$ is defined in the discussion following (\ref{eq:v(t+h)}) and successive values $x_{1,\epsilon}(\sigma_i)$, $x_{2,\epsilon}(\sigma_i)$, $y_{1,\epsilon}(\sigma_i)$ and $y_{2,\epsilon}(\sigma_i)$ are obtained by solving the fine system associated with \eqref{eq:10.3} (see \eqref{eq:10.4} in Appendix \ref{prb2:equations}) with appropriate initial conditions which is discussed in detail in Step 3 of Section \ref{impl_algo}. The computations are done when $L_1(w_1) = 0$ and $L_2(w_2) = c_2$.

\item 
To integrate the fine system (\ref{eq:10.4}), we use a modification of the velocity Verlet integration scheme to account for damping (given in \cite{Sandvik_2016_py502}). This is done so that the energy of the system does not diverge in time due to energy errors of the numerical method. 
\item
As we will show in Section \ref{prb2:res_num_case1} (where we show results for the case corresponding to the condition $\frac{k_1}{m_1} \neq \frac{k_2}{m_2}$ which we call Case 1) and Section \ref{prb2:res_num_case2} (where we show results for the case corresponding to the condition $\frac{k_1}{m_1} = \frac{k_2}{m_2}$ which we call Case 2) respectively, in Case 1, the fine evolution converges to a singleton (in the case without forcing) while in Case 2, the fine evolution generically converges to a limit set that is not a singleton (which will be shown in Case 2.2 in Section \ref{prb2:res_num_case2}), which shows the distinction between the two cases. This has significant impact on the results of average kinetic and potential energy. Our computational  scheme requires no \emph{a-priori} knowledge of these important distinctions and predicts the correct approximations of the limit solution in all of the cases considered. 

\item 
For the sake of comparison with our computational approximations, in Appendix C we provide solutions to our system \eqref{eq:10.3} corresponding to the Tikhonov framework \cite{Tikhonov_1985_DE} and the quasi-static assumption commonly made in solid mechanics for mechanical systems forced at small loading rates. We show that the quasi-static assumption does not apply for this problem. The Tikhonov framework applies in some situations and our computation results are consistent with these conclusions. As a cautionary note involving limit solutions (even when valid), we note that evaluating nonlinear functions like potential and kinetic energy on the weak limit solutions as a reflection of the limit of potential and kinetic energy along sequences of solutions of \eqref{eq:10.3} as $\epsilon \to 0$ (or equivalently $T_s \to \infty$) does not make sense in general, especially when oscillations persist in the limit. Indeed, we observe this for all results in Case 2.  

\item 
All results shown in this section are obtained from numerical calculations, with no reference to the closed-form solutions. The closed-form solution for Case 1 is derived in Appendix \ref{prb2:res_case1} while the closed-form solution for Case 2 is derived in Appendix \ref{prb2:res_case2}. 
\end{itemize}

\subsection{Results - Case 1 : $\left(\frac{k_1}{m_1} \neq \frac{k_2}{m_2}\right)$}\label{prb2:res_num_case1}

We do not use the explicit limit dynamics displayed in Appendix \ref{prb2:res_case1}. Rather, we proceed with the computations employing the kinetic and potential energies as our $H$-observables.

\par 
All simulation parameters are grouped in Table \ref{tab:simulation_details_uneql}. The total physical time over which the simulation runs is $T_0 T_s$, where $T_0$ is defined in Section \ref{sec:theory_alg} (in all computed problems here, we have chosen $T_0=1$). The PTA computations done in this section are with the load fixed while calculating $R^m_t$ and $R^m_{t-\Delta}$ using \eqref{eq:comp_impl_R1} and \eqref{eq:comp_impl_R2} respectively (by setting $\frac{dl_\epsilon}{d\sigma}=0$ in \eqref{eq:alg_fast_2}). The slow variable value  ($v(t+h)$ in \eqref{eq:v(t+h)}) is calculated using Simpson's rule as described in \textit{Remark} in Appendix \ref{prb2:res_case2}. The $PTA$ results and the {\it closed-form} results (denoted by ``$cf$" in the superscript) match for all values of $t$ ($K^{PTA}= 10^{-10} \approx 0 = K^{cf} $, $U^{PTA} = 10^{-10} \approx 0 = U^{cf}$ and $R_2^{PTA} = 10^{-5} \approx 0 = R_2^{cf}$ - note that following the discussion around \eqref{non-dim} in Appendix \ref{prb2:res_case2}, all the results presented here are non-dimensionalized). In this case, the results from the Tikhonov framework match with our computational approximations. This is because after the initial transient dies out, the whole system displays slow behavior in this particular case. However, the solution under the quasi-static approximation does not match our computational results (even though the loading rate is small), and we indicate the reason for its failure in  Appendix \ref{prb2:res_case1}. 

\begin{table}[h]
\centering
\begin{tabular}[h]{|c|c|c|}
\hline
Name &   Physical definition & Values  \\
\hline
$k_1$     &  Stiffness of left spring      &   $ 10^7 \,\mathit{N/m}$  \\
$k_2$     &  Stiffness of right spring      &   $ 10^7 \,\mathit{N/m}$  \\
$m_1$     &  Left mass      &   $1\,\mathit{kg}$  \\
$m_2$     &  Right mass      &   $2\,\mathit{kg}$  \\
$\eta$     &  Damping coefficient of dashpot      &   $ 5 \times 10^3 \,\mathit{N \, s/m}$  \\
$c_2$     &  Velocity of right wall  &   $10^{-6}$ = $\frac{0.01}{10^4}\,\mathit{m/s}$  \\
$ h $      &  Jump size in slow time scale & $ 0.25 $ \\
$ \Delta $  & Parameter used in rate calculation & $ 0.05 $ \\
\hline
\end{tabular}
\caption{Simulation parameters.}
\label{tab:simulation_details_uneql}
\end{table}

\subsubsection{Power Balance}\label{ex2_power_balance} 
It can be show from \eqref{eq:10.1} that  
\begin{align}\label{power_balance}
{1 \over T_s} \frac{d}{dt} \left( \frac{1}{ 2} m_1 y_1^2 + \frac{1}{ 2} m_2 y_2^2 + \frac{1}{ 2} m_{w_1} {v_{w_1}}^2 + \frac{1}{ 2} m_{w_2} {v_{w_2}}^2 \right)& \nonumber \\
+ {1 \over T_s} \frac{d}{dt} \left( \frac{1}{ 2} k_1 {(x_1 - w_1)}^2 + \frac{1}{ 2} k_2 {(x_2 - w_2)}^2 \right) 
&= R_1  v_{w_1} + R_2 v_{w_2} - \eta {(y_2 - y_1)}^2,
\end{align}
where $v_{w_1}= {1 \over T_s} \frac{d w_1}{dt}$ and $v_{w_2}= {1 \over T_s} \frac{d w_2}{dt}$. Equation \eqref{power_balance} is simply the statement that at any instant of time the rate of change of kinetic energy and potential energy is the external power supplied through the motion of the walls less the power dissipated as viscous dissipation. This means that the sum of the kinetic and potential energy of the system, is equal to the sum of the initial kinetic and potential energy, plus the integral of the external power supplied to the system, minus the viscous dissipation.  

The {\it fine} solution indicates that, for $c_2$ small, the dashpot kills all the initial potential and kinetic energy supplied to the system. The two springs get stretched based on the value of $c_2$. The stretches remain fixed for large times on the fast time scale and the mass $m_2$ and the right wall move with the same velocity with mass $m_1$ remaining fixed. Thus the right spring moves like a rigid body. Based on this argument and from \eqref{power_balance}, the viscous, dissipated power in the system at large fast times is $\eta {c_2}^2$, which is equal to the external power provided to the system (noting that even though mass $m_2$ moves for large times, it does so with uniform velocity in this problem resulting in no contribution to the rate of change of kinetic energy of the system).

\par
{\bf Savings in Computer time}

\begin{figure}[!h]
\centering
\includegraphics[scale = .4]{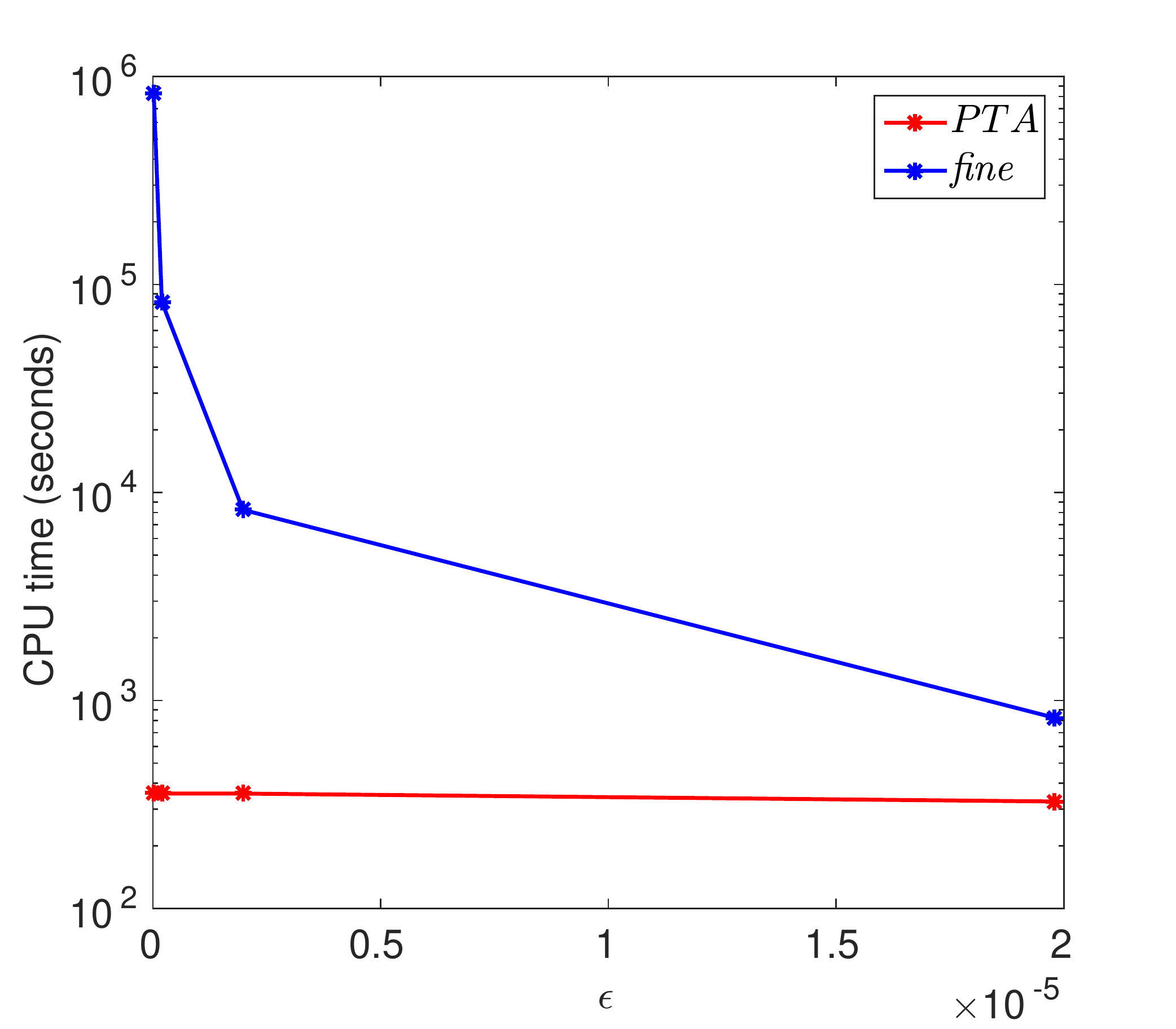}
\caption{\textit{Example II Case 1: Compute time comparison for simulations spanning $t=0.25$ to $t=0.5$}. }
\label{fig:ex2_cputime_uneql}
\end{figure}

Fig.\ref{fig:ex2_cputime_uneql} shows the comparison between the time taken by the fine and the PTA runs for simulations spanning $t=0.25$ to $t=0.5$ with $\Delta=0.05$. The speedup in compute time is given by the following polynomial: 
\begin{align}
S(\epsilon) = 2.28 \times 10^3 - 1.29 \times 10^{10} \, \epsilon + 6.46\times 10^{15} \, \epsilon^2 + 2.94 \times 10^{20} \, \epsilon^3. 
\end{align}

The function $S(\epsilon)$ is an interpolation of the computationally obtained data to a cubic polynomial. We used the datapoint of $\epsilon=1.98 \times 10^{-7}$ and obtained $S(1.98 \times 10^{-7})=231$. Results obtained in this case have accuracy of $0.0000\%$ error. As $\epsilon$ is decreased and approaches zero, the asymptotic value of $S$ is $2.28\times10^3$.  

\subsection{Results - Case 2: $\left(\frac{k_1}{m_1} =  \frac{k_2}{m_2} \right)$}\label{prb2:res_num_case2}

As already noted, in this case the quasi-static approach is not valid. The closed-form solution to this case is displayed in Appendix \ref{prb2:res_case2} (but it is not used in the computations). Recall that the computations are carried out when $L_1(w_1) = 0$ and $L_2(w_2) = c_2$. The PTA computations done in this section are with the load fixed while calculating $R^m_t$ and $R^m_{t-\Delta}$ using \eqref{eq:comp_impl_R1} and \eqref{eq:comp_impl_R2} respectively (by setting $\frac{dl_\epsilon}{d\sigma}=0$ in \eqref{eq:alg_fast_2}). The slow variable value  ($v(t+h)$ in \eqref{eq:v(t+h)}) is calculated using Simpson's rule as described in \textit{Remark} in Appendix \ref{prb2:res_case2}. All results in this section are non-dimensionalized following the discussion around \eqref{non-dim} in Appendix \ref{prb2:res_case2}. 

\begin{figure}[!h]
\centering
\includegraphics[scale = .35]{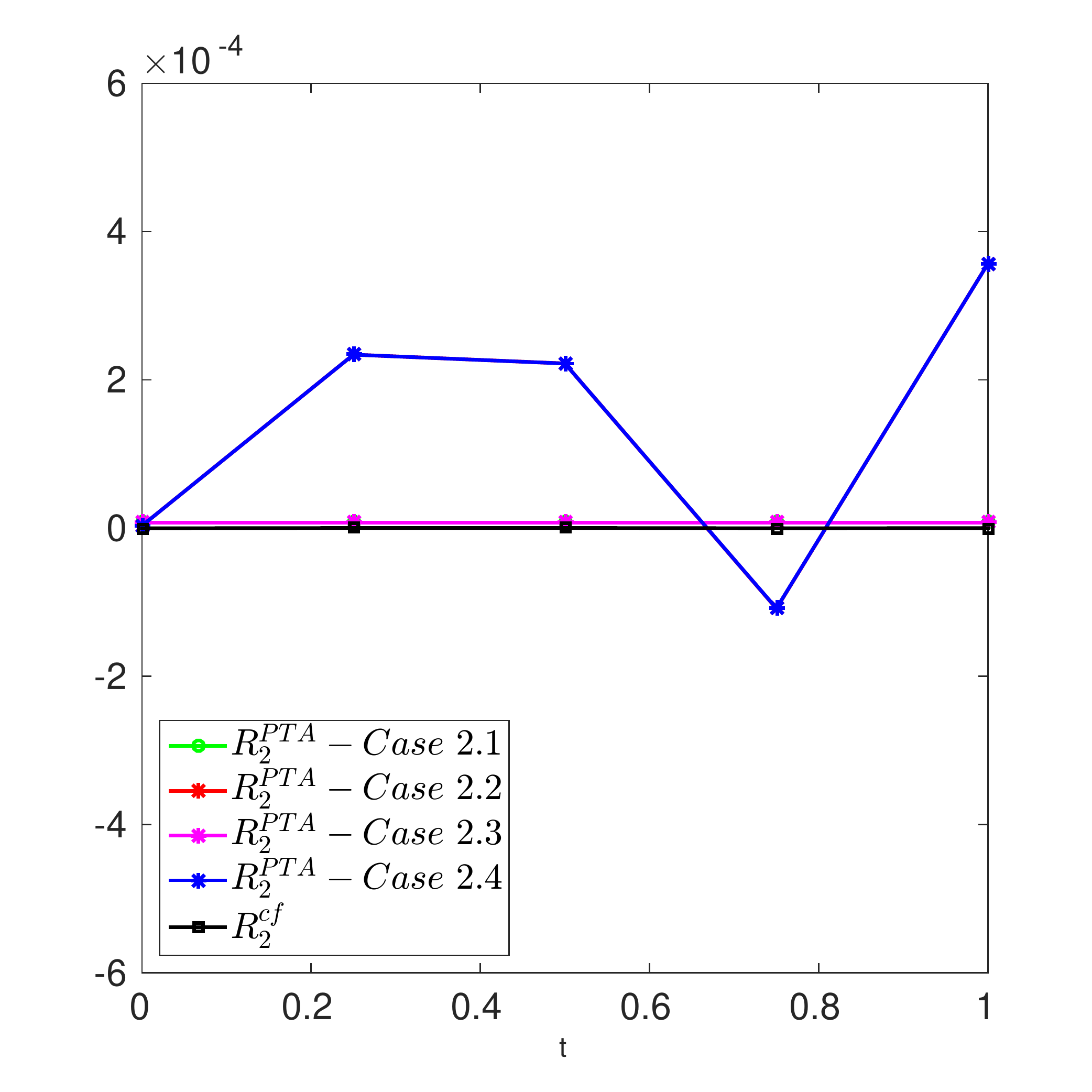}
\caption{\textit{Case 2.1 to 2.4: Comparison of $R_2^{PTA}$ and $R_2^{cf}$}. }
\label{fig:ex2_rf_pta}
\end{figure}

The following cases arise:  

\begin{itemize}
\item
{\bf Case 2.1.} When $c_2=0 $ and the initial condition does not have a component on the modes describing the dashpot being undeformed ($x_1=x_2$ and $y_1=y_2$), then the solution will go to rest. For example, the initial conditions ${x_1}^0=1.0$, ${x_2}^0=-0.5$ and ${y_1}^0 = {y_2}^0 = 0.0$ makes the solution (\ref{sol}) of Appendix \ref{prb2:res_case2} go to rest ($\kappa_3=\kappa_4=0$ in (\ref{eq:coeffs}) of Appendix \ref{prb2:res_case2}). The simulation results agree with the \textit{closed-form} results and go to zero.  

\item
{\bf Case 2.2.} When $c_2=0 $ and the initial condition has a component on the modes describing the dashpot being undeformed, then in the fast time limit the solution shows periodic oscillations whose energy is determined by the initial conditions. This happens, of course, for almost all initial conditions. One such initial condition is ${x_1}^0 = 0.5 $, ${x_2}^0 = -0.1 $ and ${y_1}^0={y_2}^0 = 0 $ ($\kappa_3=0$ but $\kappa_4=-0.4472$ in (\ref{eq:coeffs}) of Appendix \ref{prb2:res_case2}). The simulation results agree with the \textit{closed-form} results. 

This is in contrast with Case 1 where it is impossible to find initial conditions for which the solution shows periodic oscillations. 

\begin{figure}[!tbh]
\centering
\begin{minipage}{.45\textwidth}
  \centering
  \includegraphics[width=\linewidth]{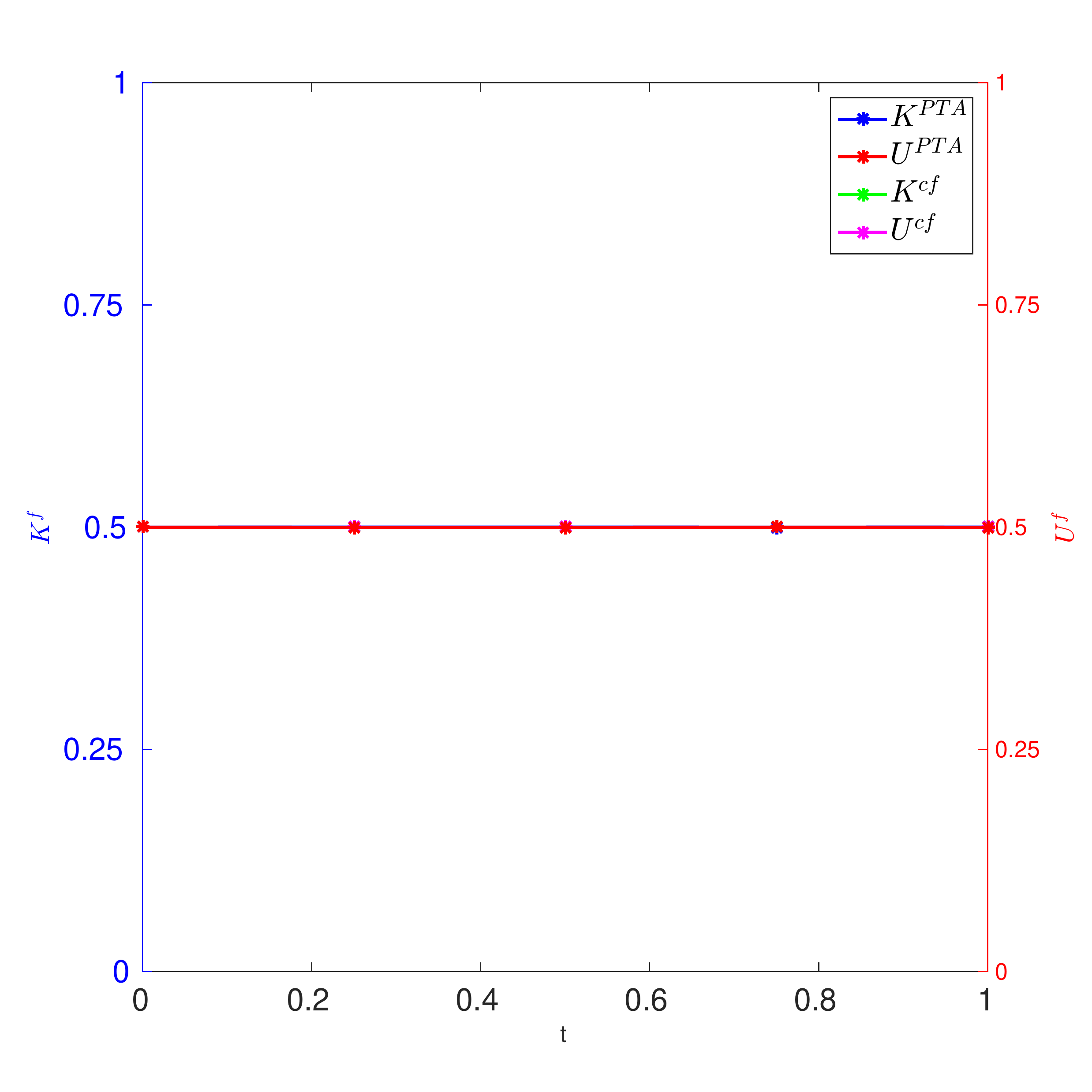}
  \caption{ \textit{Case 2.2 - Comparison of $K^{PTA}$, $U^{PTA}$, $K^{cf}$ and $U^{cf}$}. }
  \label{fig:ex2_pta_c2eql0}
\end{minipage}%
\hfill
\begin{minipage}{.45\textwidth}
  \centering
  \includegraphics[width=\linewidth]{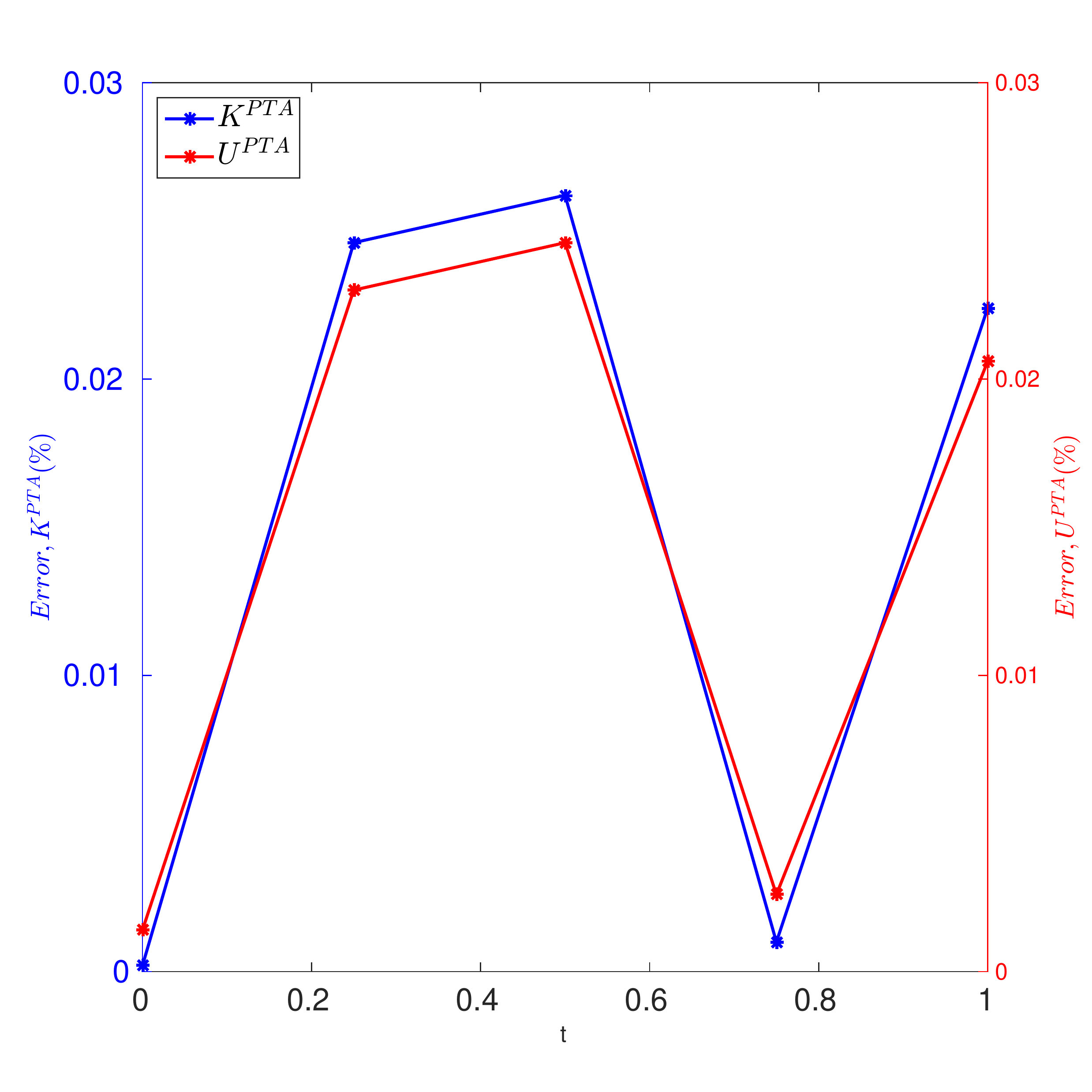}
  \caption{ \textit{Case 2.2 - Error in $K^{PTA}$ and $U^{PTA}$.}}
  \label{fig:ex2_err_c2eql0}
\end{minipage}
\end{figure}

In Fig. \ref{fig:ex2_pta_c2eql0}, we see that the PTA results are very close to the {\it closed-form} results. The error in PTA results are presented in Fig. \ref{fig:ex2_err_c2eql0}. 

\par
Oscillations persist in the limit and the potential and kinetic energies computed based on the Tikhonov framework as well as the quasi-static solution (\ref{ex2_uneql_quasistatic}) derived in Appendix \ref{prb2:res_case1} are not expected to, and do not, yield correct answers.

\item
{\bf Case 2.3.} When $c_2 \neq 0$ and the initial condition does not have a component on the modes describing the dashpot being undeformed, then the solution on the fast time scale for large values of $\sigma$ does not depend on the initial condition. One such initial condition is ${x_1}^0=1.0$, ${x_2}^0=-0.5$ and ${y_1}^0 = 0.0$ and $ {y_2}^0= 10^{-4} $. The \emph{closed-form average kinetic energy} ($K^{cf}$) and \emph{closed-form average potential energy} ($U^{cf}$) do not depend on the magnitude of the initial conditions in this case.

\item 
{\bf Case 2.4.} The initial condition has a component on the modes describing the dashpot being undeformed. But when $c_2 \neq 0 $, the dashpot gets deformed due to the translation of the mass $m_2$. The \emph{closed-form average kinetic energy} ($K^{cf}$) and \emph{closed-form average potential energy} ($U^{cf}$) depend on the initial conditions.

\begin{figure}[!tbh]
\centering
\begin{minipage}{.45\textwidth}
  \centering
  \includegraphics[width=\linewidth]{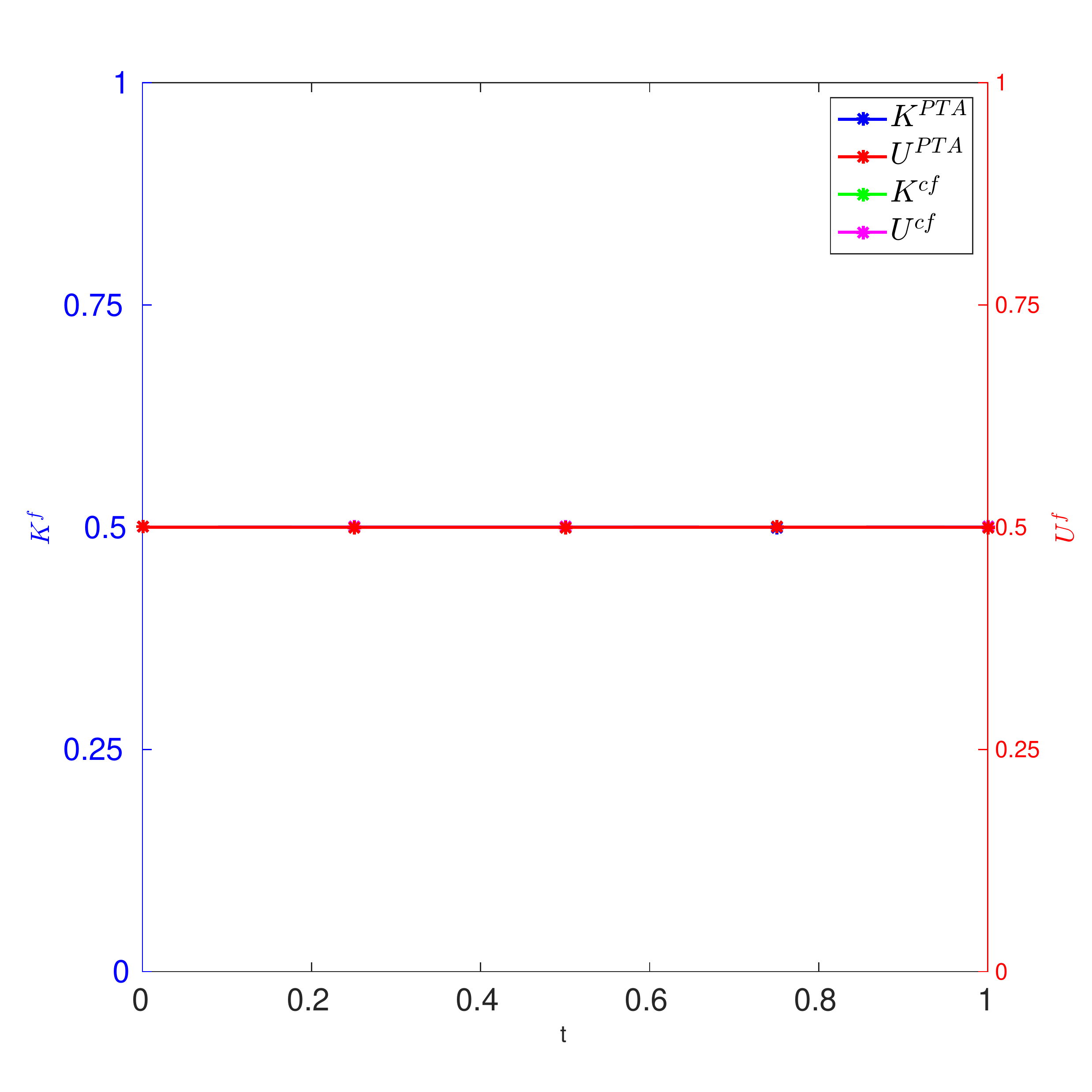}
  \caption{ \textit{Case 2.4 - Comparison of $K^{PTA}$, $U^{PTA}$, $K^{cf}$ and $U^{cf}$}. }
  \label{fig:ex2_pta}
\end{minipage}%
\hfill
\begin{minipage}{.45\textwidth}
  \centering
  \includegraphics[width=\linewidth]{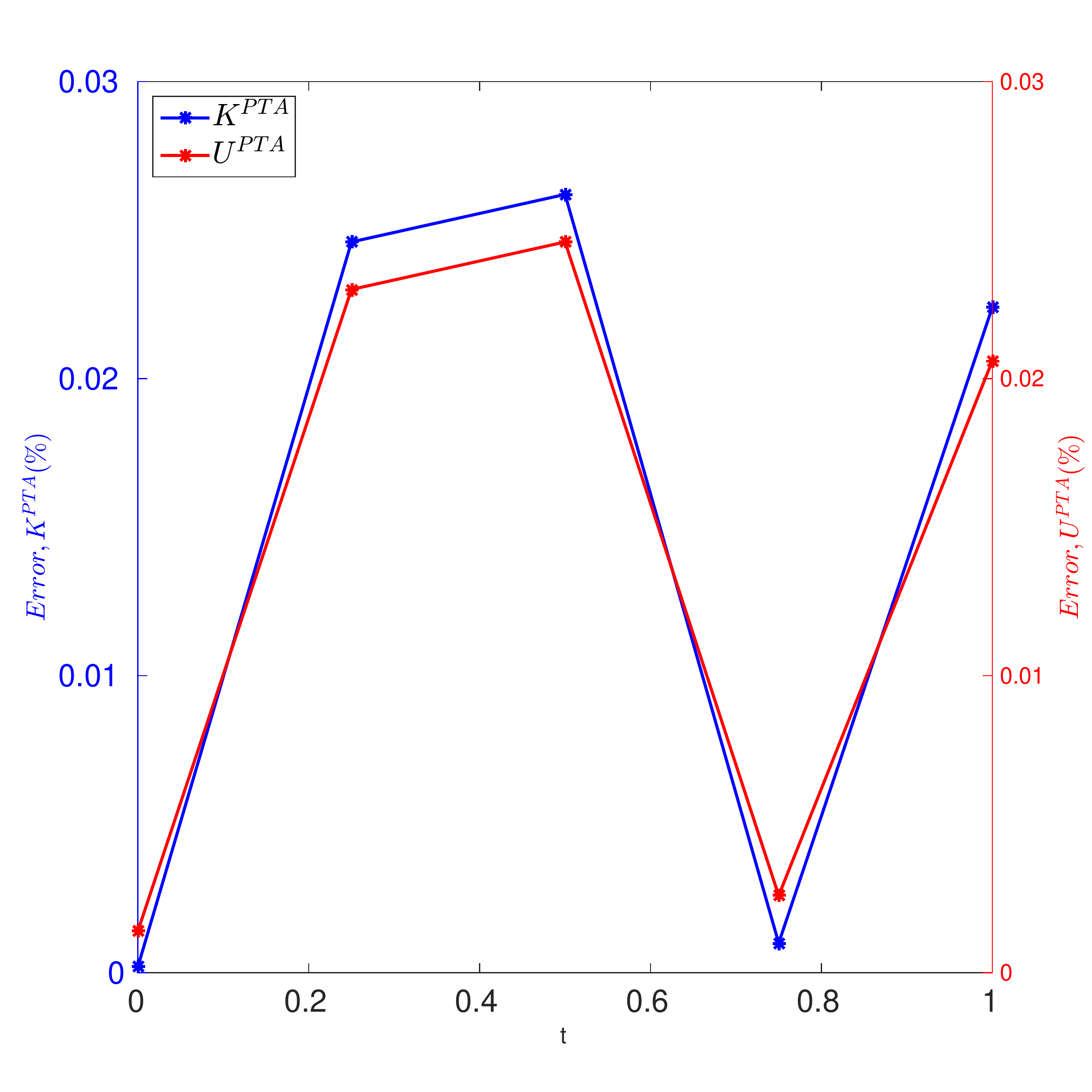}
  \caption{ \textit{Case 2.4 - Error in $K^{PTA}$ and $U^{PTA}$}. }
  \label{fig:ex2_err}
\end{minipage}
\end{figure}

In Fig. \ref{fig:ex2_pta}, we see that the PTA results are very close to the {\it closed-form} results. The errors in PTA results are presented in Fig. \ref{fig:ex2_err}. 

\item The comparison between $R_2^{PTA}$ and $R_2^{cf}$ for Case 2.1 to 2.4 is shown in Fig.\ref{fig:ex2_rf_pta}. The \textit{closed-form} result goes to zero and the PTA result becomes very small. The
displacement of the right wall ($w_2$) can be expressed as a function of time for Case 2.3 and Case 2.4 (given by \eqref{ex2:load} in Appendix \ref{prb2:res_case2}).
\item Again, oscillations persist in the limit, and the Tikhonov framework and the quasi-static approximation (see Appendix \ref{prb2:res_case1}) do not work in this case.

\item The results do not change if we decrease the value of $\epsilon$. However, the speedup changes as will be shown in  \textit{Savings in Computer time} later in this section.

\end{itemize} 

We used the same simulation parameters as in Case 1 ( shown in Table \ref{tab:simulation_details_uneql} ) but with $k_2=2 \times 10^7 N/m$ so that  $ \frac{k_1}{m_1} =  \frac{k_2}{m_2} $. Let us assume that the strain rate is $10^{-4} s^{-1} $. Then the slow time period, $T_s = {1 \over{\dot{\bar{\epsilon} } } }= 10000 s$. The fast time period is obtained as the period of fast oscillations of the spring, given by, $T_f = 2 \pi \sqrt{ \frac { m_1} { k_1} }= 0.002 s$. Thus, we find $\epsilon = \frac { T_f } {T_s} = 1.98 \times 10^{-7} $. While running the PTA code with $\epsilon = 0.002$, we have seen that the PTA scheme is not able to give accurate results and it breaks down.

\par
{\bf Power Balance} \\
The input power supplied to the system is ${1 \over T_s}R_2 \frac{d w_2}{dt} $. Since $R_2 = k_2(w_2 - x_2)$ and $\frac{d w_2}{dt}=T_s \, c_2$, the average value of input power at time $t$ is ${1 \over \Delta} \int_t^{t+\Delta} k_2 (w_2 - x_2) c_2 \, dt'$. From the results in \eqref{sol} of Appendix \ref{prb2:res_case2} and noting that the average of oscillatory terms over time $\Delta$ is approximately 0, we see that ${1 \over \Delta} \int_t^{t+\Delta}{(w_2-x_2)}\,dt'={1 \over \Delta} \int_t^{t+\Delta} \{ c_2 T_s t' - (c_2 T_s t' - \frac {\eta c_2}{k_2}) \} \, dt' = \frac {\eta c_2}{k_2}$. Hence average input power supplied is $\eta c_2^2$. The average dissipation at time $t$ is $ {1 \over \Delta} \int_t^{t+\Delta} \frac {\eta}{T_s^2} {(\frac{d x_2}{dt'}-\frac{d x_1}{dt'}) }^2 \, dt' $. Using the result from \eqref{sol} of Appendix \ref{prb2:res_case2}, and using the same argument ~that the average of oscillatory terms over time $\Delta$ is approximately 0, we can say that the dissipation is $\eta c_2^2$. Thus, the average input power supplied to the system is equal to the dissipation in the damper. A part of the input power also goes into the translation of the mass $m_2$. But its value is very small compared to the total kinetic energy of the system. 

\par
{\bf Savings in Computer time}\\ It takes the PTA run around 62 seconds to compute the calculations that start at slow time which is a multiple of $h$ (steps 1 through 5 in section \ref{impl_algo} ). It takes the fine theory run around 8314 seconds to evolve the fine equation starting at slow time $nh$ to slow time $(n+1)h$, where $n$ is a positive integer. Thus, we could achieve a speedup of 134. We expect that the speedup will increase if we decrease the value of $\epsilon$.
\begin{figure}[!h]
\centering
\includegraphics[scale = .4]{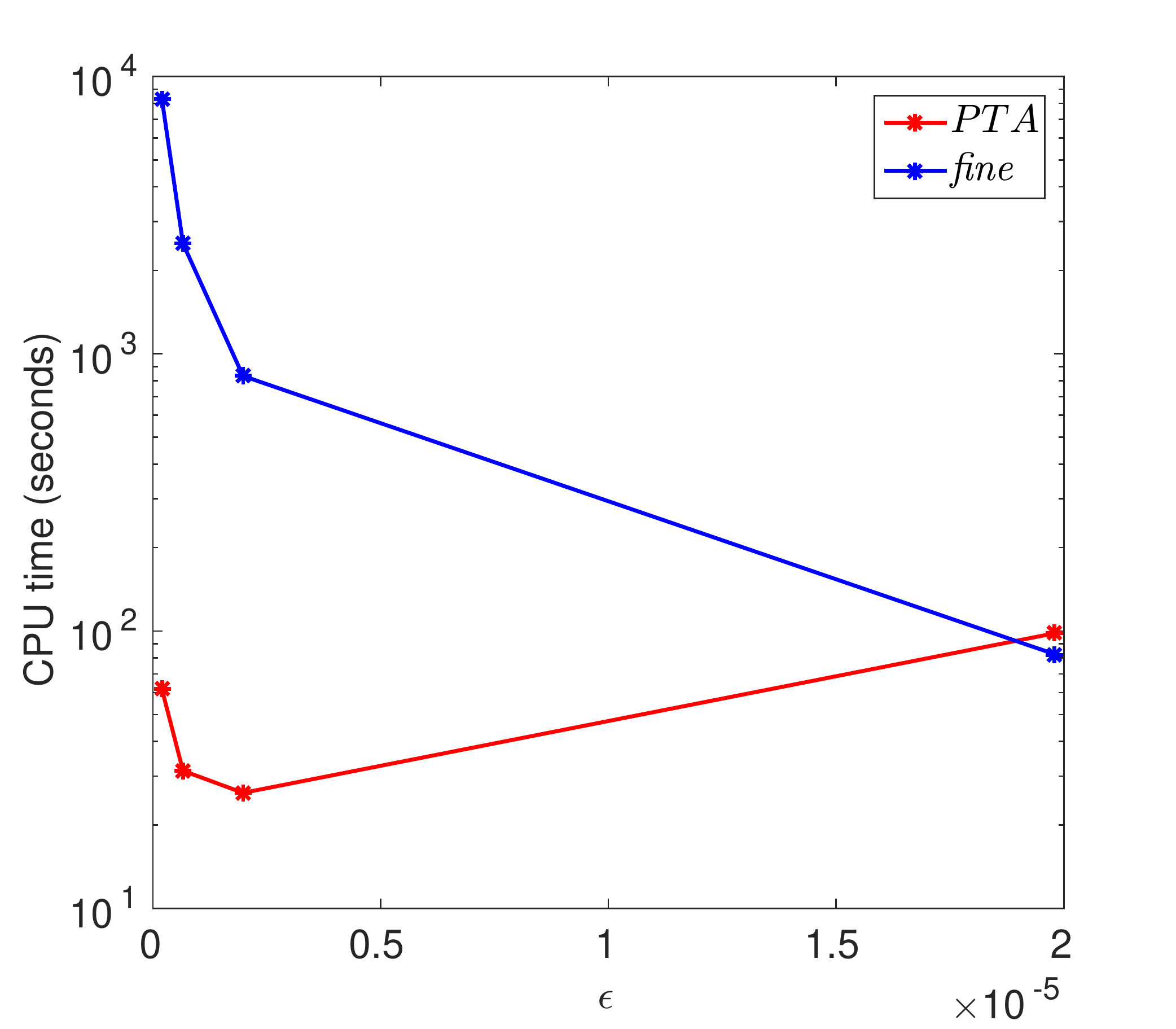}
\caption{\textit{Example II Case 2: Compute time comparison for simulations spanning $t=0.25$ to $t=0.5$}. }
\label{fig:ex2_cputime_eql}
\end{figure}

Fig.\ref{fig:ex2_cputime_eql} shows the comparison between the time taken by the fine and the PTA runs for simulations spanning $t=0.25$ to $t=0.5$ with $\Delta=0.05$. 
\par
The speedup in compute time is given by the following polynomial: 
\begin{align}
S(\epsilon) = 164.34 - 1.62\times 10^{8} \, \epsilon + 5.23 \times 10^{13} \, \epsilon^2 -2.25 \times 10^{18} \epsilon^3.
\end{align}
The function $S(\epsilon)$ is an interpolation of the computationally obtained data to a cubic polynomial. We used the datapoint of $\epsilon=1.98 \times 10^{-7}$ and obtained $S(1.98 \times 10^{-7})=134$. This speedup corresponds to an accuracy of $0.026\%$ error. As $\epsilon$ is decreased and approaches zero, the asymptotic value of $S$ is 164.  

\section{Example III: Relaxation oscillations of oscillators}
\noindent
This is a variation of the classical relaxation oscillation example(see, e.g., \cite{Artstein_2002_perturbed}). Consider the four-dimensional system
\begin{align}\label{eq:11.1}
 \frac{dx}{dt}&= z\nonumber\\
 \frac{dy}{dt}&= \frac{1}{\epsilon}(-x+y-y^{3})\\
 \frac{dz}{dt}&=
               \frac{1}{\epsilon}(w+(z - y)(\frac{1}{8}-w^{2}-(z-y)^{2}))\nonumber\\
 \frac{dw}{dt}&=
               \frac{1}{\epsilon}(-(z-y)+w(\frac{1}{8}-w^{2}-(z-y)^{2}))\nonumber.
\end{align}
\par
Notice that the $(z,w)$ coordinates oscillate around the point $(y,0)$ (in the $(z,w)$-space), with oscillations that converge to a circular limit cycle of radius ${1 \over \sqrt{8} }$. The coordinates $(x,y)$ follow the classical relaxation oscillations pattern (for the fun of it, we replaced $y$ in the slow equation by $z$, whose average in the limit is $y$). In particular, the limit dynamics of the $y$-coordinate moves slowly along the stable branches of the curve $0 = -x + y - y^{3}$,
with discontinuities at $x = -\frac{2}{3 \sqrt{3}}$ and $x = \frac{2}{3 \sqrt{3}}$. In turn, these discontinuities carry with them discontinuities of the oscillations in the $(z,w)$ coordinates. The goal of the computation is to
follow the limit behavior, including the discontinuities of the oscillations.
\par
\subsection{Discussion} The slow dynamics, or the load, in the example is the $x$-variable. Its value does not determine the limit invariant measure in the fast dynamics, which comprises a point $y$ and a limit circle in the $(z,w)$-coordinates. A slow observable that will determine the limit invariant measure is the $y$-coordinate. In particular, conditions (\ref{eq:8.1}) an (\ref{eq:8.2}) hold except at points of discontinuity, and so does the conclusion of Theorem 8.2. Notice, however, that this observable does go through periodic discontinuities.

\par
\subsection{Results} 
\begin{figure}[!h]
\centering
\begin{minipage}{.45\textwidth}
  \centering
  \includegraphics[width=\linewidth]{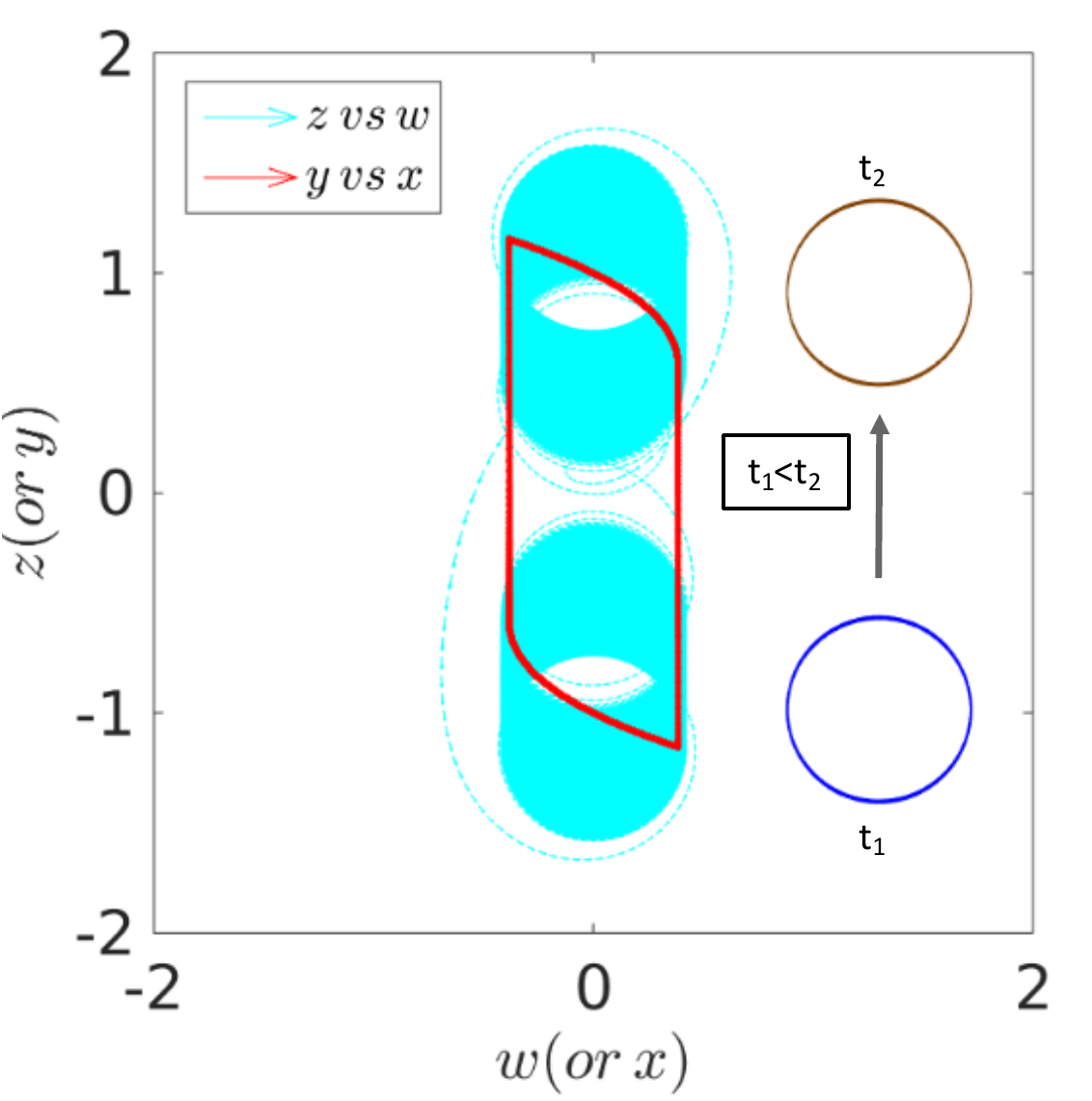}
  \caption{{\it Trajectory of \eqref{eq:11.1}. The vertical branches of the $y~vs~x$ curve correspond to very fast move on the fast time scale. The blue curve shows the portion of the phase portrait of the z vs w trajectory obtained around time $t_1$ while the brown curve shows the portion around a later time $t_2$. \textup{(}For interpretation of the references to color in this figure legend, the reader is referred to the web version of this article.\textup{)}}} 
  \label{fig:ex3_fast}
\end{minipage}%
\hfill
\begin{minipage}{.5\textwidth}
  \centering
  \includegraphics[width=\linewidth]{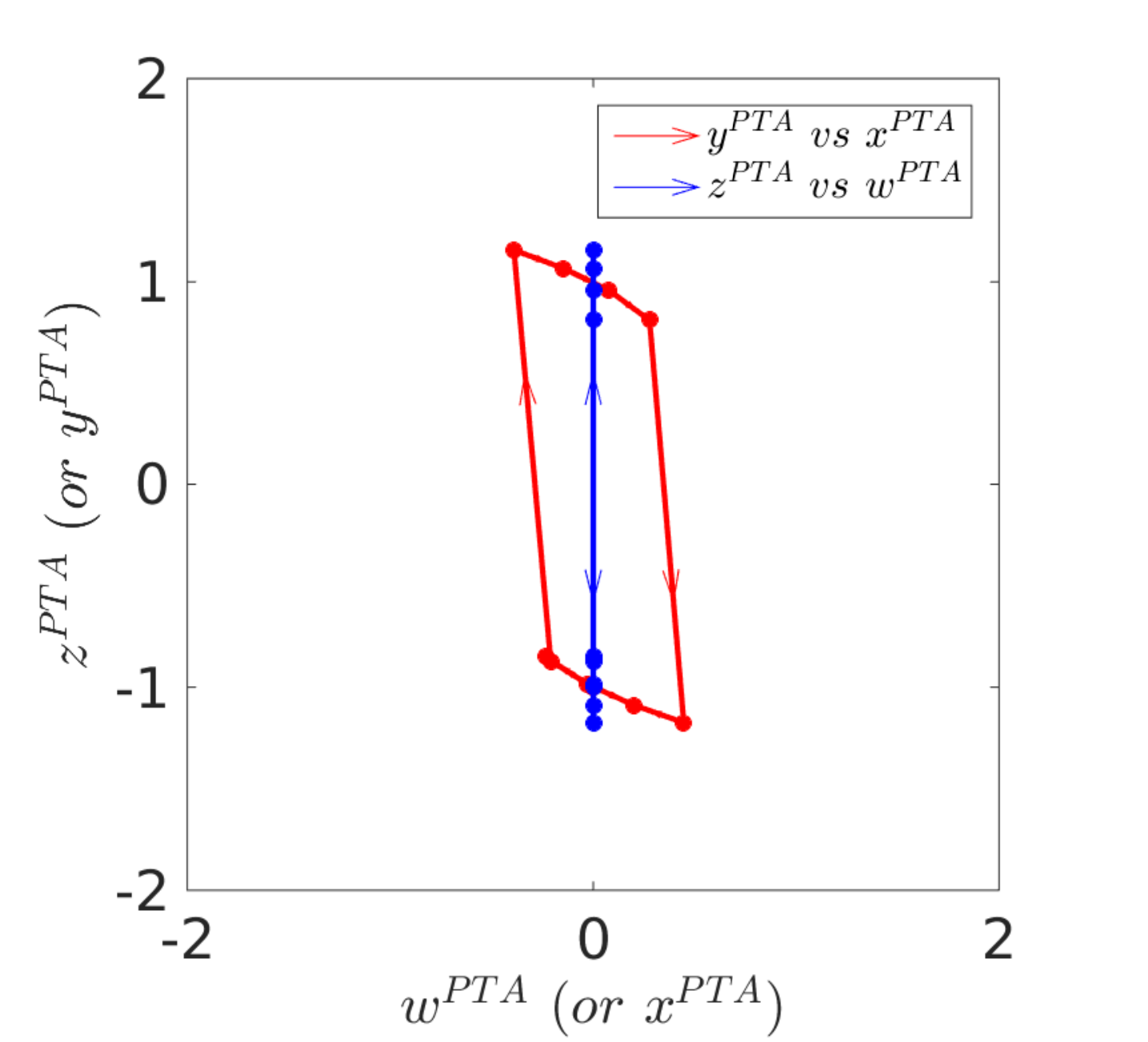}
  \caption{ {\it PTA  result. The portion with the arrows correspond to very rapid evolution on the slow time scale}.} 
  \label{fig:ex3_pta}
\end{minipage}
\end{figure}

We see in Fig. \ref{fig:ex3_fast} that the $y$-coordinate moves slowly along the stable branches of the curve $0 = -x + y - y^{3}$ which is evident from the high density of points in these branches of the curve as can be seen in Fig. \ref{fig:ex3_fast}. There are also two discontinuities at $x = -\frac{2}{3 \sqrt{3}}$ and $x = \frac{2}{3 \sqrt{3}}$.  The pair $(z,w)$ oscillates around $(y,0)$ in circular limit cycle of radius ${ 1 \over \sqrt{8} }$.  

\par 
In Fig. \ref{fig:ex3_pta}, we see that the average of $z$ and $y$ which are given by the $y$-coordinate in the plot, are the same which acts as a verification that our scheme works correctly. Also, average of $w$ is 0 as expected. 
\par 
Since there is a jump in the evolution of the measure at the discontinuities (of the Young measure), the observable value obtained using extrapolation rule is not able to follow this jump. However, the observable values obtained using the guess for fine initial conditions at the next jump could follow the discontinuity. This is the principal computational demonstration of this example. 

%\begin{figure}[!h]
%\centering
%\includegraphics[scale = .5]{figures/ex3_closest.pdf}
%\caption{ \textit{Curve 1 initiates at point 1 and is able to follow the discontinuity in the measure which lies on the opposite stable branch of the $y \, vs \, x$ curve. Similarly, curve 2 initiates at point 2 which is in the opposite stable branch and is able to follow the discontinuity in the measure} }
%\label{fig:rel_err}
%\end{figure}

\begin{figure}[h!]
\centering
\includegraphics[scale = .4]{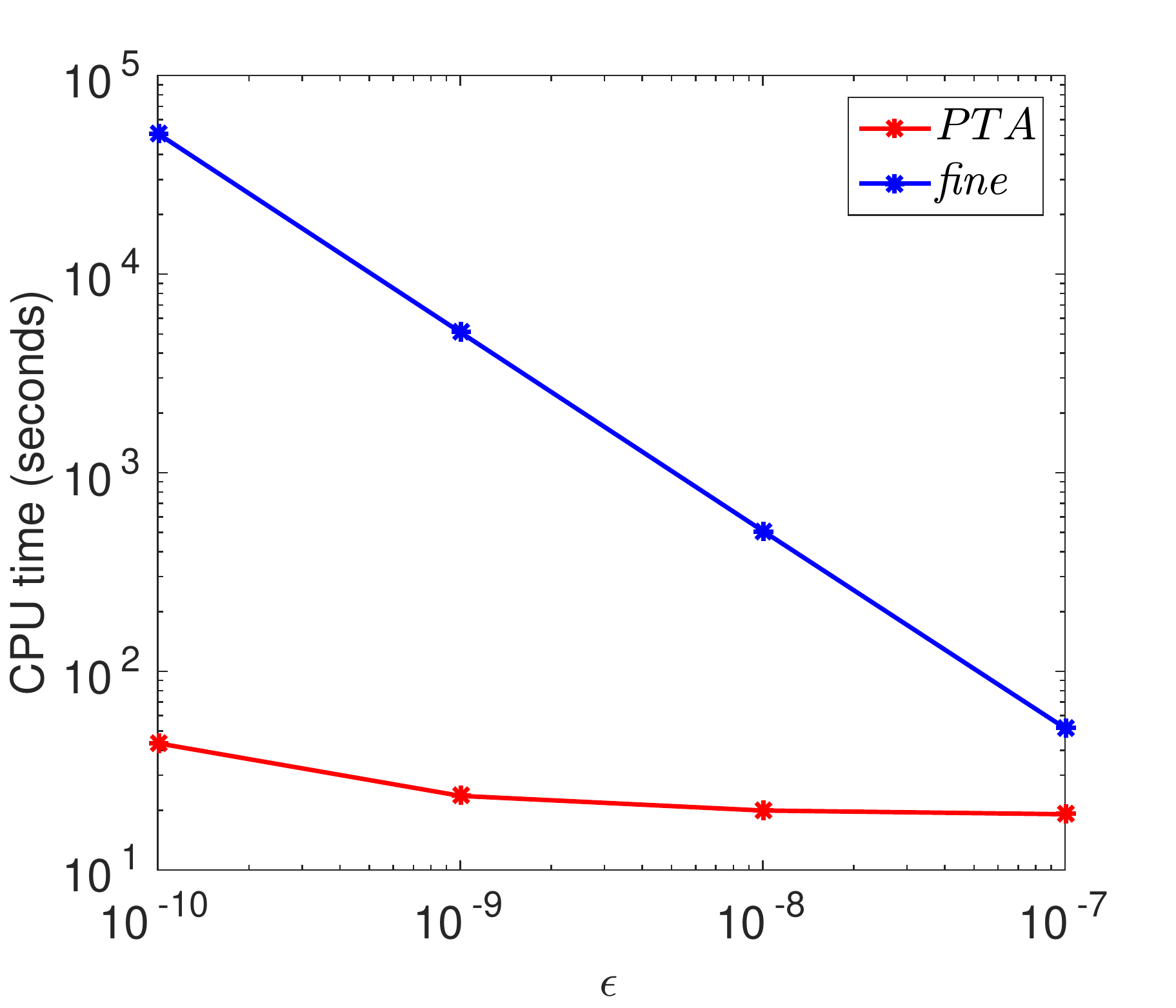}
\caption{\textit{Example III - Compute time comparison for simulations spanning $t=0.2$ to $t=0.4$}. }
\label{fig:ex3_cputime}
\end{figure}

\begin{figure*}[!h]
    \centering
    \begin{subfigure}[h]{0.5\textwidth}
        \centering
        \includegraphics[height=3.1in]{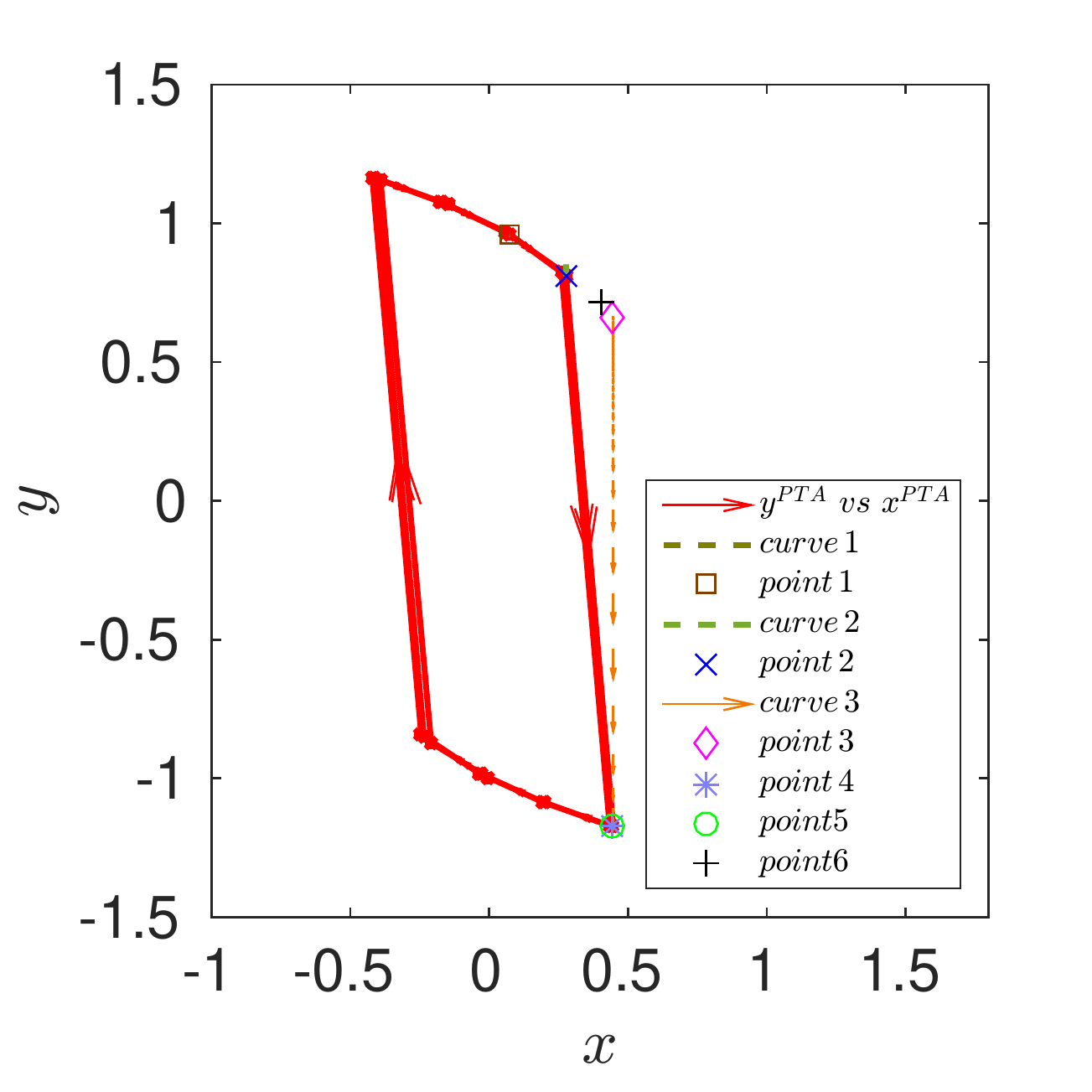}
        \caption{}
    \end{subfigure}%
    ~ 
    \begin{subfigure}[h]{0.5\textwidth}
        \centering
        \includegraphics[height=3.1in]{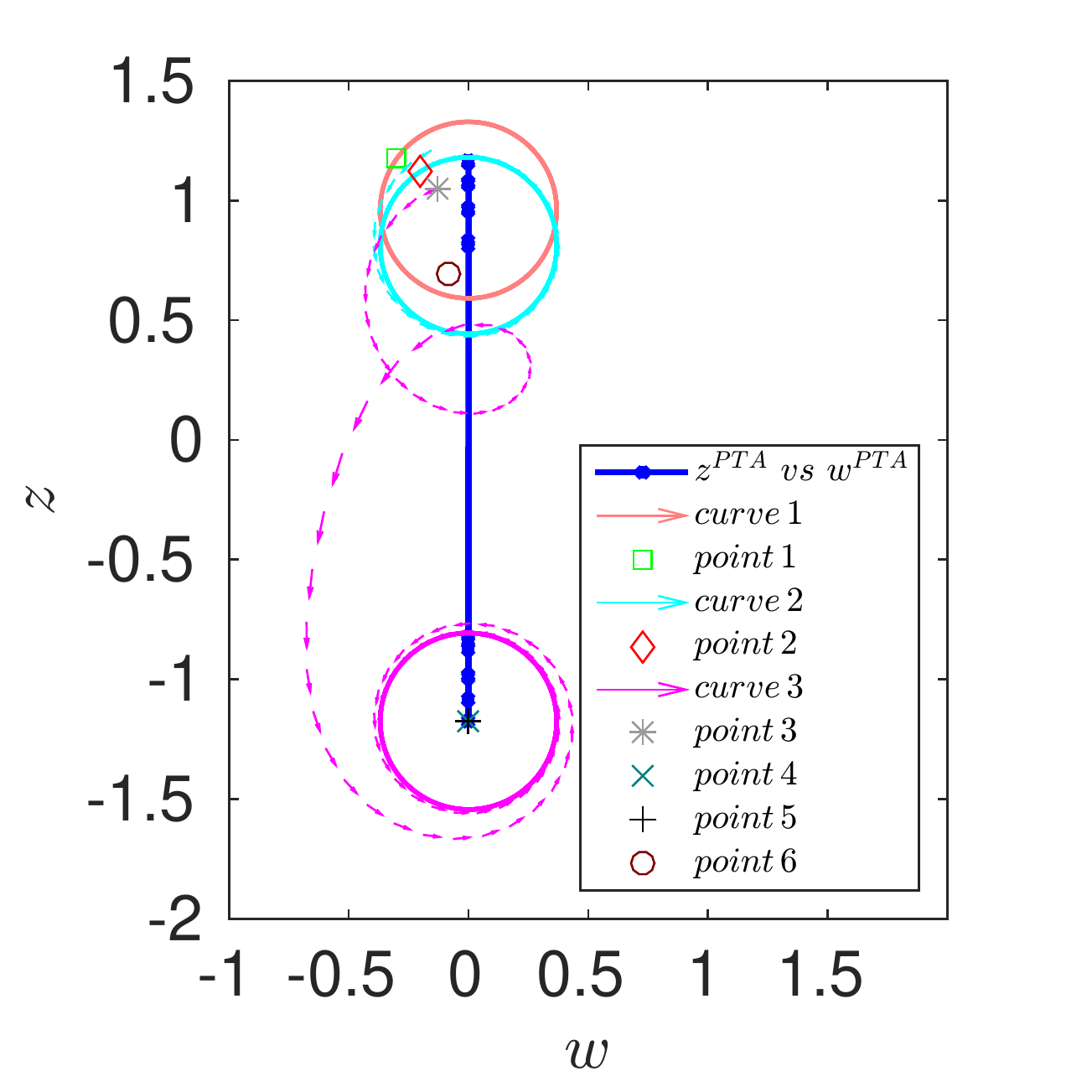}
        \caption{}
        \end{subfigure}
        \caption{\textit{ This figure shows how PTA scheme predicts the correct values of slow observables when there is a discontinuity in the Young measure. Part \textup{(}a\textup{)} shows the details for $y \, vs \, x$. The denotations of the different points mentioned here are provided in Step 3 of Section \ref{impl_algo}. Curve 1 is the set of all points in support of the measure at time $t-h-\Delta$. The point $x^{cp}_{t-h-\Delta}$ is given by point 1 in the figure (we obtain Curve 1 and point 1 using the details mentioned in Step 3 of Section \ref{impl_algo}). Curve 1 reduces to a point near point 1, so it is not visible in the figure. Curve 2 is the set of all points in support of the measure at time $t-\Delta$. The point $x^{arb}_{t-\Delta}$ is given by point 2 (we obtain Curve 2 and point 2 using the details mentioned in Step 3 of Section \ref{impl_algo}). Curve 2 reduces to a point very close to point 2, so it is not visible. Point 3 is the initial guess for time $t+h-\Delta$ which we calculate using the details in Step 3 of Section \ref{impl_algo}. Curve 3 is the set of all points in support of the measure at time $t+h-\Delta$. The slow observable value at time $t+h-\Delta$ obtained from the fine run is point 4. Point 5 is the slow observable value obtained from the PTA run using the details in Step 4 of Section \ref{impl_algo}. Point 6 corresponds to slow observable values obtained solely by using the coarse evolution equation without using the initial guess at time $t+h-\Delta$ (using Step 2 of Section \ref{impl_algo}). Part \textup{(}b\textup{)} shows the corresponding details for $ z \, vs \, w $. In this figure, we see that Curve 1 and Curve 2 do not reduce to a point.} }  
\label{fig:pta_working}
\end{figure*}

Fig. \ref{fig:pta_working} shows the working of the PTA scheme when there is a discontinuity in the Young measure. We obtain the initial guess at time $t+h-\Delta$ by extrapolating the closest point projection at time $t-h-\Delta$ of a point on the measure at time $t-\Delta$. The details of the procedure are mentioned in Step 3 of section \ref{impl_algo}. When there is a discontinuity in the Young measure, the results of the slow observables obtained using coarse evolution (Point 6) is unable to follow the discontinuity. But when the fine run is initiated at the initial guess at time $t+h-\Delta$ which is given by Point 3 in the figure, the PTA scheme is able to follow the jump in the measure and we obtain the correct slow observable values (Point 5) which is very close to the slow observable value obtained from the fine run (Point 4).

\par 
\textbf{Savings in computer time.} 
Fig.\ref{fig:ex3_cputime} shows the comparison between the time taken by the fine and the PTA runs for simulations
spanning $t=0.2$ to $t = 0.4$ with $\Delta
= 0.01$. The speedup in compute time as a function of $\epsilon$ for $\Delta=0.01$, is given by the following polynomial: 
\begin{align}
S (\epsilon) = 1.29 \times 10^3 - 1.19 \times 10^{12} \, \epsilon + 1.17 \times 10^{20} \, \epsilon ^2 - 1.05 \times 10^{27} \, \epsilon^3.
\end{align}  
The function $S (\epsilon)$ is an interpolation of the computationally obtained data to a cubic polynomial. A more efficient calculation yielding higher speedup is to calculate the slow variable $v$ using Simpson's rule instead of using \eqref{eq:v(t+h)}, by employing the procedures outlined in section \ref{prb2:res_num_case1}, section \ref{prb2:res_num_case2} and the associated \textit{Remark} in Appendix \ref{prb2:res_case2}. We used the datapoint of $\epsilon=10^{-10}$ for this problem and obtained $S(10^{-10})=1.17 \times 10^3$. This speedup corresponds to an accuracy of $4.33\%$ error. As we further decrease $\epsilon$ and it approaches zero, the asymptotic value of $S$ becomes $1.29 \times 10^3$. 

\par 
{\bf Remark.} As a practical matter, it seems advantageous to set $\epsilon = 0$ in \eqref{eq:alg_fast_2} for the computations of $R^m_{t-\Delta}$ in \eqref{eq:comp_impl_R2} and $R^m_t$ in \eqref{eq:comp_impl_R1}. Related to this, when the value of $\epsilon$ is decreased, calculating the slow variable value  ($v(t+h)$ in \eqref{eq:v(t+h)}) using Simpson's rule as described in \textit{Remark} in Appendix \ref{prb2:res_case2} reduces $T^{cpu}_{PTA}$ considerably and improves the speedup $S(\epsilon)$.

\section{Concluding remarks}
The focus of this paper has been the precise definition and demonstration of a computational tool to probe slow time-scale behavior of rapidly evolving microscopic dynamics, whether oscillatory or exponentially decaying to a manifold of slow variables, or containing both behaviors. A prime novelty of our approach is in the introduction of a general family of observables ($H$-observables) that is universally available, and a practical computational scheme that covers cases where the invariant measures may not be uniquely determined by the slow variables in play, and one that allows the tracking of slow dynamics even at points of discontinuity of the Young measure. We have solved three model problems that nevertheless contain most of the complications of averaging complex multiscale temporal dynamics. It can be hoped that the developed tool is of substantial generality for attacking real-world practical problems related to understanding and engineering complex microscopic dynamics in relatively simpler terms.

\clearpage
\begin{appendices}
\section{Verlet Integration }\label{app:verlet}

The implementation of Verlet scheme that we used in Example II to integrate the fine equation (\ref{eq:10.4}) is as follows \cite{Sandvik_2016_py502}: 
\begin{align}\label{eq:verlet}
x_1(\sigma + \Delta \sigma) &= x_1(\sigma) + \Delta \sigma \, y_1(\sigma) + \frac {1} {2} {\Delta \sigma}^2 \,  \frac {d y_1}{d \sigma} \nonumber\\
x_2(\sigma + \Delta \sigma) &= x_2(\sigma) + \Delta \sigma \, y_2(\sigma) + \frac {1} {2} {\Delta \sigma}^2 \, \frac {d y_2}{d \sigma} \nonumber\\
\hat{y_1}(\sigma + \Delta \sigma) &= y_1(\sigma) + \frac {1}{2} \Delta \sigma \left( \frac {d y_1}{d \sigma}  + a_1 \left\{x_1(\sigma + \Delta \sigma), y_1(\sigma) + \Delta \sigma \frac {d y_1}{d \sigma}, y_2(\sigma) + \Delta \sigma \frac {d y_2}{d \sigma}, w_1(\sigma) \right\} \right) \nonumber\\
\hat{y_2}(\sigma + \Delta \sigma) &= y_2(\sigma) + \frac {1}{2} \Delta \sigma \left( \frac {d y_2}{d \sigma}  + a_2 \left\{x_2(\sigma + \Delta \sigma), y_1(\sigma) + \Delta \sigma \frac {d y_1}{d \sigma}, y_2(\sigma) + \Delta \sigma \frac {d y_2}{d \sigma}, w_2(\sigma) \right\} \right) \nonumber\\
y_1(\sigma + \Delta \sigma) &= y_1(\sigma) + \frac{1}{2} \Delta \sigma \left( \frac{d y_1}{d \sigma} + a_1\left\{x_1(\sigma+\Delta \sigma), \hat{y_1} (\sigma + \Delta \sigma), \hat{y_2}(\sigma + \Delta \sigma), w_1(\sigma) \right\} \right) \nonumber\\
y_2(\sigma + \Delta \sigma) &= y_2(\sigma) + \frac{1}{2} \Delta \sigma \left( \frac{d y_2}{d \sigma} + a_2 \left\{x_2(\sigma+\Delta \sigma), \hat{y_1} (\sigma + \Delta \sigma), \hat{y_2}(\sigma + \Delta \sigma), w_2(\sigma)\right\} \right).  
\end{align}
Here $a_1(x_1,y_1,y_2,w_1)=\frac{d y_1}{d\sigma}$ and $a_2(x_2,y_1,y_2,w_2)=\frac{d y_2}{d\sigma}$ where $\frac{d y_1}{d\sigma}$ and $\frac{d y_2}{d\sigma}$ are given by (\ref{eq:10.4}). 

\section{Example II: Derivation of system of equations}\label{prb2:equations}

\noindent
Two massless springs and masses $m_{1}$ and $m_{2}$ are connected through a dashpot damper, and each is attached to two bars, or walls, that may move very slowly, compared to possible oscillations of the springs. The system is described in Fig. \ref{fig:prb2} of Section \ref{vibrate}.
\par
Let the springs constants be $k_{1}$ and $k_{2}$ respectively, and let $\eta$ be the dashpot constant. Denote by $x_{1}$ and $x_{2}$, and by $w_{1}$ and $w_{2}$, the displacements from equilibrium positions of the masses of the springs and the positions of the two walls. We agree here that all positive displacements are toward the right. We think of the movement of the two walls as being external to the system, determined by a ``slow" differential equation. The movement of the springs, however, will be ``fast", which we model as singularly perturbed. Let $m_{w_1}$ and $m_{w_2}$ be the masses of the left and the right walls respectively.  The displacements $x_1$, $x_2$, $w_1$ and $w_2$ have physical dimensions of length. The spring constants $k_1$ and $k_2$ have physical dimensions of force per unit length while $m_1$ and $m_2$ have physical dimensions of mass.
In view of the assumptions just made, a general form of the dynamics of the system is given by following set of equations:
\begin{align}\label{eq:10.1}
m_1 \frac{d^2 x_1}{d {t^*}^2} &= -k_1 (x_1 - w_1) + \eta \left( \frac{d x_2}{d t^*} - \frac{d x_1}{d t^*} \right) \nonumber\\
m_2 \frac{d^2 x_2}{d {t^*}^2} &= -k_2 (x_2 - w_2) - \eta \left( \frac{d x_2}{d t^*} - \frac{d x_1}{d t^*} \right) \\
m_{w_1} \frac{d^2 w_1}{d {t^*}^2} &= k_1 (x_1 - w_1) + R_1 \nonumber\\
m_{w_2} \frac{d^2 w_2}{d {t^*}^2} &= k_2 (x_2 - w_2) + R_2, \nonumber
\end{align}
where $R_1$ and $R_2$ incorporate the reaction forces on the left and right walls respectively due to their prescribed motion. We agree that forces acting toward the right are being considered positive. We make the assumption that $m_{w_1}=m_{w_2}=0$. The time scale $t^*$ is a time scale with physical dimensions of time. In our calculations, however, we address a simplified version, of first order equations, that can be obtained from the previous set by appropriately specifying what the forces on the system are:
\begin{align}\label{eq:10.2}
 \frac{dx_{1}}{dt^*} &= y_{1} \nonumber\\
 \frac{dy_{1}}{dt^*}&= -\frac{k_{1}}{m_{1}}(x_{1}-w_{1}) + \frac{\eta}{m_{1}}(y_{2}-y_{1})\nonumber\\
 \frac{dx_{2}}{dt^*}&= y_{2}\\
 \frac{dy_{2}}{dt^*}&= -\frac{k_{2}}{m_{2}}(x_{2}-w_{2})-\frac{\eta}{m_{2}}(y_{2}-y_{1})\nonumber\\
 \frac{dw_{1}}{dt^*}&=L_{1}(w_{1})\nonumber\\
 \frac{dw_{2}}{dt^*}&=L_{2}(w_{2})\nonumber.
\end{align}

The motion of the walls are determined by the functions $L_1(w_1)$ and $L_2(w_2)$. In the derivations that follow, we use the form $L_1(w_1) = c_1$ and $L_2(w_2) = c_2$, with $c_1=0$ and $c_2$ being a constant. The terms $c_1$ and $c_2$ have physical dimensions of velocity. 

\par
We define a coarse time period, $T_s$, in terms of the applied loading rate as $T_s=\frac{const}{L_2}$. Hence, $T_s c_2$ is a constant which is independent of the value of $T_s$. The fine time period, $T_f$, is defined as the smaller of the periods of the two spring mass systems. We then define the \emph{non-dimensional} slow and fast time scales as $t = \frac{t^*}{T_s}$ and $\sigma = \frac{t^*}{T_f}$, respectively. The parameter $\epsilon$ is given by $\epsilon = \frac{T_f}{T_s}$. Then the dynamics on the slow time-scale is given by (\ref{eq:10.3}) in Section \ref{vibrate}. The dynamics on the fast time-scale is written as:
%%

%\noindent\begin{minipage}{.5\linewidth}
%\vspace
%Slow time-scale: 
%\end{minipage}%

%\begin{minipage}{.5\linewidth}
%Fast time-scale:
\begin{align}\label{eq:10.4}
 \frac{dx_{1}}{d\sigma} &= T_f \, y_{1}\nonumber\\
 \frac{dy_{1}}{d\sigma}&= -T_f \left(\frac{k_{1}}{m_{1}}(x_{1}-w_{1}) - \frac{\eta}{m_{1}}(y_{2}-y_{1})\right)\nonumber\\
 \frac{dx_{2}}{d\sigma}&= T_f \, y_{2} \nonumber\\
 \frac{dy_{2}}{d\sigma}&= -T_f \left(\frac{k_{2}}{m_{2}}(x_{2}-w_{2}) + \frac{\eta}{m_{2}}(y_{2}-y_{1})\right)\nonumber\\
 \frac{dw_{1}}{d\sigma}&= \epsilon \, T_s \, L_{1}(w_{1})\nonumber\\
 \frac{dw_{2}}{d\sigma}&=\epsilon \, T_s \, L_{2}(w_{2}).
\end{align}

 {\bf Remark.} A special case of (\ref{eq:10.3}) is when $c_1=0$ and $c_2=0$. This represents the unforced system i.e. the walls remain fixed. Then (\ref{eq:10.3}) is modified to: 
 \begin{align}\label{eq:unforced system}
\dot{ {\bf x} } = {\bf B} {\bf x}, 
 \end{align}
 where ${\bf x}={(x_1, y_1, x_2, y_2)}^T$. The overhead dot represent time derivatives w.r.t. $t$. The matrix {\bf B} is given by
 \[
{\bf B} = T_s \begin{pmatrix}  0 & 1 & 0 & 0 \\
								 -\frac{k_1}{m_1} & -\frac{\eta}{m_1} & 0 & \frac{\eta}{m_1} \\
								0 & 0 & 0 & 1 \\
								  0 & \frac{\eta}{m_2} & -\frac{k_2}{m_2} & -\frac{\eta}{m_2} \end{pmatrix}  
 .\]

\section{Example II: Case 1 - Validity of commonly used approximations}\label{prb2:res_case1} 
The mechanical system \eqref{eq:10.3} can be written in the form
%\begin{align}\label{appC_qs}
%m \tilde{{\bf A}}_1 \frac{d^2}{d {t^*}^2} \left\{\frac{ {\bf w} }{ {\bf x} }\right\} 
%+ D \tilde{{\bf A}}_2 \frac{d}{d t^*} \left\{\frac{ {\bf w} }{ {\bf x} }\right\}
%+ k \tilde{{\bf A}}_3 \left\{\frac{ {\bf w} }{ {\bf x} }\right\}= \left\{ \frac{\bf R}{\bf x} \right\}
%\end{align}
\begin{align}\label{appC_qs}
{\left(\frac{T_i}{T_s} \right)}^2 {\bf A}_1 \frac{d^2 {\bf x}}{dt^2} &+ \left(\frac{T_\nu}{T_s} \right) {\bf A}_2 \frac{d {\bf x}}{dt} + {\bf A}_3 {\bf x} = \frac {\bf F}{k} \\
\frac{d {\bf w}}{dt} &={\bf L}({\bf w}) \nonumber,
\end{align}
where $t = \frac{t^*}{T_s}$ where $t^*$ is dimensional time and $T_s$ is a time-scale of loading defined below, ${T_i}^2=\frac{m}{k}$, $T_\nu=\frac{D}{k}$ (the mass $m$, damping $D$ and stiffness $k$ have physical dimensions of $mass$, ~~$\frac{Force \times time}{length}$ and $\frac{Force}{Length}$ respectively), {\bf x} and {\bf w} are displacements with physical units of $length$, ${\bf L}$ is a function, \emph{independent of $T_s$}, with physical units of $length$ , and ${\bf A}_1$, ${\bf A}_2$ and ${\bf A}_3$ are non-dimensional matrices. In this notation, $\frac{T_i}{T_s}=\epsilon$. In the examples considered, ${\bf L} = \tilde{\bf c} = T_s {\bf c}$, where ${\bf c}$  has physical dimensions of $velocity$ and is assumed given in the form ${\bf c} = \frac{\tilde{\bf c}}{T_s}$  thus serving to define $T_s$; $\tilde {\bf c}$ has dimensions of $length$.

Necessary conditions for the application of the Tikhonov framework are that $\frac{T_i}{T_s} \to 0$, $\frac{T_\nu}{T_s} \to 0$ as $T_s \to \infty$. Those for the quasi-static assumption, commonly used in solid mechanics when loading rates are small, are that $\frac{T_i}{T_s} \to 0$ and $\frac{T_\nu}{T_s} \approx 1$ as $T_s \to \infty$.

In our example, we have $m=1 kg$,$k=2 \times 10^7 N/m$, $D=5 \times 10^3 N s /m$ and $T_s=100 s$. Hence $\frac{T_i}{T_s}=2.24 \times 10^{-6}$ and $\frac{T_\nu}{T_s}=2.5 \times 10^{-6}$, and the damping is not envisaged as variable as $T_s \to \infty$, which shows that the quasi-static approximation is not applicable. 

Nevertheless, due to the common use of the quasi-assumption under slow loading in solid mechanics (which amounts to setting $\epsilon \frac{dy_1}{dt} =0$, $\epsilon \frac{dy_2}{dt}=0$ in \eqref{eq:10.3}), we record the quasi-static solution as well.  

{\bf The Tikhonov framework.} It is easy to see that under the conditions in Case 1, when the walls do not move, all solutions tend to an equilibrium (that may depend on the position of the walls). Indeed, the only way the energy will not be dissipated is when $y_1(t) = y_2(t)$ along time intervals, a not sustainable situation. Thus, we are in the classical Tikhonov framework, and, as we already noted toward the end of the introduction to the paper, the limit solution will be of the form of steady-state equilibrium of the springs, moving on the manifold of equilibria determined by the load, namely the walls.
Computing the equilibria in \eqref{eq:10.3} (equivalently \eqref{eq:10.2}) is straightforward. Indeed, for fixed ($w_1$,$w_2$) we get
\begin{align} \label{appC_tikhonov}
x_1 &= w_1, \quad
x_2 = w_2 \nonumber\\
y_1 &= 0, \qquad
y_2 = 0. 
\end{align}
Under the assumption that $L_1 = 0$ and $L_2 = c_2$ we get
\begin{align} \label{appC_tikhonov2}
w_1(t)  &= 0  \nonumber\\                           
w_2(t)  &= c_2 T_s t.
\end{align}
Plugging the dynamics \eqref{appC_tikhonov2} in \eqref{appC_tikhonov} yields the limit dynamics of  the springs. The real world approximation for $\epsilon$ small would be a fast movement toward the equilibrium set (i.e., a boundary layer which has damped oscillations), then an approximation of  \eqref{appC_tikhonov}-\eqref{appC_tikhonov2}. Our computations in Section \ref{prb2:res_num_case1} corroborate this claim. 

Solutions to  \eqref{eq:10.3} seem to suggest that this is one example where the limit solution \eqref{appC_tikhonov}-\eqref{appC_tikhonov2} is attained by a sequence of solutions of \eqref{eq:10.3} as $\epsilon \to 0$ or $T_s \to \infty$ in a `strong' sense (i.e. not in the `weak' sense of averages); e.g. for small $\epsilon > 0$, $y_2$ takes the value $c_2$ in the numerical calculations and this, when measured in units of slow time-scale $T_s$ (note that $y$ has physical dimensions of velocity), yields $T_s c_2$ which equals the value of the (non-dimensional) time rates of $w_2$ and $x_2$ corresponding to the limit solution \eqref{appC_tikhonov}-\eqref{appC_tikhonov2}; of course, $c_2 \to 0$ as $T_s \to \infty$, by definition, and therefore $y_2 \to 0$ as well. Thus, the kinetic energy and potential energy evaluated from the limit solution \eqref{appC_tikhonov}-\eqref{appC_tikhonov2}, i.e. 0 respectively, are a good approximation of the corresponding values from the actual solution for a specific value of small $\epsilon > 0$, as given in Sec. \ref{prb2:res_num_case1}. 

\par
{\bf The quasi-static assumption.} We solve the system of equations:
\begin{align} \label{eq:10.5}
-k_1 (x_1 - w_1) + {\eta}( {y_2} - {y_1} ) &= 0 \nonumber\\
-k_2 (x_2 - w_2) - {\eta}( {y_2} - {y_1} ) &= 0 \nonumber\\
y_1 = {1 \over T_s} \frac{d x_1}{dt} \nonumber\\
y_2 = {1 \over T_s} \frac{d x_2}{dt}.
\end{align}

We assume the left wall to be fixed and the right wall to be moving at a constant velocity of magnitude $c_2$, so that $w_1=0$ and $w_2=c_2 \, T_s \, t$. This results in
\begin{align}
\frac{d x_1}{dt} + \frac{k_1 T_s}{\eta (1+ \frac{k_1}{k_2} ) } x_1 = \frac {c_2 T_s} {1+ \frac{k_1}{k_2}}. 
\end{align} 
\par 
Solving for $x_1$ and using \eqref{eq:10.5}, we get the following solution
\begin{align}\label{ex2_uneql_quasistatic}
x_1 &= \frac{c_2 \eta}{k_1} + \alpha \, e^{- \beta T_s t }, \quad
y_1 = -\alpha \, \beta e^{-\beta T_s t}, \nonumber\\
x_2 &= c_2 T_s t - \frac{c_2 \eta} { k_2} - \frac{k_1}{k_2} \alpha e^{-\beta T_s t }, \quad
y_2 = c_2 - \frac{k_1}{k_2} \alpha \beta e^{-\beta T_s t}, 
\end{align}
where $\alpha$ is a constant of integration and $\beta=\frac{k_1}{\eta (1 + \frac{k_1}{k_2})}$.

%% For large $t$, these expressions reduce to
%%\begin{align}\label{ex2_uneql_sol}
%%x_1 = \frac {c_2 \eta} {k_1}, \quad
%%y_1 = 0, \quad
%%x_2 = c_2 \, T_s \, t - \frac{c_2 \eta}{k_2}, \quad
%%y_2 = c_2.
%%\end{align}
 
{\bf Remark.} The solution of the unforced system given by (\ref{eq:unforced system}) in Appendix \ref{prb2:equations} is of the form 
\begin{align}\label{eq:sol_complex}
{\bf x} (t) = \sum_{i=1}^4 Q_i e^{\lambda_i t}  {\bf V}_i,  
\end{align}
where ${\lambda_i}$ and ${\bf V}_i$ are the eigenvalues and eigenvectors of {\bf B} respectively. Using the values provided in Table \ref{tab:simulation_details_uneql} to construct {\bf B}, we find that ${\lambda_i}$ and ${\bf V}_i$ are complex. The general real-valued solution to the system can be written as $\sum_{i=1}^4 \psi_i e^{\gamma_i t}  {\bf M}_i(t)$ where 
\begin{table}[H]
\resizebox{1\textwidth}{!}{
\centering
\begin{tabular}[h]{|c|c|c|c|c|}
\hline
$\xi=2.58 \times 10^5$ & $j=1$ & $j=2$ & $j=3$ & $j=4$ \\
\hline
$\gamma_j$     &  $-6.17 \times 10^5$  &   $-1.22 \times 10^5 $ & $-5.63 \times 10^3$ & $-5.63 \times 10^3$ \\
\hline
${\bf M}_{1,j}(t)$ & $0.0001$ & $-0.8732$ & $-0.0001$ & $0.4873$ \\
\hline
${\bf M}_{2,j}(t)$ & $0.0006$ & $-0.7486$ & $-0.0005$ & $0.6630$ \\
\hline
${\bf M}_{3,j}(t)$ & $-0.0001~ cos(\xi t) - 0.0003 ~ sin (\xi t)$ & $-0.6812 ~cos( \xi t )+ 0.1818 ~sin( \xi t)$ & $-0.0003 ~sin( \xi t)$ & $-0.7092 ~cos( \xi t)$ \\
\hline
${\bf M}_{4,j}(t)$ & $0.0003~ cos(\xi t) - 0.0001 ~ sin (\xi t)$ & $-0.1818 ~cos( \xi t )+ 0.6812 ~sin( \xi t)$ & $0.0003 ~cos( \xi t)$ & $-0.7092 ~sin( \xi t)$ \\
\hline
\end{tabular}
}
%\caption{Solution.}
\label{tab:uneql_eigval}
\end{table}
We see that all the real (time-dependent) modes ${\bf M}_i(t)$ are decaying. Moreover, none of the modes describe the dashpot as being undeformed i.e. ${\bf M}_{i,1}(t)={\bf M}_{i,3}(t) $ and ${\bf M}_{i,2}(t)={\bf M}_{i,4}(t)$ (where ${\bf M}_{i,j}(t)$ is the $j^{th}$ row of the mode ${\bf M}_i(t)$). Therefore, solution {\bf x}(t) goes to rest when $t$ becomes large and the initial transient dies.

\section{Example II: Case 2 - Closed-form Solution}\label{prb2:res_case2}

We can convert \eqref{eq:10.3} to the following:
\begin{equation}
	 \begin{pmatrix} \ddot{x_1} \\  \ddot{x_2} \end{pmatrix} 
	+ \alpha {T_s}^2 \begin{pmatrix} x_1 \\ x_2 \end{pmatrix} 
	+ {T_s} \begin{pmatrix} \frac {\eta} {m_1} & -\frac {\eta} {m_1} \\ -\frac {\eta} {m_2} & \frac {\eta} {m_2} \end{pmatrix} \begin{pmatrix} \dot{x_1} \\  \dot{x_2} \end{pmatrix} 	   
	= \begin{pmatrix} \alpha \, w_1 \\  \alpha \, w_2 \end{pmatrix},
\label{sys2}
\end{equation}
 where $\alpha = \frac {k_1} {m_1} = \frac{ k_2} {m_2} $ and $w_1$ and $w_2$ are defined as : 
\begin{align}\label{ex2:load}
w_1(t)=0 \nonumber\\
w_2(t)=c_2 \, T_s \, t .
\end{align}
The overhead dots represent time derivatives w.r.t. $t$.
The above equation is of the form: 
\begin{equation}
\ddot{ {\bf x} } + \alpha  \, {T_s}^2 \, {\bf x} + T_s { \bf A} \dot{ {\bf x} }  =  {\bf g} (t), 
\end{equation}
with the general solution
\begin{subequations} \label{sol}
\begin{equation}
x_1= C_1 cos (\sqrt{\alpha} T_s t) + C_2 sin (\sqrt{\alpha} T_s t) - \frac {m_2} {m_1} C_3 e^{- p_1 T_s t} - \frac {m_2} {m_1} C_4 e^{- p_2 T_s t} + \frac { \eta c_2 } {k_1},
\label{sol1}
\end{equation} 
\begin{equation}
x_2= C_1 cos (\sqrt{\alpha} T_s t) + C_2 sin (\sqrt{\alpha} T_s t) + C_3 e^{- p_1 T_s t} + C_4 e^{- p_2 T_s t} + c_2~T_s~t - \frac { \eta c_2 } {k_2},
\label{sol2}
\end{equation}
\end{subequations}
where $ p_1 = \frac { \eta(m_1 + m_2) + \sqrt{ \eta^2 (m_1 + m_2)^2 - 4 \alpha m_1^2 m_2^2 } } { 2 m_1 m_2} $ and $ p_2 = \frac { \eta(m_1 + m_2) - \sqrt{ \eta^2 (m_1 + m_2)^2 - 4 \alpha m_1^2 m_2^2 } } { 2 m_1 m_2}  $. In the computational results in Section \ref{prb2:res_num_case1} and \ref{prb2:res_num_case2}, we found that $p_1,p_2 > 0$.

Imposing initial conditions  ${x_1}^0$ and ${x_2}^0$ on displacement and ${v_1}^0$ and ${v_2}^0$ on velocity of the two masses $m_1$ and $m_2$ respectively, we obtain
\begin{table}[H]
\centering
\resizebox{0.65\textwidth}{!}{
\begin{tabular}[h]{|c|c|}
\hline
$C_1 $ & $\frac {m_1 {x_1}^0 + m_2 {x_2}^0  } { m_1 + m_2}$ \\[5pt]
\hline 
$C_2$     & $ \frac {m_1 {v_1}^0 + m_2 {v_2}^0 - c_2 m_2 } { \sqrt{\alpha} (m_1 + m_2) }$ \\[5pt]
\hline
$C_3$ & $\frac { c_2 \eta m_1 p_2 + c_2 \eta m_2 p_2 - \alpha c_2 m_1 m_2 - \alpha m_1 m_2 p_2 {x_1}^0 + \alpha m_1 m_2 p_2 {x_2}^0 } { \alpha m_2 (m_1 + m_2) (p_1 - p_2) }$ \\[5pt]
\hline
$C_4$ & $\frac {c_2 \eta m_1 p_1 + c_2 \eta m_2 p_1 - \alpha c_2 m_1 m_2 - \alpha m_1 m_2 p_1 {x_1}^0 +  	 \alpha m_1 m_2 p_1 {x_2}^0 } {\alpha m_2 (m_1 + m_2) (p_1 - p_2)}$\\[5pt]
\hline
\end{tabular}
}
\label{tab:ex2_coeff}
\end{table}

The \textit{closed-form average kinetic energy} ($K^{cf}$), \textit{closed-form average potential energy} ($P^{cf}$)and \textit{closed-form average reaction force} ($R_2^{cf}$) are
\begin{align}\label{eq:analytical_sol}
K^{cf}(t) &= {1 \over \Delta} \int_{t-\Delta}^{t} \left({1 \over 2} m_1 y_1(s)^2 + {1 \over 2} m_2 y_2(s)^2 \right)ds \nonumber \\
P^{cf}(t) &= {1 \over \Delta} \int_{t-\Delta}^{t} \left( {1 \over 2} k_1 x_1(s)^2 + {1 \over 2} k_2 \left(x_2(s)-w_2(s)\right)^2 \right) ds \\
R_2^{cf}(t) &= {1 \over \Delta} \int_{t-\Delta}^{t}  \left( -k_2 \left((x_2(s)-w_2(s)\right) \right) ds \nonumber,
\end{align} 
where $x_1(s)$ and $x_2(s)$ can be substituted from \eqref{sol} and $y_1(s)={1\over T_s} \frac{dx_1}{ds}$ and $y_2(s)={1\over s} \frac{dx_2}{ds}$. 
\par 
{\bf Non-dimensionalization.} Let us denote 
\[
m_{max}=\max \limits_i ~m(x_{\epsilon}(\sigma_i),l_{\epsilon}(\sigma_i)) ,
\]
where $m(x_{\epsilon}(\sigma_i),l_{\epsilon}(\sigma_i))$ is given in \eqref{eq:comp_impl_R1} and \eqref{eq:comp_impl_R2}. Please note that $i$ is chosen such that there is no effect of the initial transient. Using \eqref{ex2_statefn} and computational results in Case 2.4 in Section \ref{prb2:res_num_case2}, we find that
\begin{align}
K_{max}={1 \over 2} C_1^2 (k_1 + k_2), \nonumber ~~
P_{max}={1 \over 2} C_1^2 (k_1 + k_2), \nonumber ~~
R_{2,max}= C_1 k_2. \nonumber 
\end{align}
We introduce the following non-dimensional variables:
\begin{align}\label{non-dim}
\tilde{K}^{cf}=\frac {K^{cf}}{K_{max}},  ~~ 
\tilde{P}^{cf}=\frac {P^{cf}}{P_{max}},
 ~~
\tilde{R_2}^{cf}=\frac {R_2^{cf}}{R_{2,max}}.
\end{align}
Henceforth, while referring to the dimensionless variables, we drop the overhead tilde for simplicity. 
\begin{figure}[!tbh]
\centering
\begin{minipage}{.45\textwidth}
  \centering
  \includegraphics[width=\linewidth]{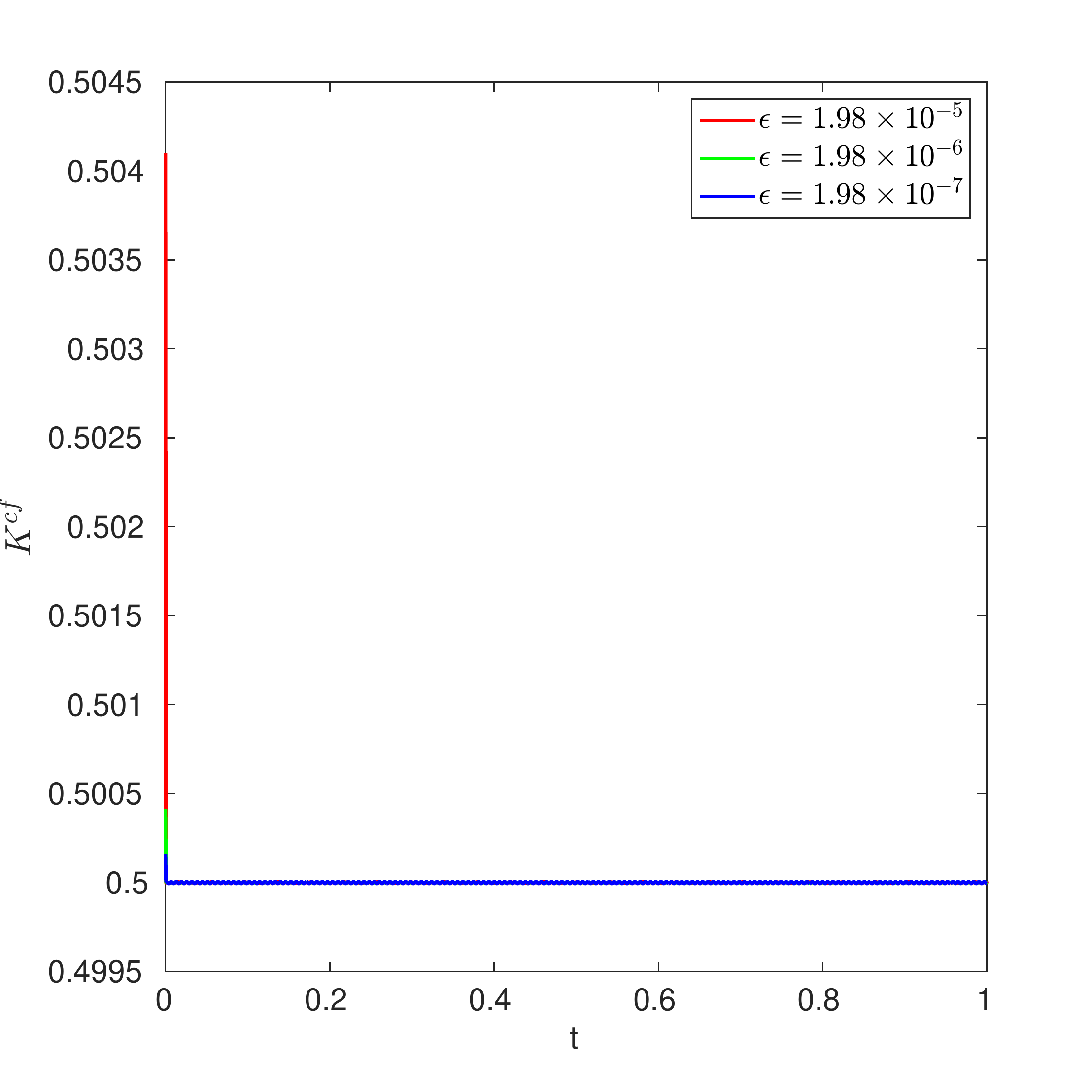}
  \caption{ \textit{$K^{cf}$ as a function of $t$.} }
  \label{fig:ex2_K_analytical}
\end{minipage}%
\hfill
\begin{minipage}{.45\textwidth}
  \centering
  \includegraphics[width=\linewidth]{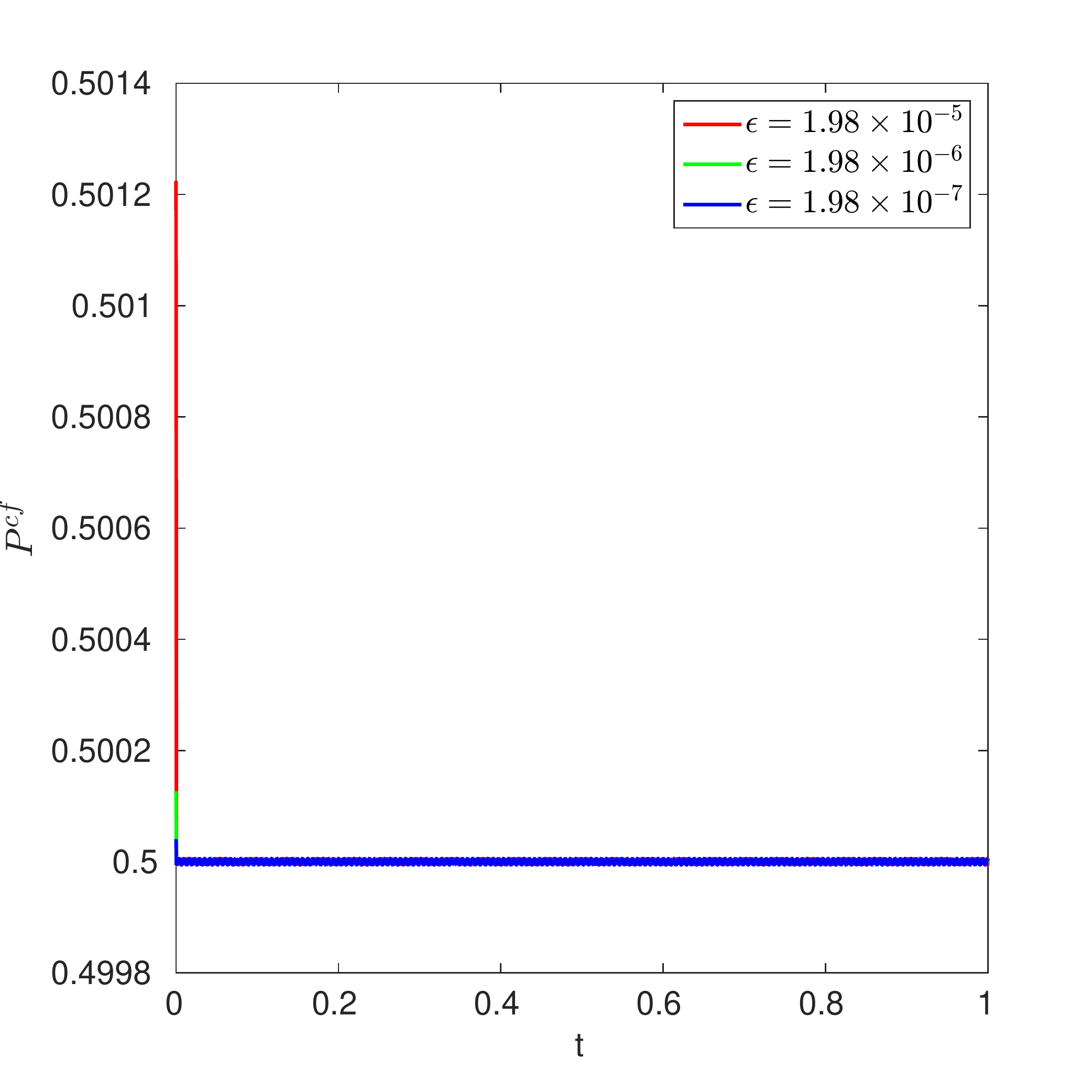}
  \caption{ \textit{$P^{cf}$ as a function of $t$}. }
  \label{fig:ex2_P_analytical}
\end{minipage}
\end{figure}

\begin{figure}[!h]
\centering
\includegraphics[scale = .45]{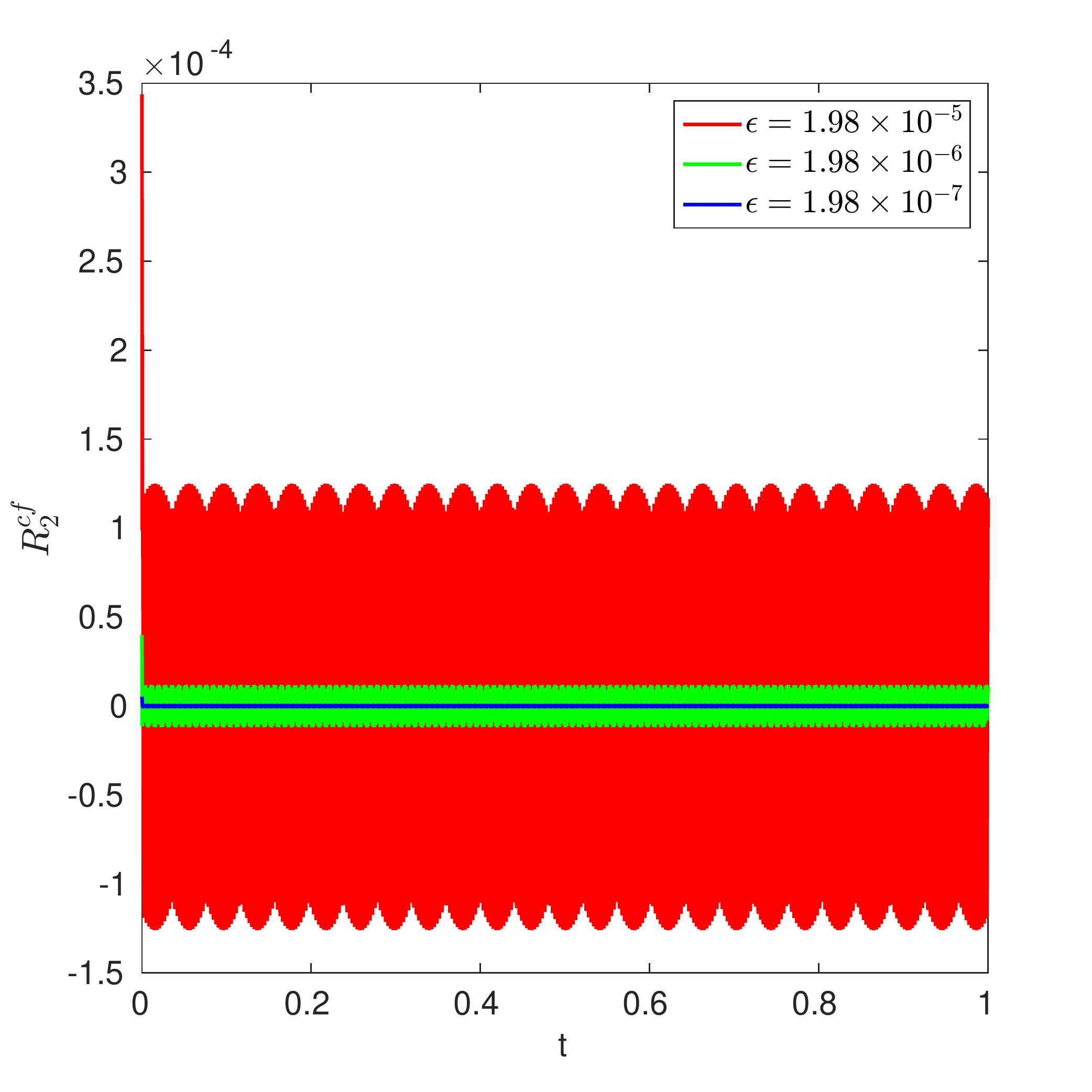}
\caption{ \textit{$R_2^{cf}$ as a function of $t$.} }
  \label{fig:ex2_R_analytical}
\end{figure}

We evaluated \eqref{eq:analytical_sol} numerically (since the analytical expressions become lengthy) with $2 \times 10^6$ integration points over the interval $[t-\Delta,t]$ using Simpson's rule and we used $2000$ discrete points for $t$. We substituted the values provided in Table \ref{tab:simulation_details_uneql}. We repeated the calculations with $4 \times 10^6$ integration points over the interval $[t-\Delta,t]$ and $4000$ discrete points for $t$, and found the results to be the same. We see in Fig. \ref{fig:ex2_K_analytical} and Fig.\ref{fig:ex2_P_analytical} that $K^{cf}$ and $P^{cf}$ oscillate around 0.5 with very small amplitude for different values of $\epsilon$ (recall that $\epsilon = \frac{T_f}{T_s}$ and we think of $T_f$ being fixed with $T_s \rightarrow \infty$ to effect $\epsilon \rightarrow 0$). However, as we decrease $\epsilon$, the amplitude of oscillations of $R_2^{cf}$ decreases and it goes to zero for $\epsilon=1.98 \times 10^{-7}$ as we see in Fig. \ref{fig:ex2_R_analytical}. We use these `closed-form' results to compare with the PTA results in Section \ref{prb2:res_num_case1} and Section \ref{prb2:res_num_case2}.  

\par
{\bf Remark.} 
The criteria for convergence of $R^m_t$ (as mentioned in the discussion following \eqref{eq:alg_fast_2} in Section \ref{impl_algo}) is discussed as follows. Let us denote 
\[
m_{I,t}={1 \over I} \sum_{i=1}^{I} m(x_\epsilon(\sigma_i),l_\epsilon(\sigma_i)),
\]
where $x_\epsilon(\sigma_i),l_\epsilon(\sigma_i)$ is defined in \eqref{eq:alg_fast_2}, $I \in \mathbb{Z_+}$ and $m(x_\epsilon(\sigma_i),l_\epsilon(\sigma_i))$ is non-dimensionalized $\forall i \in [1,I]$. Then, we say that $m_{I,t}$ has converged if 
\begin{align}\label{eq:conv_crit}
|\frac{m_{I-jk,t} ~~ - ~~ m_{I-pk,t}}{m_{I-pk,t}}| \leq tol_1^m,  
\end{align}
$\forall j \in [0,p]$ where $p,k \in \mathbb{Z_+}$ and $p,k<I$ and $tol_1^m$ is a specified value of tolerance (which is a small value generally around $10^{-2}$). We declare $N_t=I$ (where $N_t$ is defined in the discussion following \eqref{eq:alg_fast_2}) and $R^m_t=m_{I,t}$. In situations where $m_{I-pk,t}$ becomes very small so that the convergence criteria in \eqref{eq:conv_crit} cannot be practically implemented, we say that $m_{I,t}$ has converged if
\begin{align}
|m_{I-jk,t}| \leq tol_2^m,
\end{align}
$\forall j \in [0,p]$ and $tol_2^m$ is a specified value of tolerance (generally around $10^{-5}$).
\par
{\bf Remark.} We used the Simpson's rule of numerical integration to obtain the value of the slow variable instead of using \eqref{eq:v(t+h)} to obtain the results in Section \ref{prb2:res_num_case1} and \ref{prb2:res_num_case2}. The Simpson's rule of numerical integration for any function $f(t)$ over the interval $[a,b]$ where $a,b \in \mathbb{R}$ is \[
\int_a^b f(t)dt \approx \frac{\Delta t}{n} [f(t_0)+4f(t_1)+2f(t_2)+...+2f(t_{n-2})+4f(t_{n-1})+f(t_n)],
\]
where $\Delta t=\frac{b-a}{n}$, $n \in \mathbb{Z_+}$, $n$ is even and $t_i=t_0+i\Delta t$ for $i \in \mathbb{Z}$ and $i \leq n$. All the function evaluations at points with odd subscripts are multiplied by 4 and all the function evaluations at points with even subscripts (except for the first and last) are multiplied by 2. Since we are calculating the value of the slow variable $v$ given by \eqref{coarse_obs_impl}, the function $f(t)$ is given by $R^m_t$ which is defined in \eqref{eq:comp_impl_R1}, with $\epsilon = 0$ in \eqref{eq:alg_fast_2}. We chose $n=2$ and calculated the values of $R^m_{t-\Delta+\frac{i}{n}\Delta}$ using \eqref{eq:comp_impl_R1} for $i = 0,1,2$. We obtained the fine initial conditions 
$x_{guess}(t-\Delta+\frac{i}{n}\Delta)$ as  
\[
x_{guess}(t-\Delta+\frac{i}{n}\Delta)=x^{arb}_{t-\Delta} + \frac {\left( x^{arb}_{t-\Delta} - x^{cp}_{t-h} \right)}{\left(h-\Delta \right)} \, \frac{i\Delta}{n},
\]
where $x^{arb}_{.}$ and $x^{cp}_{.}$ are defined in Step 3 in Section \ref{impl_algo}.

Please note that we also used Simpson's rule to evaluate \eqref{eq:analytical_sol} with $n=2\times10^6$ and the details of the calculation are mentioned in \textit{Non-dimensionalization} in this Appendix above. 

\par
{\bf Remark.} To see the effect of the initial condition on the solution, we solve (\ref{eq:unforced system}) in a particular case using the values provided in Table \ref{tab:simulation_details_uneql} but with $k_2=2 \times 10^7 N/m$ and $c_2=0$. The general real-valued solution to the system can be written as $\sum_{i=1}^4 \psi_i e^{\gamma_i t}  {\bf M}_i(t)$ where
\begin{table}[H]
\resizebox{1\textwidth}{!}{
\centering
\begin{tabular}[h]{|c|c|c|c|c|}
\hline
$\omega=3.16 \times 10^5 $ & $j=1$ & $j=2$ & $j=3$ & $j=4$ \\
\hline
$\gamma_j$     &  $-5.76 \times 10^5$  &   $-1.73 \times 10^5 $ & $0$ & $0$ \\
\hline
${\bf M}_{1,j}(t)$ & $0.0002$ & $-0.8944$ & $-0.0001$ & $0.4472$ \\
\hline
${\bf M}_{2,j}(t)$ & $0.0005$ & $-0.8944$ & $-0.0003$ & $0.4421$ \\
\hline
${\bf M}_{3,j}(t)$ & $\quad \quad 0.0002~ sin(\omega t)$ & $\quad \quad -0.7071 ~cos( \omega t )$ & $\quad \quad0.0002 ~sin( \omega t)$ & $\quad \quad0.7071 ~cos( \omega t)$ \\
\hline
${\bf M}_{4,j}(t)$ & $\quad \quad-0.0002~ cos(\omega t)$ & $\quad \quad0.7071 ~sin( \omega t )$ & $\quad \quad-0.0002 ~cos( \omega t)$ & $\quad \quad0.7071 ~sin( \omega t)$ \\
\hline
\end{tabular}
}
\label{tab:uneql_eigval}
\end{table}
\vspace{-0.5cm}
We see that while ${\bf M}_1(t)$ and ${\bf M}_2(t)$ are decaying real (time-dependent) modes, ${\bf M}_3(t)$ and ${\bf M}_4(t)$ are the non-decaying real (time-dependent) modes. Moreover, both ${\bf M}_3(t)$ and ${\bf M}_4(t)$ describe the dashpot as being undeformed i.e. ${\bf M}_{i,1}(t)={\bf M}_{i,3}(t) $ and ${\bf M}_{i,2}(t)={\bf M}_{i,4}(t)$ (where ${\bf M}_{i,j}(t)$ is the $j^{th}$ row of the mode ${\bf M}_i(t)$). The solution ${\bf x}(t)$ described in (\ref{eq:sol_complex}) of Appendix \ref{prb2:res_case1} can be written as 
\begin{align}\label{eq:coeffs}
{\bf x}(t) = \sum_{i=1}^4 \kappa_i {\bf M}_i(t) ,
\end{align}
where the coefficients $\kappa_i$ are obtained from the initial condition ${\bf x}_0 $ using 
\[
\kappa_i = {\bf x}_0 \cdot {\bf M}_i^d (0) 
\]
where ${\bf M}_i^d (t)$ are the dual basis of ${\bf M}_i (t)$. 
\end{appendices}

\section*{Acknowledgment}
We thank the anonymous referees for their valuable comments that helped improve the presentation of our paper. S. Chatterjee acknowledges support from NSF grant NSF-CMMI-1435624. A. Acharya acknowledges the support of the Rosi and Max Varon Visiting Professorship at The Weizmann Institute of Science, Rehovot, Israel, and the kind hospitality of the Dept. of Mathematics there during a sabbatical stay in Dec-Jan 2015-16. He also acknowledges the support of the Center for Nonlinear Analysis at Carnegie Mellon and grants NSF-CMMI-1435624, NSF-DMS-1434734, and ARO W911NF-15-1-0239. Example I is a modification of an example that appeared in a draft of [AKST07], but did not find its way to the paper.
\bibliography{testbibliog}
\bibliographystyle{alpha}
\end{document}